%% file: ise.tex
\documentclass[10pt,a4paper,reqno]{amsart} 
\usepackage{latexsym,floatflt}
\usepackage{epsfig}
\usepackage{rotating}
\usepackage{amssymb}
\usepackage[T1]{fontenc}% pour les accents français
\usepackage{afterpage}
\usepackage{color}

\usepackage{colordvi}

\catcode`\@=11
\def\section{\@startsection{section}{1}%
 \z@{.7\linespacing\@plus\linespacing}{.5\linespacing}%
 {\normalfont\bfseries\scshape\centering}}
\def\subsection{\@startsection{subsection}{2}%
  \z@{.5\linespacing\@plus\linespacing}{.5\linespacing}%
  {\normalfont\bfseries\scshape}}

\def\subsubsection{\@startsection{subsubsection}{3}%
 \z@{.5\linespacing\@plus\linespacing}{-.5em}%{.5\linespacing}%
  {\normalfont\bfseries\itshape}}
\catcode`\@=12

\def\cqfd{$\hfill{\vrule height 3pt width 5pt depth 2pt}$}

\newcommand{\GK}{\mathbb{K}}

\newcommand{\ns}{\mathbb{N}}
\newcommand{\zs}{\mathbb{Z}}

\newcommand{\qs}{\mathbb{Q}}
\newcommand{\rs}{\mathbb{R}}

\newcommand{\cs}{\mathbb{C}}

\newcommand{\A}{\mathcal A}
\newcommand{\I}{\mathcal I}
\newcommand{\D}{\mathcal D}
\newcommand{\E}{\mathcal E}
\newcommand{\C}{\mathcal C}
\newcommand{\T}{\mathcal T}
\newcommand{\J}{\mathcal J}
\renewcommand{\S}{\mathcal S}
\newcommand{\Sn}{\ensuremath{\mathcal S}} 

\newcommand{\Ha}{\mathcal H}

\newcommand{\EE}{\mathbb{E}}
\newcommand{\PP}{\mathbb{P}}

\newcommand{\Ref}[1]{(\ref{#1})}
\newcommand{\beq}{\begin{equation}}
\newcommand{\eeq}{\end{equation}}
\newcommand{\gf}{generating function}
\newcommand{\gfs}{generating functions}

\newcommand{\fps}{formal power series}
\def\cqfd{\par\nopagebreak\rightline{\vrule height 4pt width 5pt depth 1pt}
\medbreak}
\def\emm#1,{{\em #1}}
\newcommand{\deq}{\stackrel{\mathrm{d}}{=}}
\newcommand{\dconv}{\stackrel{\mathrm{d}}{\rightarrow}}
\font\dsrom=dsrom10 scaled 1400
\def \1{\textrm{\dsrom{1}}}

\newtheorem{Theorem}{Theorem}
\newtheorem{Propo}[Theorem]{Proposition}

\newtheorem{Lemma}[Theorem]{Lemma}

\newtheorem{Conj}[Theorem]{Conjecture}

\begin{document}

%%%%%%%%%%%%%%%%%%%%%%%%%%%%%%%%%%%%%%%%%%%%%%%%%%%%%%%%%%%%
\title[Limit laws for embedded trees]{Limit laws for embedded trees.\\
 Applications to the integrated
  superBrownian excursion}

\author{Mireille Bousquet-M\'elou}
\address{CNRS, LaBRI, Universit\'e Bordeaux 1, 351 cours de la Lib\'eration,
  33405 Talence Cedex, France}
\email{mireille.bousquet@labri.fr}
\thanks{MBM  was partially supported by the European Commission's IHRP
  Programme, grant HPRN-CT-2001-00272, ``Algebraic Combinatorics in
  Europe''} 
%\subjclass{Primary: subject; Secondary: subject}
\date{January 17, 2005}

\begin{abstract}
We study three families of labelled  plane
trees. In all these trees, the root is labelled $0$, and the labels of
two adjacent nodes differ by $0, 1$ or $-1$. 

One part of the paper is devoted to enumerative results. For each
family, and for all $j\in \ns$, we obtain closed form expressions for the
following three generating functions: the generating function of trees having no label larger than
$j$; the (bivariate) generating function  of trees, counted by the  number of
edges and the number of nodes labelled $j$; and finally the
(bivariate) generating function  of trees, counted by the  number of edges
 and the number of nodes labelled \emm at least, $j$. Strangely
enough, all these series turn out to be algebraic, but we have no
combinatorial intuition for this algebraicity.

The other part of the paper is devoted to deriving limit laws from
 these enumerative results. In each of our families of
trees, we endow the trees of size $n$ with the uniform
distribution, and study the following random variables: $M_n$, the
largest label occurring in a (random) tree; $X_n(j)$, the number of
nodes labelled $j$; and $X_n^+(j)$, the number of nodes labelled $j$
or more. We obtain limit laws for scaled versions of these random
variables.

Finally, we translate the above limit results into statements dealing
with  the
integrated superBrownian excursion (ISE). In particular, we describe
the law of the supremum of its support (thus recovering some earlier
results obtained by Delmas), and the law of its
distribution function at a given point. We also conjecture the law of
its density (at a given point).
\end{abstract}

\maketitle
%%%%%%%%%%%%%%%%%%%%%%%%%%%%%%%%%%%%%%%%%%%%%%%%%
\section{Introduction}
\label{section-intro}
We study in this paper three families of labelled  plane
trees.  In all these trees, the root is labelled $0$, and the labels of
two adjacent nodes differ by $0, 1$ or $-1$. 

\begin{figure}[hbt]
\begin{center}
\input{tree-examples.pstex_t}
\end{center}
\caption{A labelled plane tree with increments $\pm1$. --- A labelled tree
  with increments $0, \pm 1$. --- A naturally embedded binary tree.}
\label{fig-three-families}
%\hrule
\end{figure}
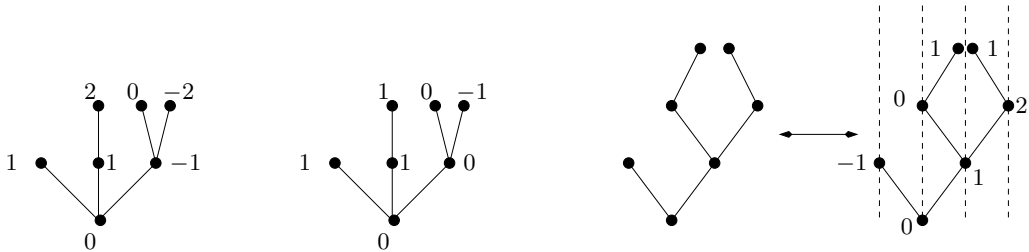

 More precisely, the first family we consider is the set of plane
 trees, and the increments of 
the labels along edges are constrained to be $\pm 1$. In the closely
related second family, these increments can be $0, \pm 1$. The third
family is a bit different. It is simply the set of (incomplete) binary
trees, in 
which the nodes are labelled in a deterministic way: the label
of a node is the difference between the number of right steps and the
number of left steps occurring in the path that yields from the root
 to the node under consideration. See
Figure~\ref{fig-three-families} for an illustration. We call this
labelling the \emm natural labelling, of the binary
tree. Note that the label of each node is simply its abscissa, if we
draw the tree in the plane in such a way the right (resp.~left) son of
a node lies one unit to the right (resp.~left) of its father. For this
reason, we will sometimes call these labelled binary trees \emm
naturally embedded binary trees., More generally, for \emm any,\/
plane labelled tree,  we may consider that the label of each node tells
where to embed it  in  $\zs$; hence the title of the paper.

In each of these three families, we endow the set of trees having a
given size (say, $n$ edges) with the uniform distribution.  We address
(via \gfs) 
the following three questions:
\begin{enumerate}
\item
What is the maximal label that occurs in the tree? This label is in
fact a random variable $M_n$. We prove that $M_n/n^{1/4}$ converges in
distribution to a random variable $N$ having a density. We give this
density explicitly. We also compute the moments of $N$ and prove the
convergence of the moments of $M_n/n^{1/4}$  
 to those  of $N$.
\item
How many nodes of the tree have label $j$? Let $X_n(j)$ denote the
corresponding random variable. If $j$ is fixed, and $n$ goes to
infinity, then the answer to this question is independent of $j$. We
prove that for any $j \in \zs$, the variable  $X_n(j)/n^{3/4}$
converges in distribution to  $cT^{-1/2}$, where 
$c$ is a constant depending on which family of trees we consider, and
$T$ follows a unilateral stable law of parameter $2/3$. 

Given that the maximal label grows like $n^{1/4}$, we get a better
insight on the label distribution by asking how many nodes in a tree
of size $n$ have label $\lfloor  \lambda n^{1/4} \rfloor$. We prove
that, for any $\lambda \in \rs$, the random variable $X_n(\lfloor
\lambda n^{1/4} \rfloor)/n^{3/4}$ converges in distribution to a limit
variable $Y(\lambda)$. This variable  admits a Laplace transform,
which we give explicitly. The convergence of the Laplace transform, and
of the moments,  hold as well. We say we have obtained a \emm
local limit law, for embedded trees, because we look at \emm one, value of
the labels only.
\item
Finally, we also obtain a \emm global limit law, by studying the
variable  $X_n^+(j)$ that gives the number of nodes having label $j$
at least. Remarkably, we prove that $X_n^+(0)/n$, the (normalized)
number of nodes having a non-negative label, converges to the uniform
distribution on $[0,1]$. More generally, for $\lambda \in \rs$, the
variable $X_n^+(\lambda n^{1/4})/n$ converges in distribution to a
variable $Y^+(\lambda)$.  This variable  admits a Laplace transform,
which we give explicitly. Once again, the convergence of the Laplace
transform, and of the moments,  hold as well. 
\end{enumerate}
The laws of $N$, $Y(\lambda)$ and $Y^+(\lambda)$ naturally depend on
which family of trees we consider, but only by a simple normalization
factor. 

\subsection{Embedded trees and the integrated superBrownian excursion}
Why should one study such labelled trees? 

The first two classes of
trees we consider have a close connection with certain families of
planar
maps~\cite{bdg-statistics,chassaing-schaeffer,cori-vauquelin}. In
particular, the 
diameter of a random quadrangulation having $n$ faces is distributed
like the largest 
label in \emm non-negative, random trees of our second
family. Moreover, once scaled by $n^{1/4}$, this diameter has the same limit law
as $(M_n-m_n)n^{-1/4}$, where $M_n$ (resp.~$m_n$) is the largest
(resp.~smallest)  label occurring in
a random tree of our second family~\cite{chassaing-schaeffer}.

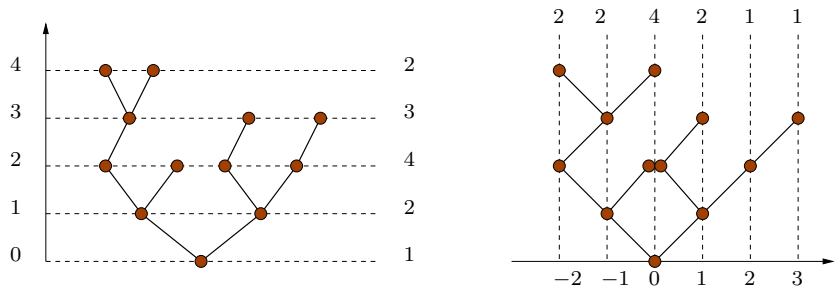
\begin{figure}[t]
\begin{center}
\input{profils.pstex_t}
\end{center}
\caption{An (incomplete) binary tree having horizontal profile
  $[1,2,4,3,2]$ and   vertical profile $[2,2,4; 2,1,1]$.}
\label{fig-profils}
%\hrule
\end{figure}

The third class we study
is the good old family of binary trees, and this may suffice to motivate
its study! More seriously, the three questions addressed above
have, for binary trees, a natural geometric formulation. 
The random variable $M_n$ (the 
maximum label) tells us about the ``true width'' of a binary tree (as
opposed to the maximal number of nodes lying at the same
level, which is known to grow like $\sqrt n$). More generally, the
variables $X_n(j)$ 
tell us about the \emm vertical 
profile, of the tree (as opposed to the \emm horizontal profile, which
describes the repartition of nodes by
level~\cite{drmota-profile}). See Figure~\ref{fig-profils}. 

We may also invoke an \emm a posteriori,\/ justification to the study
of these trees: the form of the \gfs\ we obtain is  remarkable,
whatever family of trees we consider, and suggests that there must be
some beautiful hidden combinatorics in these problems, which should be
explored further.

\begin{figure}[b]
\begin{center}
\includegraphics[height=3cm,width=3cm]{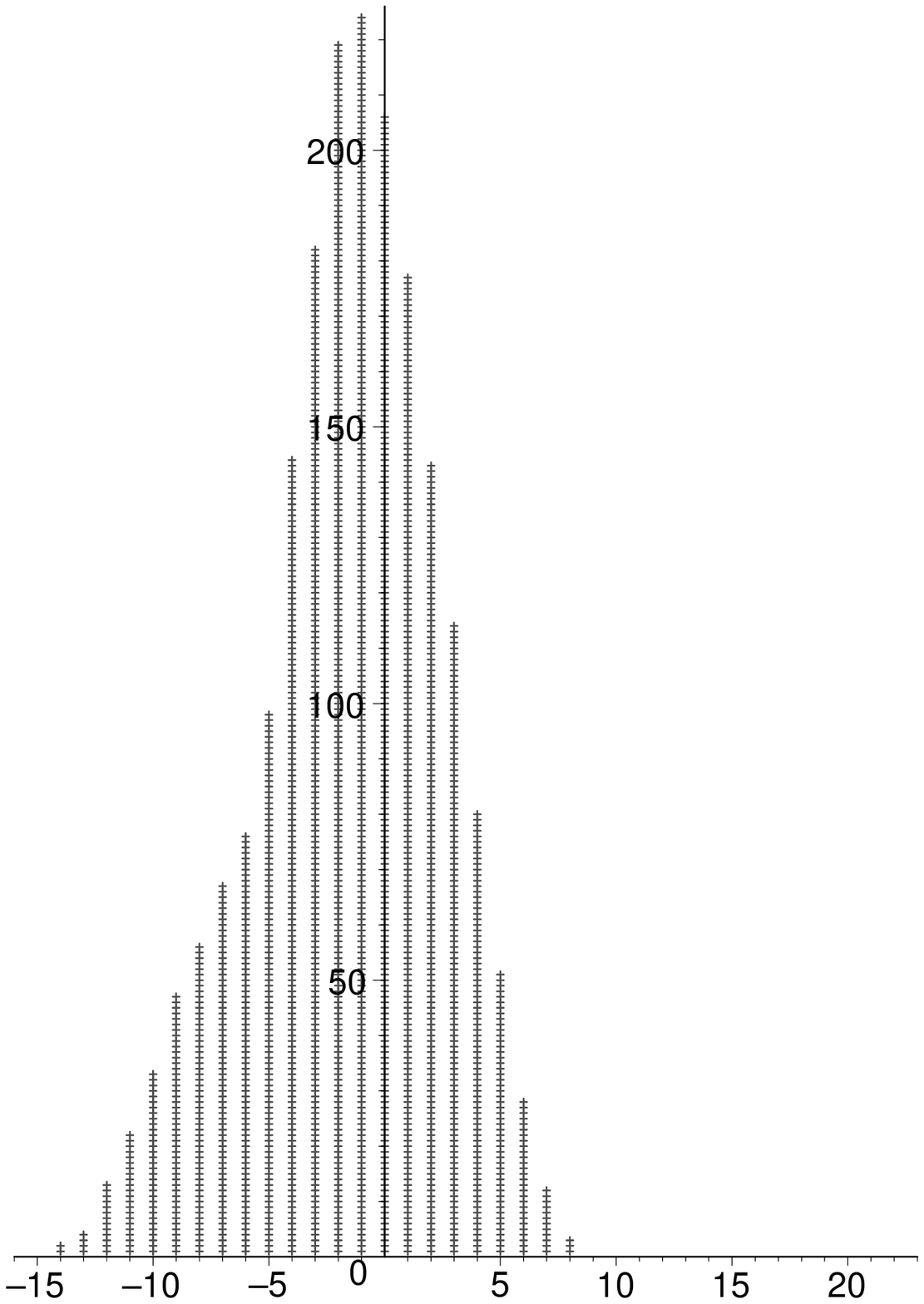}
\hskip 2mm
\includegraphics[height=3cm,width=3cm]{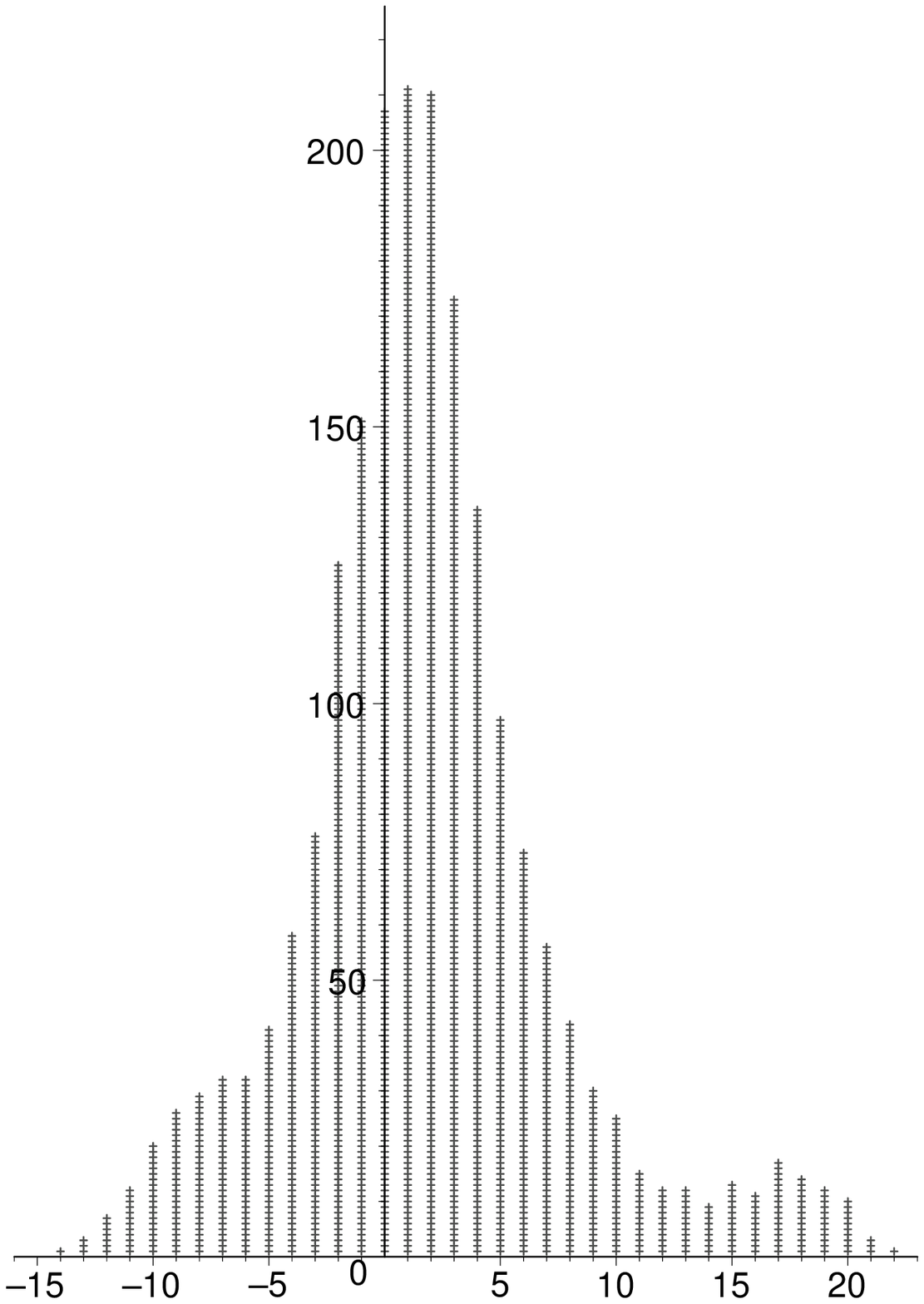}
\hskip 2mm
\includegraphics[height=3cm,width=3cm]{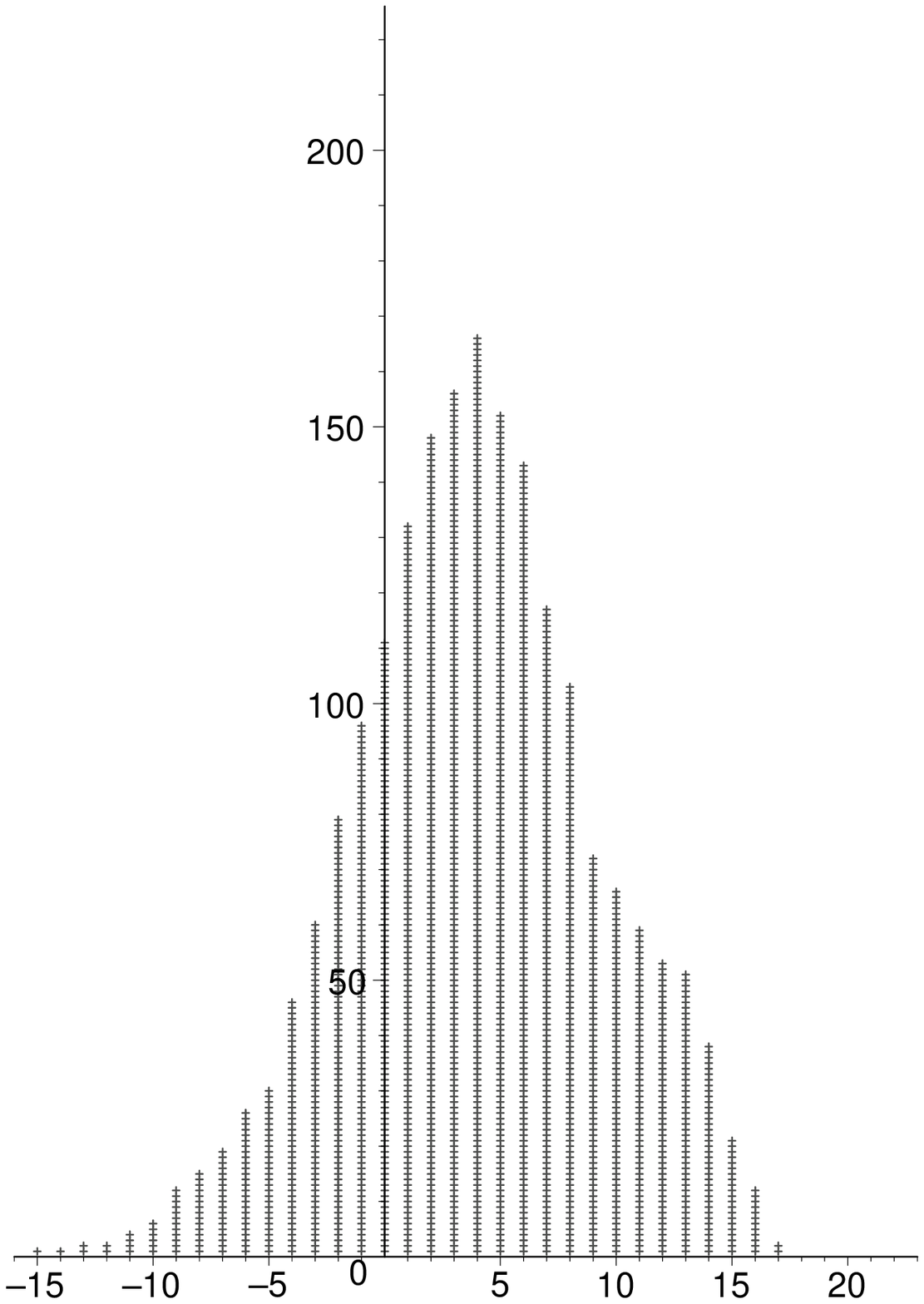}
\hskip 2mm
\includegraphics[height=3cm,width=3cm]{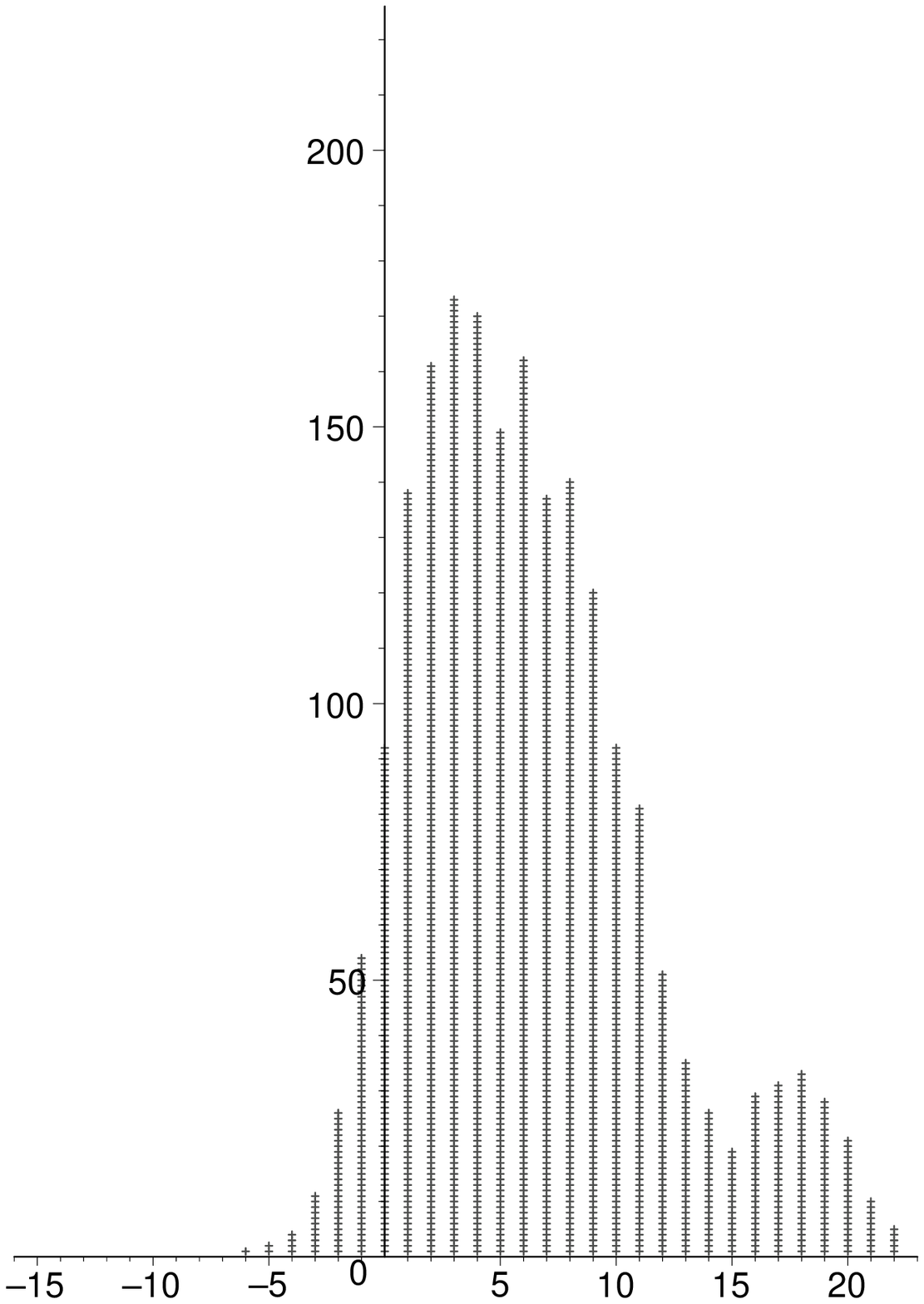}

\includegraphics[height=3cm,width=3cm]{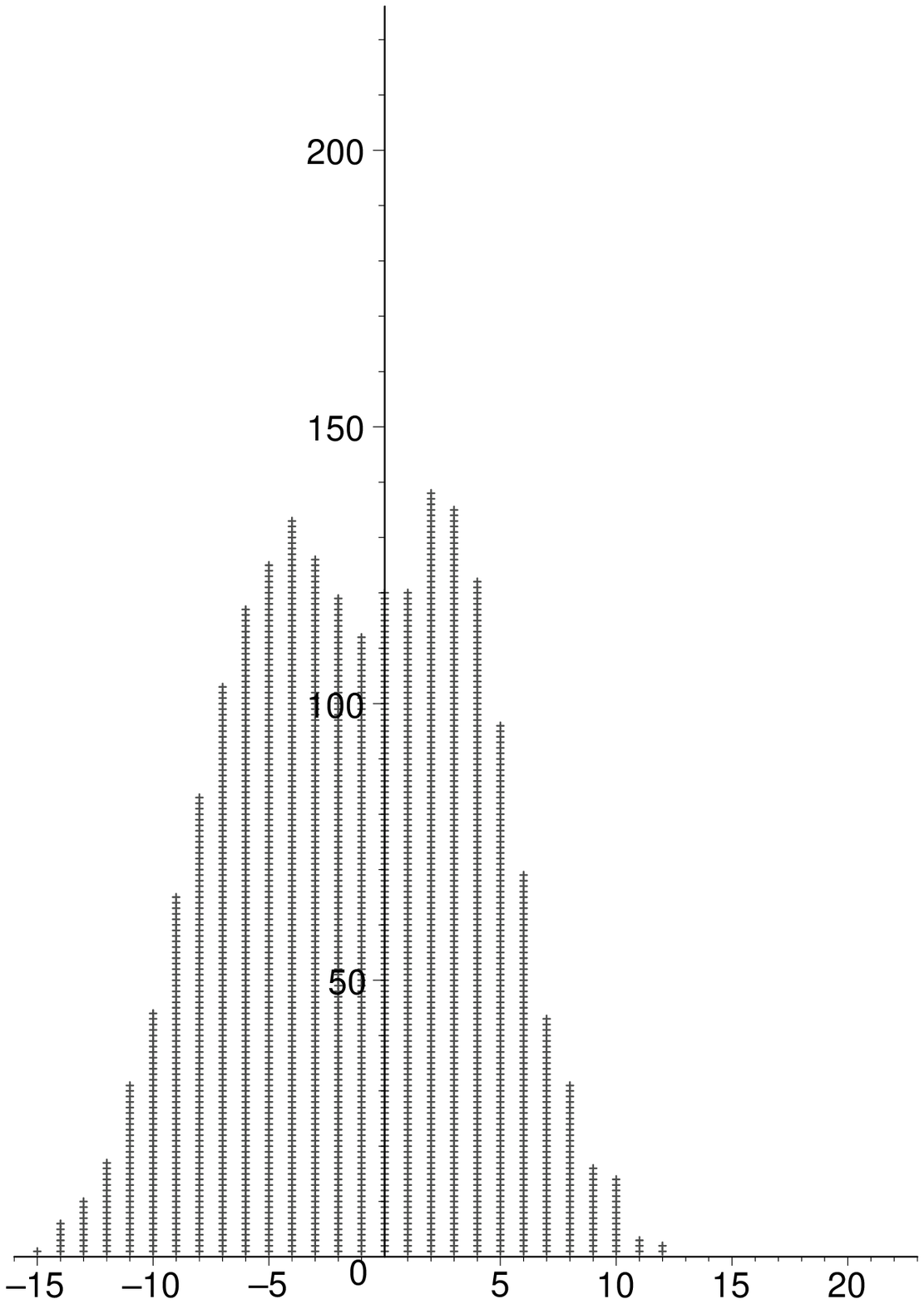}
\hskip 2mm
\includegraphics[height=3cm,width=3cm]{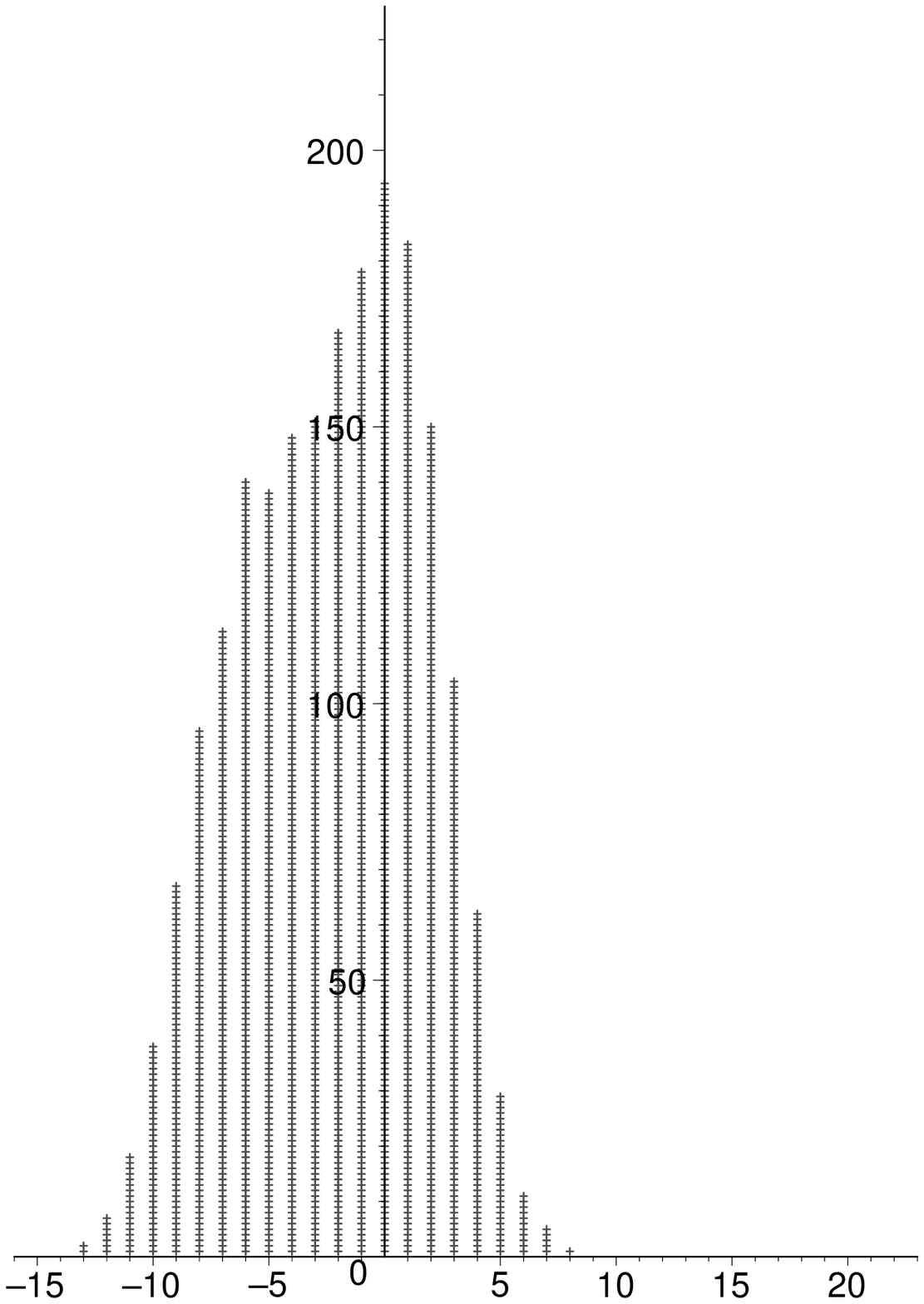}
\hskip 2mm
\includegraphics[height=3cm,width=3cm]{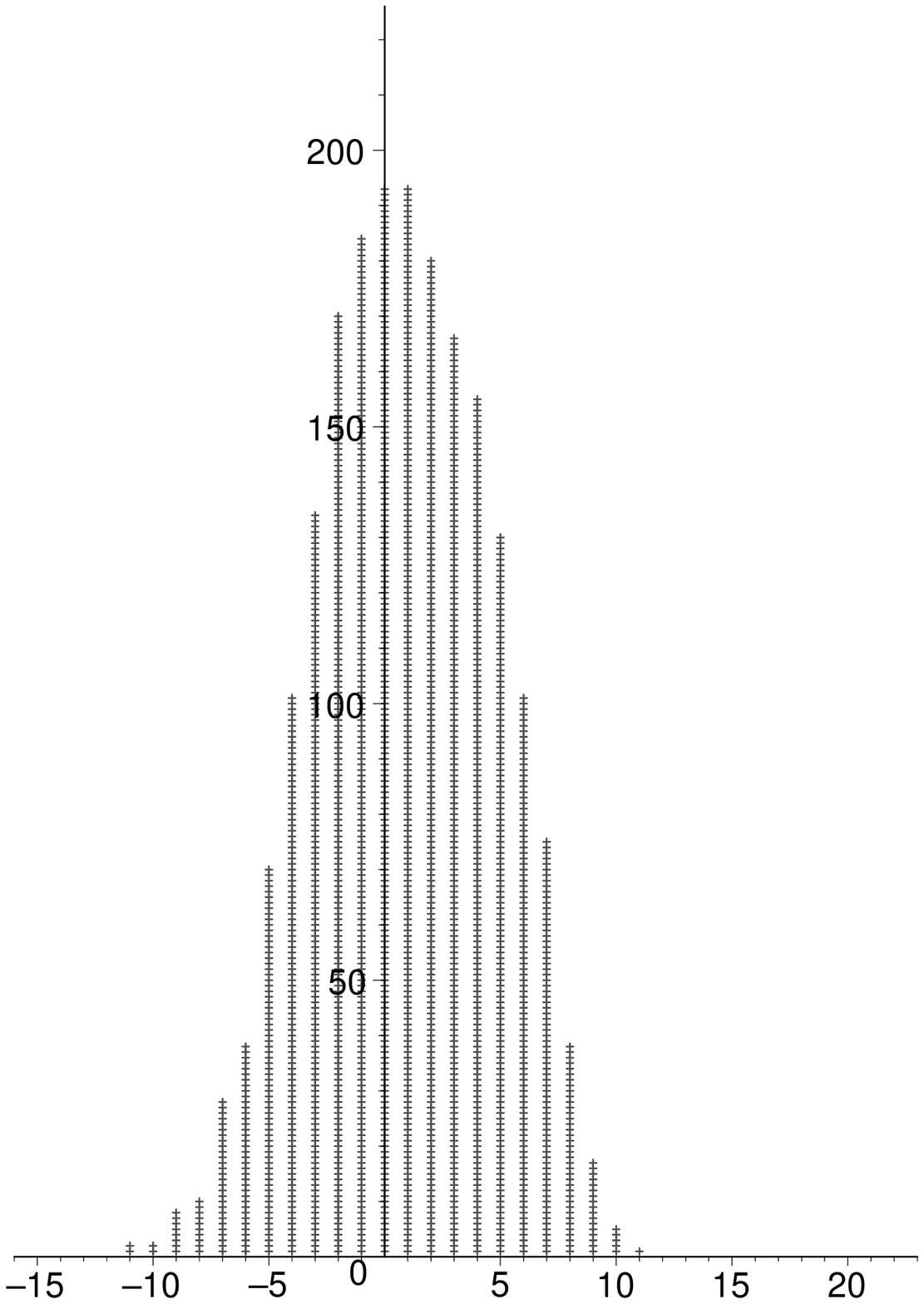}
\hskip 2mm
\includegraphics[height=3cm,width=3cm]{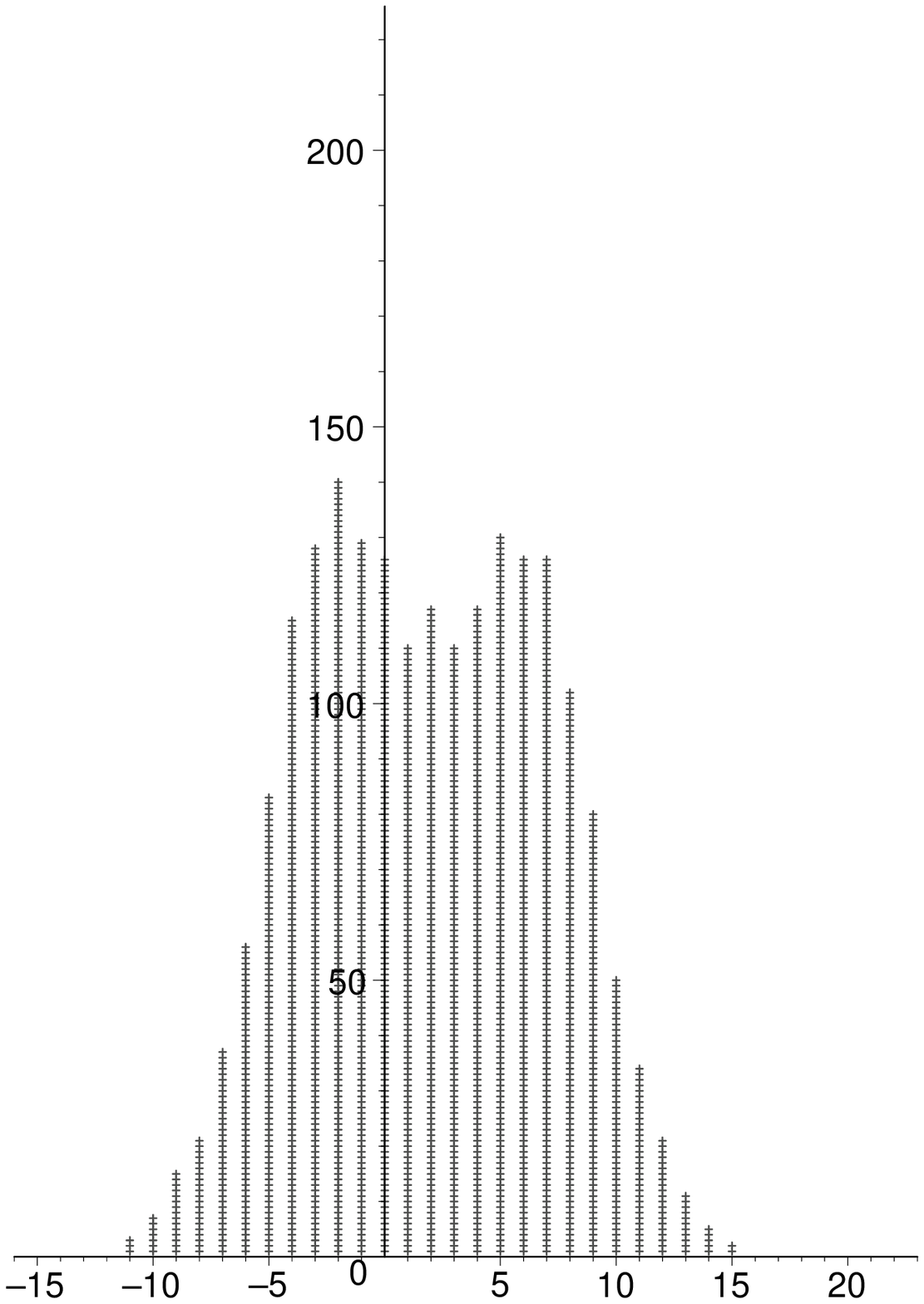}
\end{center}
\caption{The plot of $X_n(j)$ vs. $j$ for random binary trees with
  $n=1000$ nodes.}
\label{fig-simulations}
\end{figure}

However, the main motivation for this work is the connection between
embedded trees and the \emm integrated superBrownian excursion,
(ISE). Choose one of the three families of trees, and consider the
following \emm random, probability distribution on $\rs$:
\beq
\label{random-measure}
\mu_n = \frac 1 {n+1} \sum_{j\in \zs} X_n(j) \delta_{cjn^{-1/4}},
\eeq
where $X_n(j)$ is the (random) number of nodes labelled $j$,
$\delta_x$ denotes the Dirac measure at $x$, and the constant $c$ 
equals $\sqrt 2$ for the first family, $\sqrt 3$ for the second one and
$1$ for the family of binary trees. Then $\mu_n$ is
known to converge weakly
to a limiting random probability distribution called the
ISE~\cite{aldous,marckert-mokka-snake,jf-homothetie,jf-janson}. See
Figure~\ref{fig-simulations} for simulations of $\mu_n$.

Our limit results provide some information  about the law of the
ISE. For instance, we prove  that $cM_nn^{-1/4}$,  the largest point
having a positive weight under $\mu_n$, converges in law to $N_{\hbox{\sc ise}}$, the
supremum of the support of the ISE. We denote this by
$$
cM_nn^{-1/4}\dconv N_{\hbox{\sc ise}}.
$$
 The results we obtain for the limit law of $M_nn^{-1/4}$ thus
translate into expressions of the moments, distribution function and
density of the supremum of the ISE. Note that the moments were already
obtained by Delmas~\cite{delmas}.
Our second limit result deals with the random variables
$X_n(\lfloor  \lambda n^{1/4} \rfloor)$. Observe that
\beq
\label{limit-density}
\mu_n(c\lambda-cn^{-1/4}, c\lambda]= \frac 1 {n+1} X_n(\lfloor  \lambda n^{1/4} \rfloor) .
\eeq
 This leads us to \emm conjecture,\/
that the random variable $Y(\lambda)$ involved in our local
limit law 
%(it describes the limit law of  $n^{-3/4} X_n(j)$) 
satisfies
\beq
\label{conv-densite}
Y(\lambda) \deq cf_{\hbox{\sc \small ise}}(c\lambda)
 \eeq
where $f_{\hbox{\sc \small ise}}$ is the (random) density of
the ISE. Similarly, 
$$
\mu_n[c\lambda, +\infty)= \frac 1 {n+1} X_n^+(\lceil  \lambda n^{1/4} \rceil), 
$$
 and we \emm prove,\/ that  the
random variable $Y^+(\lambda)$ involved in our global limit law
satisfies
$$
Y^+(\lambda) \deq g_{\hbox{\sc \small ise}}(c\lambda)
$$
where $g_{\hbox{\sc \small ise}}$ is the
(random) tail distribution function of the ISE. The results we obtain
about the laws of $Y(\lambda)$ and $Y^+(\lambda)$ thus translate into
formulas for the Laplace transforms of $f_{\hbox{\sc \small
    ise}}(\lambda)$ and  $g_{\hbox{\sc \small ise}}(\lambda)$ (the
formula for $f_{\hbox{\sc \small
    ise}}(\lambda)$ being conjectural).

Our conjecture on $f_{\hbox{\sc \small
    ise}}$  is naturally supported
  by the fact that the law of $Y(\lambda/c)/c$ is 
  independent of the tree family we start from. This is one of the
  reasons why we consider as many as three families of trees. The
  other reasons 
  involve the connections with planar maps, the remarkable form of the
  \gfs\ we obtain, and our unshakeable
  interest in binary trees.  The details of the calculations
  are only given for the first of the three families
  (Sections~\ref{section-snake-gf} 
  to~\ref{section-limit-global}), while  the results are merely stated
  for the other   two families   (Section~\ref{section-universality}).

Let us finally mention that the moments of the \emm center of mass,\/
of the ISE have recently been determined by two different
approaches~\cite{chassaing-janson,janson}. In our discrete setting,
this boils down to studying the convergence of the variable
$$
\frac 1{n^{5/4} }\sum _{j\in \zs} jX_n(j).
$$

\subsection{Overview of the paper}
The starting point of our approach is a series of exact enumerative
results dealing with our first class of trees: plane trees in which
the labels of  adjacent nodes differ by $\pm 1$. These
results are gathered in the next section. We obtain for 
instance an explicit expression for the bivariate \gf \ of labelled trees,
counted by the number of edges and the number of nodes labelled $j$
(for $j$ fixed). This section includes, and owes a lot to, some results
recently obtained by Bouttier, Di Francesco and
Guitter~\cite{bdg-geodesic,bdg-statistics} on the 
enumeration of  trees having no label greater than $j$. This
part of our work raises a number of challenging combinatorial
questions --- why are these expressions so simple? --- which are not
addressed in this paper. 

The limit behaviours of the random variables $M_n$, $X_n(\lfloor
\lambda n^{1/4}\rfloor )$ and
$X_n^+(\lambda n^{1/4})$ are respectively established in the next three
sections (Sections~\ref{section-max} 
  to~\ref{section-limit-global}). The main technique that we use is
  the ``analysis of singularities'' of Flajolet and
Odlyzko~\cite{flajolet-odlyzko}. It permits to extract 
the asymptotic behaviour of the coefficients of a \gf . This technique
has already proved useful in numerous occasions, in particular for
proving limit theorems that are similar in flavour to the ones
obtained in this paper: these theorems deal with the height of simply
generated trees and their profile, which are known to be related to the 
height of the Brownian excursion and its local
time~\cite{flajolet-binary-trees,drmota-profile}. This 
technique is carefully exemplified in Section~\ref{section-max} (which
is devoted to the maximal label) before the more difficult questions
of the local and global limit laws are attacked
(Sections~\ref{section-limit-local}
and~\ref{section-limit-global}). 

Finally, two other families of trees are briefly studied in
Section~\ref{section-universality}: trees with increments $0, \pm 1$
and naturally embedded binary trees.  The emphasis is put on their
enumerative properties, which turn out to be as remarkable and
surprising as those of our first family of trees. The limit laws we
obtain are (up to a scalar) the same as for the first family.

\bigskip
Let us conclude with some notation and a few definitions on \fps \ and
generating functions. 
Let $\GK$ be a field. We denote by $\GK[t]$ the ring of polynomials in
$t$ with coefficients in $\GK$, and by $\GK(t)$ the field of rational
functions in $t$ with coefficients in $\GK$. We denote by $\GK[[t]]$
the ring of \fps\ in $t$ 
with coefficients in $\GK$. If $A(t) \in \GK[[t]]$ and  $n\in\ns$,
the notation $[t^n]A(t)$ stands for the coefficient of $t^n$ in
$A(t)$.
The series $A(t)$ is said to be \emm algebraic over $\GK(t)$, if it satisfies a
non-trivial 
polynomial equation of the form $P(t,A(t))=0$, where $P$ is a
bivariate polynomial with coefficients in $\GK$. In this case, the
\emm degree,\/ of $A(t)$ 
is the smallest possible degree of $P$ (in its second variable).

Let $\A$ be a set of
discrete objects, equipped with a \emm size, that takes 
nonnegative integer values.  Assume
that for all $n \in \ns$, the number of objects of $\A$ of size $n$ is
finite, and denote this number by $a_n$. The {\em \gf\ of the objects of
$\A$, counted by their size,}\/ is the formal power series
$$
A(t)=\sum_{n\ge 0} a_n t^n.
$$
The above notions generalize in a straightforward way to multivariate
power series. Such series  arise naturally when enumerating objects
according to several parameters.

%%%%%%%%%%%%%%%%%%%%%%%%%%%%%%%%%%%%%%%%%%%%%%%%%%
%%%%%%%%%%%%%%%%%%%%%%%%%%%%%%%%%%%%%%%%%%%%%%%%%
\section{Enumerative results}
\label{section-snake-gf}
We consider in this section (and  in the three following ones)
our first family of labelled plane trees: the root is labelled $0$,
and the labels of two adjacent nodes differ by $\pm 1$. 

\subsection{Trees with small labels}
The first enumerative problem we address has already been studied by
Bouttier, Di Francesco and
Guitter~\cite{bdg-geodesic,bdg-statistics}. It deals with the largest
label occurring in a tree.
For $j\in \ns$, let  $T_j\equiv T_j(t)$ be the \gf \ of labelled 
trees in which all labels are less than or equal to $j$. The
indeterminate $t$ keeps track of the number of edges. Let  $T\equiv
T(t)$ be the \gf \ of all labelled trees. Clearly, $T_j$ converges to
$T$ (in the space of \fps \ in $t$) as $j$ goes to infinity. It is
very easy to describe an infinite set of equations that completely
defines the collection of series $T_j$.

\begin{figure}[hbt]
\begin{center}
\input{tree-decomp.pstex_t}
\end{center}
\caption{The decomposition of plane labelled trees.}
\label{fig-decomp}
%\hrule
\end{figure}
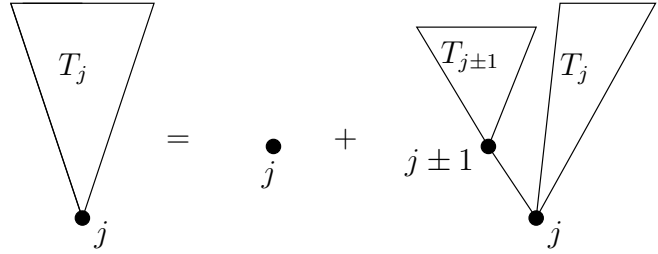

\begin{Lemma} The series $T$ satisfies
\label{lemma-Tj}
\beq
\label{eqT}
T=1+2tT^2.
\eeq
More generally, for $j\ge 0$,
%\beq\label{eqTj}
$$
T_j=1+t(T_{j-1}+T_{j+1})T_j
%\eeq
$$
while $T_j=0$ for $j<0$.
\end{Lemma}
\noindent
{\bf Proof.} The two ingredients of the proof will be useful for the 
other enumerative problems we  address below.
Firstly, replacing each label $k$ by $j-k$ shows that $T_j$ is also the
\gf\ of trees 
\emm rooted at $j$, and having only non-negative labels (we say that a tree
is rooted at $j$ if its root has label $j$).
Secondly, consider such a tree and assume it is not reduced to a
single node. The root has a leftmost child, which  is the root of a 
labelled subtree, rooted at $j\pm 1$ and having only non-negative labels.
   Deleting this subtree leaves a smaller tree rooted at
$j$, having only non-negative labels (see Figure~\ref{fig-decomp}). The
   result follows.
\cqfd

The above lemma shows that the series $T$, counting labelled trees
by edges, is algebraic, and the short proof we have given provides a simple
combinatorial explanation for this property. What  is far less clear
--- but nevertheless true --- is that 
each of the series $T_j$ is algebraic too, as stated in the
proposition below, which we borrow
from~\cite{bdg-geodesic,bdg-statistics}. These series will be 
expressed in terms of the 
series $T\equiv T(t)$ and of the  unique formal power series
$Z\equiv Z(t)$, with constant term $0$, satisfying
\beq
\label{eqZt}
Z=t\, \frac{(1+Z)^4}{1+Z^2}.
\eeq
Observe that $T$ and $Z$ are related by:
\beq
\label{eqTZ}
T=\frac{(1+Z)^2}{1+Z^2}.
\eeq

\begin{Propo}[{\bf Trees with small labels \cite{bdg-geodesic,bdg-statistics}}]
\label{propo-snake-bounded}
Let $T_j\equiv T_j(t)$ be the \gf\ of trees having no label greater
than $j$. Then $T_j$ is algebraic of degree (at most) $2$. In particular,
$$
T_0=1-11\,t-{t}^{2}
+4\,t \left( 3+2\,t \right) {T_0}
-16\,{t}^{2}{{ T_0}}^{2}.
$$
Moreover,  for all $j \ge -1$,
\beq
\label{product-formTj-small}
T_j= T\,\frac{(1-Z^{j+1})(1- Z^{j+5})}{(1- Z^{j+2})(1- Z^{j+4})},
\eeq
where $Z\equiv Z(t)$ is given by~\eqref{eqZt}.
\end{Propo}
\noindent
{\bf Proof.} It is very easy to check, using~(\ref{eqZt}--\ref{eqTZ}),
that the above values of $T_j$ satisfy the recurrence relation of
Lemma~\ref{lemma-Tj} and the initial condition $T_{-1}=0$.
How to \emm discover, such a formula is another story, which is told
in~\cite{bdg-geodesic}. 
The remarkable product form of $T_j$ still awaits a 
combinatorial explanation. 

The equation satisfied by $T_0$ is obtained by eliminating $T$ and $Z$
from the case $j=0$ of~\eqref{product-formTj-small}. Then an induction
of $j$, based on Lemma~\ref{lemma-Tj}, implies that each $T_j$ is
quadratic (at most) over $\qs(t)$. \cqfd

\noindent
{\bf Remarks}\\
{\bf 1.} The product form~\eqref{product-formTj-small}, combined with
the facts that $T$ is quadratic over $\qs(t)$ and $Z$ is 
quadratic over $\qs(T)$, shows that $T_j$ belongs to an extension of
$\qs(t)$ of degree $4$. This is true,  but not optimal, since $T_j$ is
actually quadratic over $\qs(t)$. Hence this product form does not
give the best possible information on the degree of $T_j$.
\\
{\bf 2.} The trees counted by $T_0$ (equivalently, the trees having only
non-negative labels) are known to be in bijection with certain planar
maps called \emm Eulerian triangulations,~\cite{bdg-statistics}. Through
this bijection, the number of edges of the tree is sent to the number
of black faces of the triangulation. These triangulations are nothing
but the dual maps of the \emm bicubic, (that is, bipartite and
trivalent) maps, which were first enumerated by Tutte~\cite{census-maps}.
In particular, the
coefficients of $T_0(t)$ are remarkably simple:
$$
T_0(t)= 
\frac{(1-8t)^{3/2}-1+12t+8t^2}{32t^2}
=1+\sum_{n\ge 1} \frac{3. 2^{n-1}}{(n+1)(n+2)}{{2n}\choose n}t^n.
$$

\subsection{The number of nodes labelled $j$}
Let us now turn our attention to a bivariate counting problem.
For $j \in \zs$, let $S_j \equiv S_j(t,u)$ be the \gf\ of labelled trees,
counted by the number of edges (variable $t$) and the number of
nodes labelled $j$ (variable $u$).  Clearly, $S_j(t,1)=T(t)$ for
all $j$. Moreover, an obvious symmetry entails that $S_j=S_{-j}$. 
\begin{Lemma} For $j\not = 0$,
\label{lemma-Sj-rec}
\beq
\label{eqSj}
S_j=1+t(S_{j-1}+S_{j+1})S_j
\eeq
while for $j=0$,
\beq
\label{eqS0}
S_0= u + t(S_{-1}+S_{1})S_0 =  u +2 tS_{1}S_0.
\eeq
\end{Lemma}
\noindent
{\bf Proof.}
Observe that $S_j \equiv S_j(t,u)$ is also the \gf\ of labelled trees
rooted at $j$,
counted by the number of edges and the number of
nodes labelled $0$. The  decomposition of trees illustrated in Figure~\ref{fig-decomp} then provides the  lemma. The only difference between the cases
   $j=0$ and $j\not =0$ lies in the \gf\ of the tree reduced to a single node.
\cqfd

Again, the series $S_j(t,u)$ turn out to be algebraic, for reasons
that currently remain mysterious (from the combinatorics
viewpoint). They can be expressed in terms of the 
series $T$ and $Z$ given by~(\ref{eqZt}--\ref{eqTZ}).
\begin{Propo}[{\bf The number of nodes labelled $j$}]
\label{propo-snake-gf}
For any $j\in \zs$, the \gf \ $S_j\equiv S_j(t,u)$ that counts
labelled trees by the number of edges and the number of nodes labelled $j$ is
algebraic over $\qs(T,u)$ of degree at most $3$ (and hence has degree
at most $6$ over $\qs(t,u)$). More precisely,
\beq
\label{eqS0T}
\frac{(T-S_0)^2}{(u-1)^2}= 1- \frac{2(1-T^2)}{2+S_0-S_0T},
%
%{S_0}^{3} \left( T-1 \right) -2\,{S_0}^{2} \left( {T}^{2}-T+1 \right) 
%+ S_0 \left((u-1)^2(1-T)+T^3-T^2+4T\right)
%+2\,u{T}^{2} \left(u -2 \right)=0.
\eeq
and all the $S_j$ belong to $\qs(t,u,S_0)$.
Moreover,  for all $j \ge 0$,
\beq
\label{product-formSj}
S_j= T\,\frac{(1+\mu Z^j)(1+\mu Z^{j+4})}{(1+\mu Z^{j+1})(1+\mu Z^{j+3})},
\eeq
where $Z\equiv Z(t)$ is given by~\eqref{eqZt}
and $\mu\equiv \mu(t,u)$ is the unique \fps\ in $t$
satisfying
\beq
\label{eqMZ}
\mu= (u-1) \frac{(1+Z^2)(1+\mu Z)(1+\mu Z^2)(1+\mu Z^3)}
{(1+Z)(1+Z+Z^2)(1-Z)^3(1-\mu Z^2)}.
\eeq
The series $\mu(t,u)$  has polynomial coefficients in $u$,
and satisfies $\mu(t,1)=0$. It has degree $3$ over $\qs(Z,u)$ and
$12$ over $\qs(t,u)$. 
%The first terms in its expansion are:
%$$
%\mu(t,u)=(u-1)\Big(
%1+ut+u \left( u+6 \right) {t}^{2}+u \left( {u}^{2}+14\,u+34 \right)
%{t}^{3}+O(t^4)\Big) .
% u \left( {u}^{3}+24\,{u}^{2}+127\,u+200 \right) {t}^{4}\Big)+O(t^5).
%$$
\end{Propo}
At some point, we will need a closed form expression for $\mu$ in
terms of $Z$. Here is one.
\begin{Propo}
\label{propo-mu-closed-form}
Write
$$
v= \frac{(u-1) Z (1+Z^2)}{(1+Z)(1+Z+Z^2)(1-Z)^3}.
$$
Then the algebraic series $\mu$ involved in the
expression~\eqref{product-formSj} of $S_j$, and defined by~\eqref{eqMZ}, is
$$
\mu(t,u)= \frac 1 {Z^2} \left(    \frac 2
{1+v(1-Z)^2/3+2/3 \sqrt{3+v^2(1-Z)^4}\cos (\phi/3)}
-1 \right)
$$
where
$$
\phi= \arccos \left( \frac{-9v(1+4Z+Z^2)+v^3(1-Z)^6}{(3+v^2(1-Z)^4)^{3/2}}\right).
$$
\end{Propo}
\noindent
    {\bf Proof of Propositions~\ref{propo-snake-gf}
      and~\ref{propo-mu-closed-form}.}\footnote{All the 
      calculations in this paper have been 
      done using {\sc Maple}. We do not recommend  the reader to
      check them by hand.}
First, observe that the family of series $S_0, S_1, S_2, \ldots$ is
completely determined by~\eqref{eqSj} (taken for $j>0$) and the second
part of~\eqref{eqS0}. The fact that for any series $\mu \in \qs(u)[[t]]$, the
expression~\eqref{product-formSj} satisfies~\eqref{eqSj} for all $j>
      0$ is a straighforward 
verification, once $t$ and $T$ have been expressed in terms of $Z$
(see~\eqref{eqZt} and~\eqref{eqTZ}). The form
of~\eqref{product-formSj} is 
borrowed from~\cite{bdg-geodesic}. In order for~\eqref{product-formSj}
to be the correct expression of $S_j$, it remains to
satisfy the second part of~\eqref{eqS0}. This last condition provides a polynomial
equation relating $\mu$, $T$, $Z$, $t$  and $u$. In this equation,
replace $t$ and 
$T$ by their expressions  in terms of $Z$ (given
by~(\ref{eqZt}--\ref{eqTZ})). This gives exactly~\eqref{eqMZ}.
It
can be easily checked that $\mu$ has degree $6$ over $\qs(T,u)$ and
degree $12$ over $\qs(t,u)$.
 
\medskip

The  equation~\eqref{eqS0T} satisfied by $S_0$ is obtained  by
eliminating  $\mu$ and $Z$ (using~\eqref{eqMZ} and~\eqref{eqTZ}) from
the expression~\eqref{product-formSj} of $S_0$.  This equation  gives
 an equation of degree $6$ over $\qs(t,u)$ if one eliminates $T$
thanks to~\eqref{eqT}.

Now the equations~\eqref{eqS0}, \eqref{eqSj} and~\eqref{eqT}, combined with an
induction on $j$, imply that for $j\ge
1$, the series $S_j$ belongs to the field  $\qs(T,u,S_0)$,
which has just been proved to be an extension of $\qs(T,u)$ of degree
  $3$. This concludes the proof of  Proposition~\ref{propo-snake-gf}.

\smallskip
Let us finally prove Proposition~\ref{propo-mu-closed-form}. The
equation~\eqref{eqMZ} that defines $\mu$ can be rewritten
$$
\mu= \frac v Z \frac{(1+\mu Z)(1+\mu Z^2)(1+\mu Z^3)}
{1-\mu Z^2}.
$$
Hence $\mu $ is the unique \fps \ in $v$ (with rational coefficients
in $Z$) that satisfies the above
equation and equals $0$ when $v$ is $0$. It is not hard to check that
the closed form expression we give satisfies these two conditions.
\cqfd

\noindent{\bf Remarks}\\
{\bf 1.} The product form~\eqref{product-formSj} of
Proposition~\ref{propo-snake-gf}  refines the
product form~\eqref{product-formTj-small} that deals with trees with small
labels. Indeed, when $u=0$, Eq.~\eqref{eqMZ} gives $\mu=-1$, and the
expression of $S_j(t,0)$ coincides, as it should, with the expression
of $T_{j-1}(t)$ given by Proposition~\ref{propo-snake-bounded}.\\
{\bf 2.} There exists an alternative way to
derive an equation for  $S_0$ from the system of Lemma~\ref{lemma-Sj-rec}.
 As was observed in~\cite[p.~645]{bdg-statistics} for the 
problem of counting trees with bounded labels, 
Eq.~\eqref{eqSj} implies that for $j\ge 1$,
$$
I(S_{j-1},S_j)=I(S_{j},S_{j+1})
$$
where the ``invariant'' function $I$ is given by
$$
I(x,y)= xy(1-tx)(1-ty)+txy-x-y.
$$
But $S_j$ converges to $T$ as $j$ goes to infinity, in the set of
\fps\ in $t$. This implies
$$
I(S_{0},S_1)=I(T,T).
$$
Eliminating $S_1$ between the above equation and~\eqref{eqS0} gives an
equation between $S_0$, $T$ and $t$.
%%%%%%%%%%%%%%%%%%%%%%%%%%%%%%%%%%%%%%%%%%%%%%%%%%

\subsection{The number of nodes labelled $j$ or more}
Let us finally study our third and last enumeration problem. For $j \in
\zs$, let $R_j \equiv R_j(t,u)$ be the \gf\ of labelled trees, 
counted by the number of edges (variable $t$) and the number of
nodes labelled $j$ at least (variable $u$). 
\begin{Lemma}
\label{lemma-Rj-rec}
The set of series $R_0, R_1, R_2, \ldots$ is completely determined by the
following equations: for $j\ge 1$,
\beq
R_j=1+tR_j(R_{j-1}+R_{j+1} ) 
\label{rec-Rj}
\eeq
and
\beq
\label{R0-R1}
 R_0(t,u)=uR_1(tu,1/u).
\eeq
More generally, for all $j \in\zs$, one has:
\beq
\label{eq-symmetry}
 R_{-j}(t,u)=uR_{j+1}(tu,1/u).
\eeq
\end{Lemma}
\noindent {\bf Proof.} For all $j\in \zs$, the series  $R_j\equiv
R_j(t,u)$ is also the \gf \ of 
trees rooted at $j$, counted by their number of edges
and the number of nodes having a non-positive label. The equation
satisfied by $j$, for $j \ge 1$, 
follows once again from the decomposition of trees illustrated in
Figure~\ref{fig-decomp}. It remains to prove the symmetry
relation~\eqref{eq-symmetry}. For any tree $\tau $, let $n_{\le 0}(\tau )$
denote the number of nodes of $\tau $ having a non-positive label. We 
use similar notations for the number of nodes having label at most
$j$, etc. Let $\T_{j,n}$ denote the set of trees rooted at
$j$ and having $n$ edges.
As observed above,
$$
R_{-j}(t,u)= \sum_{n\ge 0} t^n \sum_{\tau  \in \T_{-j,n}} u^{n_{\le 0}(\tau )}
= \sum_{n\ge 0} t^n \sum_{\tau  \in \T_{-j,n}} u^{n+1-n_{> 0}(\tau )},
$$
because a tree with $n$ edges has a total of $n+1$ nodes. A
translation of all labels by $-1$ gives
$$
R_{-j}(t,u)=u \sum_{n\ge 0} (tu)^n \sum_{\tau  \in \T_{-j-1,n}} u^{-n_{\ge
    0}(\tau )},
$$
while replacing each label $k$ by $-k$ finally gives
$$
R_{-j}(t,u)=u \sum_{n\ge 0} (tu)^n \sum_{\tau  \in \T_{j+1,n}} u^{-n_{\le
    0}(\tau )}=u R_{j+1}(tu,1/u).
$$
\cqfd

Again, the series $R_j$ are algebraic, and admit a closed form
expression in terms of $T$ and $Z$.
\begin{Propo}[{\bf  The number of nodes labelled $j$ or more}]
\label{propo-snake-gf-rep}
Let $j \in \zs$. The \gf\ $R_j(t,u)\equiv R_j$
that counts labelled trees by the number of edges and the number of
nodes labelled $j$ or more
is  algebraic of degree at most $2$ over $\qs(T(t),T(tu))$. Hence it
has degree at most $8$ over $\qs(t,u)$. More precisely, it belongs to
the extension of $\qs(T(t),T(tu))$ generated by 
$$
%\sqrt\Delta_1= 
\sqrt{(T+\tilde T)^2-4T\tilde T (T-1)(\tilde T -1)}
$$
where $T\equiv T(t)$ and $\tilde T \equiv T(tu)$. \\
Moreover,  for all $j \ge 0$,
\beq
\label{product-formRj}
R_j= T\,\frac{(1+\nu Z^j)(1+\nu Z^{j+4})}{(1+\nu Z^{j+1})(1+\nu Z^{j+3})},
\eeq
where $Z\equiv Z(t)$ is given by~\eqref{eqZt}
and $\nu\equiv \nu(t,u)$ is a \fps\ in $t$, with polynomial
coefficients in $u$, which is algebraic
%  of degree $2$ over $\qs(Z(t),Z(tu))$, 
of degree $4$ over $\qs(u,Z)$, and  of degree $16$
over $\qs(t,u)$. 
This series satisfies $\nu(t,1)=0$. 
The first terms in its expansion are:
$$
\nu(t,u)=(u-1)\Big(
1+2\,ut+ \left(7\,u+ 6\,{u}^{2} \right) {t}^{2}
+ \left( 32\,u+36\,{u}^{2}+23\,{u}^{3} \right) {t}^{3}+
O(t^4)\Big) .
$$
\end{Propo}
Before we prove this proposition, let us give something like a closed form for
$\nu$. Since $\nu$ has degree 4 over $\qs(u,Z)$, and $Z$ has degree 4
over $\qs(t)$, the series $\nu$ is in theory expressible in terms of
radicals... It turns that this expression is less terrible than one
could fear.
\begin{Propo}
\label{propo-nu}
Define the following four \fps\ in $t$ with polynomial coefficients
in $u$:
$$
\delta\equiv\delta(t,u)= 1-8(u-1) \frac{Z(1+Z^2)}{(1-Z)^4} 
= \frac{1-8tu}{1-8t},
$$
$$
V\equiv V(t,u) = \frac{1-\sqrt \delta} 4 =
 \frac{1-\sqrt {\frac{1-8tu}{1-8t}}} 4,
$$
$$
\Delta\equiv \Delta(t,u)= (1-V)^2 - \frac{4Z V^2}{(1+Z)^2},
$$
and 
$$
P=(1+Z) \frac{1-V-\sqrt \Delta}{2VZ}.
$$
Then $P$ has degree $16$ over $\qs(t,u)$, degree $2$ over $\qs(V,Z)$,
and satisfies the following ``Lagrangian'' equation:
$$
P=\frac V{1+Z} (1+P)(1+ZP).
$$
Moreover,  the algebraic series $\nu$ involved in the
expression~\eqref{product-formRj}  of $R_j$ is
$$
\nu= \frac P Z \frac {1-P(1+Z)-P^2(1+Z+Z^2)}{1+Z+Z^2+PZ(1+Z)-P^2Z^2}.
$$
\end{Propo}
\noindent
{\bf Proof of Proposition~\ref{propo-snake-gf-rep}.} We have already
checked, in the proof of 
Proposition~\ref{propo-snake-gf}, that for any \fps \ $\nu$ in $t$, the series
defined by~\eqref{product-formRj} for $j\ge 0$ satisfy  the recurrence
relation~\eqref{rec-Rj} for $j\ge 1$. It remains to prove that one can
choose $\nu$ so as to satisfy~\eqref{R0-R1}. For any \fps\
$A$ in $t$  having rational coefficients in $u$, we denote by
$\tilde A$ the series $\tilde A(t,u) =A(tu,1/u)$. Observe that $\tilde{\tilde
  A}=A$. With this notation, if $R_j$ is of the generic
form~\eqref{product-formRj}, 
the relation~\eqref{R0-R1} holds if and only if
\beq
\label{nu-nutilde}
 1+\nu =
u\, \frac{\tilde T} T\,\frac
{(1+\nu Z)(1+\nu Z^{3})(1+ \tilde \nu \tilde Z)(1+\tilde \nu \tilde Z^{5})}
{(1+\nu Z^{4})(1+\tilde \nu \tilde Z ^{2})(1+\tilde \nu \tilde Z^{4})}.
\eeq
Let $\rs_m[u]$ denote the space of polynomials in $u$, with real
coefficients, of degree at most $m$. 
Let $\rs_n[u][[t]]$ denote the set of \fps \ in $t$ with polynomial
coefficients in $u$ such that for all $m \le n$, the coefficient of
$t^m$ has degree at 
most $m$. Observe
that this set of series in stable under the usual operations on
series: sum, product, and quasi-inverse.
Write  $\nu=\sum_{n\ge 0} \nu_n(u) t^n$. We are going to prove, by
induction on $n$,  that~\eqref{nu-nutilde} determines uniquely each
coefficient $\nu_n(u)$, and that this coefficient belongs to $\rs_{n+1}[u]$.

First, observe that for any \fps\ $\nu$, the right-hand side of~\eqref{nu-nutilde}  is $u+
O(t)$. This implies $\nu_0(u)=u-1$. Now assume that our induction hypothesis
holds for all $m<n$. Recall that $Z$ is a multiple of $t$:
this implies that  $ \nu  Z$ belongs to  $\rs_n[u][[t]]$. The
induction hypothesis also implies that the coefficient 
of $t^m$ in $u\tilde \nu$ belongs to $\rs_{m+1}[u]$, for all $m<n$. 
Note that $\tilde Z = Z(tu)=tu+O(t^2)$ is a
multiple of $t$ and $u$ and also belongs to  $\rs_n[u][[t]]$. This
implies that  $\tilde \nu \tilde Z$ 
belongs to $\rs_n[u][[t]]$ too.
 The same is true for all the other series
occurring in the right-hand side of~\eqref{nu-nutilde}, namely $T,
\tilde T, Z, \tilde Z$. Given the closure properties of the set
$\rs_n[u][[t]]$, we conclude that the right-hand side
of~\eqref{nu-nutilde}, divided by $u$, belongs to this set. Moreover,
the fact that $Z$ and $\tilde Z$ are multiples of $t$ guarantees that
the coefficient of $t^n$ in this series only involves the $\nu_i(u)$
for $i<n$.
  By extracting the coefficient
of $t^n$ in~\eqref{nu-nutilde}, we conclude that $\nu_n(u)$ is
uniquely determined and belongs to  $u\rs_n[u]\subset\rs_{n+1}[u] $. 

This completes the proof of the existence and uniqueness of the series
 $\nu$  satisfying~\eqref{nu-nutilde}. 
 Also, setting $u=1$ (that is, $\tilde T=T$ and $\tilde Z=Z$) in this
 equation shows that $\nu(t,1)=0$.

Let us now replace $t$ by $tu$ and $u$ by $1/u$
in~\eqref{nu-nutilde}. This gives:
\beq
\label{nutilde-nu}
 1+\tilde \nu =\frac 1
u\, \frac{ T} {\tilde T} \,\frac
{(1+\tilde \nu \tilde Z)(1+\tilde \nu \tilde Z^{3})(1+ \nu  Z)(1+ \nu  Z^{5})}
{(1+\tilde \nu \tilde Z^{4})(1+ \nu  Z ^{2})(1+ \nu  Z^{4})}.
\eeq
In the above two equations, replace $T$ by its expression~\eqref{eqTZ}
in terms of $Z$. Similarly, replace $\tilde T$ by its expression
in terms of $\tilde Z$. Finally, it follows from~\eqref{eqZt} and from
the fact that $\tilde Z = Z(tu)$ that
\beq
\label{uZZtilde}
u= \frac {\tilde Z} Z \, \frac{(1+Z)^4(1+\tilde Z^2)}
{(1+\tilde Z)^4(1+Z^2)}.
\eeq
Replace $u$ by this expression in~\eqref{nu-nutilde} and~\eqref{nutilde-nu}.
Eliminate $\tilde \nu$ between the resulting
two equations: this gives a polynomial equation that relates $\nu, Z$
and $\tilde Z$, of degree $2$ in $\nu$. The elimination of $\tilde Z$
between this quadratic equation and~\eqref{uZZtilde} provides an
equation of degree $4$ in $\nu$ that relates $\nu$ to $Z$ and $u$.
Finally, the elimination of $ Z$ shows that $\nu$ is algebraic of
degree $16$ over $\qs(t,u)$.

\medskip
Let us now focus on the first part of the proposition.
From the form~\eqref{product-formRj}, and the fact that $\nu$ has
degree $4$ over $\qs(u,Z)$ and $Z$ has degree $4$ over $\qs(t)$, we
conclude that the degree of $R_j$ over $\qs(t,u)$ is a divisor of
$16$. Let us prove that is is actually a divisor of $8$. The proof
goes as follows:
\begin{enumerate}
\item  Using the generic form~\eqref{product-formRj}, and the equations
satisfied by $T, Z$ and $\nu$, we obtain a polynomial equation of
degree $8$ over $\qs(t,u)$ for $R_0$.

\item  Using~\eqref{eqT} to express $t$ in terms of $T$, and 
$$
u=\frac {T^2}{\tilde T ^2} \, \frac{1-\tilde T}{1-T},
$$
(which also follows from~\eqref{eqT}),
we convert the equation satisfied by $R_0$ into a polynomial equation
(still of degree $8$ in $R_0$) relating $R_0$  to $T$ and $\tilde
T$. This equation factors into four quadratic polynomials in
$R_0$. The factor that actually vanishes is identified by setting $u=1$
(in which case $\tilde T=T=R_0$).

\item  From this equation, we conclude that  $R_0$ belongs
to the extension of $\qs(T,\tilde T)$ generated by
$$
\sqrt\Delta_1= \sqrt{(T+\tilde T)^2-4T\tilde T (T-1)(\tilde T -1)}.
$$  
Observe that this extension of $\qs(t,u)$ is left invariant by the
transformation $A\mapsto \tilde A$.

\item  From the fact that $R_1=u \tilde{R_0}$ (see~\eqref{eq-symmetry}),
 we conclude that $R_1$ also belongs to $\qs(T,\tilde
 T,\sqrt\Delta_1 )$.

\item  The recurrence relation~\eqref{rec-Rj} on the $R_j$ allows us to
extends this to all $R_j$, for $j\ge 0$.

\item  Finally,~\eqref{eq-symmetry} shows that our algebraicity result
actually holds for all $R_j$, for $j \in \zs$.
\end{enumerate}
\cqfd

\noindent{\bf Proof of Proposition~\ref{propo-nu}.}
In the course of the proof of Proposition~\ref{propo-snake-gf-rep}, we
have obtained a polynomial equation $P(\nu,Z,u)=0$, of degree 4 in
$\nu$, relating the series $\nu(t,u), Z(t)$, and the variable
$u$. This equation is not written in the paper (it is a bit too big),
but it follows from~\eqref{nu-nutilde} and~\eqref{nutilde-nu}. In this
equation, replace $u$ by its expression in terms of $\delta$ and
$Z$. Then replace $\delta$ by its expression in terms of $V$:  the
resulting equation factors into two terms! Each of them is quadratic
in $\nu$. In order to decide which of these factors cancels, one uses
the fact that when $u=1$ (that is, $V=0$), the series $\nu$ must be
$0$.
It remains to solve a quadratic equation in $\nu$. Its discriminant is
found to be $\Delta$, and one may find convenient to introduce the
series $P$ which is Lagrangian in $V$.
\cqfd

\noindent{\bf Remark.} Again, the product form~\eqref{product-formRj}
of Proposition~\ref{propo-snake-gf-rep} includes as a special case
the enumeration of trees with labels at most $j-1$, obtained when
$u=0$. Indeed,~\eqref{nu-nutilde} shows that $\nu=-1$ when $u=0$,
and~\eqref{product-formRj} 
 then reduces to~\eqref{product-formTj-small}.

%%%%%%%%%%%%%%%%%%%%%%%%%%%%%%%%%%%%%%%%%%%%%%%%%
\section{The largest label, and the support of the ISE}
\label{section-max}
%%%%%%%%%%%%%%%%%%%%%%%%%%%%%%%%%%%%%%%%%%%%%%%%%
Let $\mathcal T _0$ denote the set of labelled trees (rooted at $0$),
and let $\mathcal T _{0,n}$ denote the subset of $\mathcal T _0$ 
formed by trees having $n$ edges.
We  endow $\mathcal T _{0,n}$  with the
uniform distribution. In other words, any of its elements occurs with
probability
$$
\frac 1{2^n C_n}
$$
where $C_n= \frac{1}{n+1} {{2n} \choose n}$ is the $n$th Catalan
number, and is well-known to be the number of (unlabelled) plane trees
with $n$ edges.

\begin{figure}[phtb]
  \begin{center}
    \includegraphics[height=3cm,width=3cm]{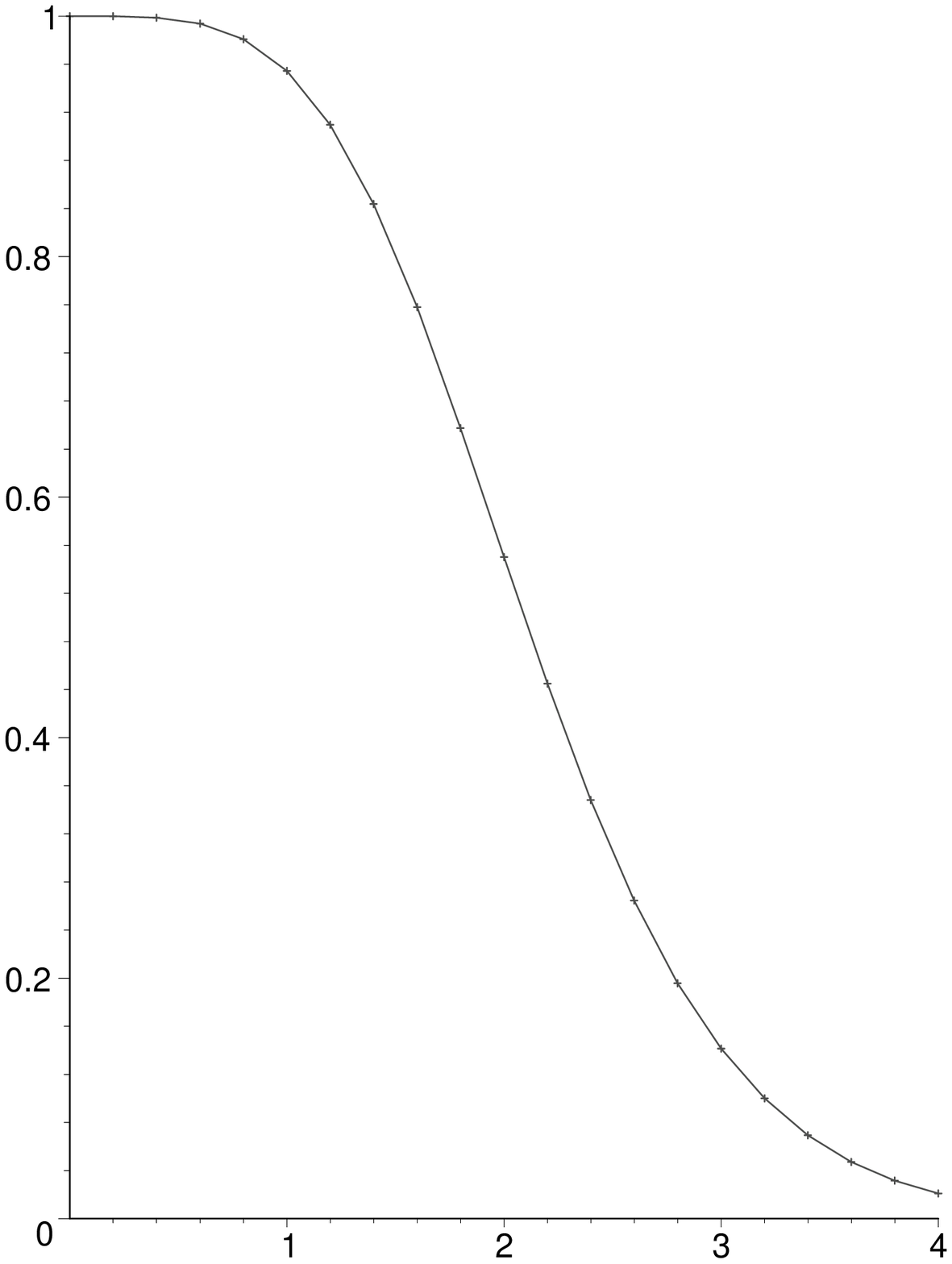}
\hskip 28mm
\includegraphics[height=3cm,width=3cm]{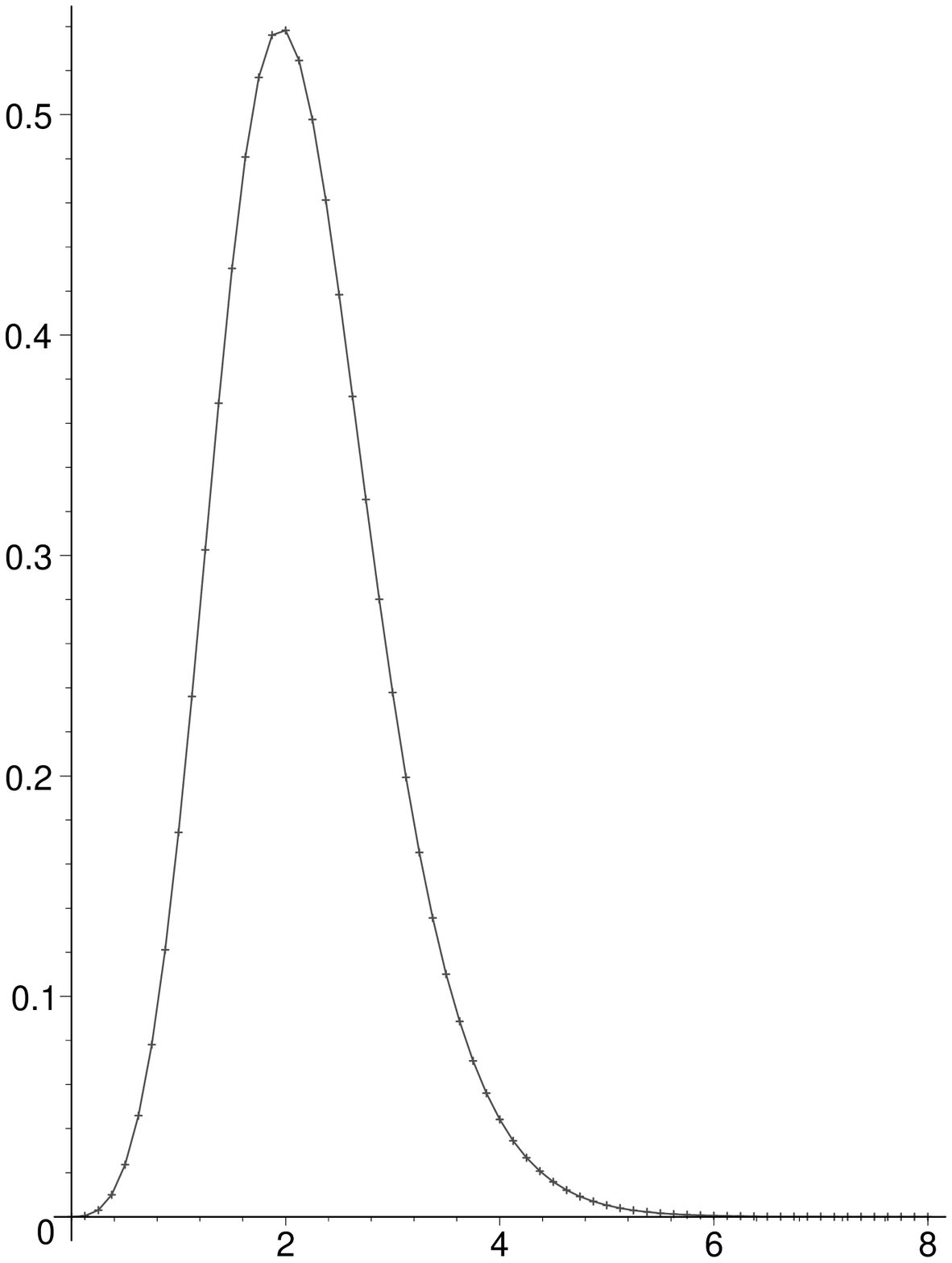}
    \caption{The tail distribution function $G$ and the
    density $f$ of the    limit distribution $N$.} 
    \label{figure-densite}
  \end{center}
\end{figure}

 Let $M_n$ denote the random variable equal to the
largest label occurring in a random tree of  $\mathcal T _{0,n}$. 
The law of $M_n$ is related to the series $T_j$ studied in
Proposition~\ref{propo-snake-bounded}:
$$
\PP\left(M_n\le j\right)=  \frac{[t^n] T_j}{2^n C_n}.
$$
Let us define a normalized version of $M_n$ by
$$
N_n= \frac{M_n}{n^{1/4}}.
$$
The aim of this section is to prove the convergence of $N_n$ in
distribution\footnote{The above convention will be used throughout the
  paper: if a random variable depending on $n$ is denoted by some
  letter of the alphabet, then its suitably normalized version is
  denoted by \emm the next letter, of the alphabet.}.
\begin{Theorem}
\label{thm-max}
As $n$ goes to infinity, the random variable $N_n$ converges
 in distribution to a non-negative random
variable $N$. The tail distribution function of $N$, defined by
$ G(\lambda)=\PP(N > \lambda )$, satisfies
$$
G(\lambda)=
\frac {12}{i\sqrt \pi}\int_ \Gamma\frac{v^5 e^{v^4}}{\sinh^2(\lambda v) }
dv
=\frac{6}{\sqrt \pi \lambda ^6}\int_0^\infty\frac{1-\cos u \cosh
  u}{(\cosh u -\cos u)^2} u^5 e^{-u^4/(4\lambda^4)}du
$$
where the contour $\Gamma$ is formed of two half-lines:
$$
\Gamma= \{1-te^{-i\pi/4}, t\in (\infty,0]\} \cup \{1+te^{-i\pi/4},
t\in [0,\infty)\}. 
$$
Equivalently, the variable $N$ has density
$$
f(\lambda)=\frac {24}{i\sqrt \pi}\int_ \Gamma\frac{\cosh(\lambda v)v^6
  e^{v^4}}{\sinh^3(\lambda v) } dv
=\frac{6}{\sqrt \pi \lambda ^{11}}\int_0^\infty\frac{1-\cos u \cosh
  u}{(\cosh u -\cos u)^2} u^5(6\lambda^4-u^4) e^{-u^4/(4\lambda^4)}du
$$
with respect to the Lebesgue measure on $\rs_+$.
The moments of $N$ are finite, and admit simple expressions:
$$
\EE(N)=\frac{3\sqrt \pi}{2\Gamma(3/4)},\quad
\EE(N^2)=3\sqrt \pi,
$$
 and for  $k\ge 3$, 
$$
\EE(N^k)=\frac{24 \sqrt \pi k! \zeta(k-1)}{2^k \Gamma ((k-2)/4)}.
$$
Finally, the moments of $N_n=M_n/n^{1/4}$ converge to the moments of $N$.
\end{Theorem}
\noindent The functions $G$ and $f$ are plotted in Figure~\ref{figure-densite}. 

\noindent The proof of this theorem will be split into four
subsections (Sections~\ref{section-distribution-function}
to~\ref{section-moments-Nn}). In view of the following proposition,
this theorem gives 
the density, distribution function and moments of the supremum of  the
support of the ISE. 

\begin{Propo}[{\bf The supremum of the support of the ISE}]
\label{coro-support}
Let $N_{{\hbox{\sc \small ise}}}$ denote the supremum of the support of
the ISE
$$
N_{{\hbox{\sc \small ise}}}= \sup\{y: \mu_{\hbox{\sc \small ise}}(y,
\infty)>0\}.
$$
Then 
$ N_{{\hbox{\sc \small ise}}}$ has the same law as the random variable
$\sqrt 2N$ described in  Theorem~{\em\ref{thm-max}.} 
\end{Propo}

\noindent{\bf Remark.} The moments of $N_{\hbox{\sc \small ise}}$ 
are thus
$$
\EE(N_{{\hbox{\sc \small ise}}})=\frac{3\sqrt \pi}{\sqrt 2\Gamma(3/4)},\quad
\EE(N_{{\hbox{\sc \small ise}}}^2)=6\sqrt \pi,
$$
and for  $k\ge 3$, 
$$
\EE(N_{{\hbox{\sc \small ise}} }^k)=\frac{24 \sqrt \pi k! \zeta(k-1)}{\sqrt 2^k \Gamma ((k-2)/4)}.
$$
They were
already obtained by Delmas~\cite{delmas} using a completely different
(and continuous) approach.  The expressions he gives actually differ
from ours by a factor $2^{k/4}$, due to a different choice of
normalization. Note that the zeta function also appears in the moments
of the maximum of the Brownian excursion, which follows a theta
law~\cite{chung}. This law is known to describe the limiting
normalized height of  simple trees~\cite{flajolet-binary-trees}.
Finally, let us mention that another,  more complicated
 expression of the density of the limiting
 variable $N$ was obtained
 in~\cite{bdg-geodesic} (maybe in a slightly less rigorous
 way). Proposition~\ref{coro-support} is proved in
 Section~\ref{section-supremum}.

\subsection{Convergence of the distribution function}
\label{section-distribution-function}
We prove in this section that the tail distribution function of $N_n$
converges pointwise.
Let $\lambda \ge 0$ and $j=\lfloor \lambda n^{1/4} \rfloor$. 
The probability we are interested in is
\beq
\label{GU}
\PP(N_n > \lambda )=\PP(M_n > \lambda n^{1/4})=\PP(M_n >j)=\frac{[t^n] U_j(t)}{2^n C_n},
\eeq
where
\beq
\label{Uj-expr}
U_j(t)\equiv U_j = T- T_j=  
\frac { \left( 1+Z \right) ^{2}{Z}^{j+1} \left( 1+Z+{Z}^{2} \right)  
\left(1- Z \right) ^{2}}
{ \left( 1+{Z}^{2} \right) 
 \left( 1-{Z}^{j+2} \right)  \left( 1-{Z}^{j+4} \right) }
 \eeq
is the \gf  \ of trees having at least one label greater than
$j$. This algebraic series has a positive radius of
convergence\footnote{So do all  algebraic power series}, 
and by Cauchy's formula,
\begin{eqnarray}
[t^n] U_j& =&\frac 1 {2i\pi} \int_\C U_j(t)
\frac{dt}{t^{n+1}} \nonumber\\
&=&\frac {1} {2i\pi} \int_\C 
\frac { \left( 1+Z \right) ^{2}{Z}^{j+1} \left( 1+Z+{Z}^{2} \right)  
\left(1- Z \right) ^{2}}
{ \left( 1+{Z}^{2} \right) 
 \left( 1-{Z}^{j+2} \right)  \left( 1-{Z}^{j+4} \right) }
\frac{dt}{t^{n+1}},
\label{cauchy}
\end{eqnarray}
for any contour $\mathcal C$ included in the analyticity domain of
$U_j$ and enclosing positively the origin.

This leads us to study  the singularities of $U_j$, and therefore
 those of  $Z$. 
 We  gather in the following lemma a few  properties of this series.

\begin{Lemma}[{\bf Analytic properties of $Z$}]
\label{lemma-Z}
Let $Z\equiv Z(t)$ be the unique \fps \ in $t$ with constant term $0$
satisfying~\eqref{eqZt}.
This series has non-negative integer coefficients.
It has radius of convergence $1/8$, and can be continued analytically on 
the domain $\D=\cs\setminus[1/8, +\infty)$. In the neighborhood of
  $t=1/8$, one has
\beq
\label{Z-asympt}
Z(t)=1-2(1-8t)^{1/4} + O(\sqrt{1-8t}).
\eeq
Moreover,  $|Z(t)|<1$ on the domain $\D$.
More precisely, the only roots of unity 
that are  accumulation points of the set $Z(\D)$ are $1$ and
$-1$, and they are only approached by $Z(t)$ when $t$ tends to $1/8$
and when $|t|$ tends to $\infty$, respectively.
\end{Lemma}
\noindent
{\bf Proof.}
In order to establish the first statement, we observe that
$$
Z=W(1+Z)^2
$$
where $W\equiv W(t)$ is the only \fps\ in $t$ with constant term zero
satisfying
\beq
\label{eqW}
W=t+2W^2.
\eeq
These equations imply that both $W$ and $Z$ have non-negative integer
coefficients.

\medskip
The general approach for studying the singularities of
algebraic series (see for instance~\cite{fla-poly-alg}) gives
the second part of the lemma (up to~\eqref{Z-asympt}).  The polynomial
equation defining $Z(t)$ has leading coefficient $t$ and
discriminant  $4(1-8t)^3$, so that the
only possible singularity of $Z$ is $1/8$. Alternatively,
one can exploit the following closed form expression:
\beq
\label{Z-expr}
Z(t)=\frac{\sqrt{1-4t+\sqrt{1-8t}}\left(\sqrt{1-4t+\sqrt{1-8t}} -\sqrt
  2 (1-8t)^{1/4}\right) }{4t}. 
\eeq

Let us now come to the third part of the lemma, and prove that
  $|Z(t)|$ never reaches $1$ on the
  domain $\D$. Assume $Z(t)=e^{i\theta}$, with
  $\theta\in[-\pi,\pi]$. From~\eqref{eqZt}, one has 
$$
t=t_\theta \quad \hbox{ where } \quad t_\theta=
 \frac{\cos\theta}{8\cos^4(\theta/2)} \quad \hbox{ and } \quad
 \theta\in(-\pi,\pi). 
$$
This shows that $t$ is real, and belongs to $(-\infty, 1/8)$. 
But the expression~\eqref{Z-expr} of $Z(t)$ shows that $Z(t)$ is real,
which contradicts the 
hypothesis $Z(t)=e^{i\theta}$, unless $\theta=0$. But then $t=1/8$
and does not belong to the domain $\D$. Hence the modulus of $Z$ never
reaches $1$ on $\D$.
One can actually prove that, for $\theta \in (-\pi, 0)$,
$$
Z(t_\theta)= \frac{1+ \sin \theta}{\cos \theta},
$$
but we do not need so much precision here.
% Conversely, let us try to evaluate $Z(t_\theta)$ for
% $t=\theta\in(-\pi,0)$ (by parity, this is sufficient for $t$ to
% cover the whole interval $(-\infty, 1/8)$). Can it be equal to $
% e^{i\theta}$? 
%The four roots of the equation (in $Z$)
%$$
%Z=t-\theta \frac{(1+Z)^4}{1+Z^2}$$
%are
%$$
%Z_{1,2}= e^{\pm i \theta} \quad\hbox{ and } \quad
%Z_{3,4}= \frac{1\pm\sin \theta}{\cos \theta}.
%$$
%%
%The expression~\eqref{Z-expr} shows that $Z(t)$ is real when $t\in
%(-\infty, 1/8)$. Therefore the solutions $Z_{1,2}$ cannot coincide
%with $Z(t)$. Moreover, for $|t|<1/8$, the value of $Z(t)$ is given by
%a \fps\ with non-negative coefficients, so that $Z(t)$ is an \emm
%increasing, function of $t$ on $[0,1/8)$. This rules out  the
%possibility that $Z(t)=Z_4$ (the solution with a ``minus'' sign),
%andwe conclude that $$ Z(t_\theta)=  \frac{1+ \sin \theta}{\cos
%\theta}.$$This function of $\theta$ increases from $-1$ to $1$ as
%$\theta$ goesfropm $-\pi$ to $0$, and thus has modulus less than $1$
%for $t\in(-\infty, 0)$. Consequently, $Z(t\theta)\not = e^{i \theta}$
%and more generally $|Z(t)|<= 1$ for $t \in \cs\setminus[1/8, +\infty)$.

Finally, if a sequence $t_n$ of $\D$ is such that $Z(t_n) \rightarrow
e^{i\theta}$ as $n \rightarrow \infty$, with $\theta\in(-\pi, \pi]$,
then either $\theta=\pi$ and, by~\eqref{eqZt}, the
sequence $|t_n|$  tends to $\infty$, or $\theta\in (-\pi, \pi)$
  and $t_n$ converges to $t_\theta$. But then by continuity, $Z(t_n)$
actually converges to $Z(t_\theta)$, which, as argued above, only
coincides with $e^{i\theta}$ when $\theta=0$, that is, $t_\theta=1/8$. In this
case, $Z(t_n)\rightarrow 1$.
\cqfd

\begin{figure}[hbt]
\begin{center}
\input{contour.pstex_t}
\end{center}
\caption{The integration contour $\C_n$.}
\label{fig-contour}
%\hrule
\end{figure}
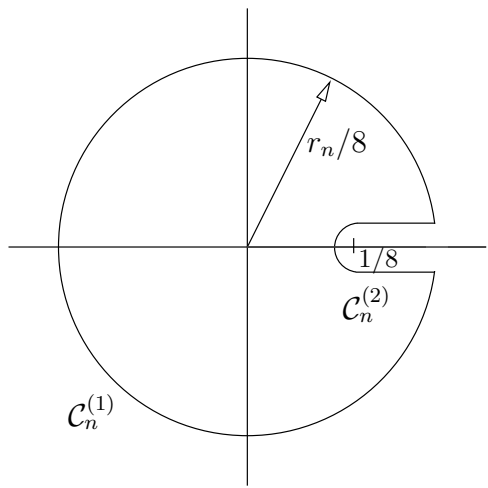

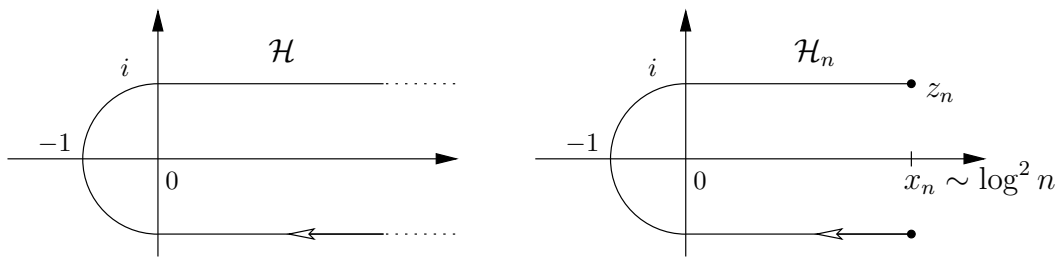
\begin{figure}[b]
\begin{center}
\input{hankel.pstex_t}
\end{center}
\caption{The Hankel contour $\Ha$ and its truncated version $\Ha_n$.}
\label{fig-hankel}
%\hrule
\end{figure}

Let us now go back to the evaluation of the tail distribution function
of $N_n$ via the integral~\eqref{cauchy}.
%The  argument is now routine for amateurs of singularity analysis.
 We choose a contour $\C={\mathcal C}_n$ that depends on $n$ and consists of
 two parts ${\mathcal C}_n^{(1)}$ and ${\mathcal C}_n^{(2)}$ (see
 Figure~\ref{fig-contour}): 
\begin{itemize}
\item[$\bullet$ ] ${\mathcal C}_n^{(1)}$ is an arc  of radius
$r_n/8=(1+ \log^2 n /n)/8$, centered at the origin; note that its
  radius tends to $1/8$ as $n$ goes to infinity,
\item[$\bullet$ ] ${\mathcal C}_n^{(2)}$ is a Hankel contour around $1/8$, at distance
$1/(8n)$ of the real axis, 
which meets  ${\mathcal C}_n^{(1)}$ at both ends; this contour shrinks
around $1/8$  as $n$ goes to infinity; more precisely, as $t$ runs
along $\C_n^{(2)}$, the variable $z$ defined by
$$
t=\frac 1 8 \left(1+ \frac z n\right)
$$
runs over the truncated Hankel contour $\Ha_n$ shown on the right
 of Figure~\ref{fig-hankel}:
$$
{\Ha_n}= \left\{  {x-i}, \ x \in [0, x_n] \right\} \cup
\left\{-{e^{i\theta}}, \  \theta \in [-\pi/2, \pi/2]\right\}
\cup
\left\{ {x+i}, \ x \in [0, x_n] \right\} 
$$
where $(1+x_n/n)^2+1/n^2=r_n^2$, so that $x_n \le \log ^2 n$ 
and $x_n=\log^2n+O(1/n)$.
\end{itemize}
We denote by $z_n=x_n+i$ the top right end of $\Ha_n$. This point
tends to infinity as $n$ does.

 The integral~\Ref{cauchy} on $\C=\mathcal C _n$ is the sum of the
contributions of the 
contours $\C_n^{(1)}$ and $\C_n^{(2)}$. We shall see that the dominant
contribution is that of ${\mathcal C}_n^{(2)}$, because of the
vicinity of the singularity at $t=1/8$. 

Let us first bound carefully $Z(t)$ for $t \in \C_n$. Let $t_n \in
\C_n$ be such that
$$
|Z(t_n)|=\max_{t\in \C_n} |Z(t)|.
$$ 
By Lemma~\ref{lemma-Z}, $|Z(t_n)|$ tends to $1$ as $n$
grows. Moreover, every accumulation point $a$ of the sequence $t_n$
satisfies $|a|\le 1/8$ and $|Z(a)|=1$. This forces $a=1/8$, and we
conclude that $t_n\rightarrow 1/8$.  Write $t_n=(1-u_n)/8$. Then $u_n
\rightarrow 0$, but $|u_n|\ge 1/n$. By~\eqref{Z-asympt},
$$
Z(t_n)=1-2u_n^{1/4} \left( 1 +o(1)\right).
$$
Let us write, for short, $v_n=1-Z(t_n)$. Then $v_n \rightarrow 0$ but
\beq
\label{1-Z-bound}
|v_n| =2|u_n|^{1/4}\left( 1 +o(1)\right)\ge n^{-1/4}
\eeq
for $n$ large enough. Moreover,
$$
|\arg v_n| =\frac 1 4 |\arg(u_n)| +o(1) \le \frac \pi 4 + o(1),
$$
so that
$$
\cos (\arg v_n) \ge \frac 1 {\sqrt 2} +o(1).
$$
Finally,
$$
|Z(t_n)|^2= |1-v_n|^2 = 1-2|v_n| \cos (\arg v_n) + |v_n|^2
\le 1 -\sqrt 2 |v_n|\left( 1 +o(1)\right),
$$
that is,
$$
|Z(t_n)|\le\displaystyle 1- \frac 1 {\sqrt 2} |v_n|\left( 1 +o(1)\right)
\le \displaystyle 1- \frac 1 2 n^{-1/4} .
$$
The latter inequality follows from~\eqref{1-Z-bound}, and holds for
$n$ {large enough}. Finally, for $t \in \C_n$,
\beq
\label{Z-bound}
1-|Z(t)| \ge \frac 1 2 n^{-1/4}.
\eeq

Let us now consider the integral on  the contour ${\mathcal C}_n^{(1)}$.
By Lemma~\ref{lemma-Z}, the quantity
$$
\frac { \left( 1+Z \right) ^{2}{Z}^{j+1} \left( 1+Z+{Z}^{2} \right)  
\left(1- Z \right) ^{2}}
{  1+{Z}^{2}   }
$$
 is uniformly bounded on this contour by some constant
 $c$, independant of $n$ and $t$. Moreover, 
$$
|1-Z^{j+2}|\ge 1 -|Z|^{j+2} \ge 1-|Z|\ge \frac 1 2 n^{-1/4}
$$
by~\eqref{Z-bound}.
The same bound holds for the term $1-Z^{j+4}$.
Therefore the modulus of the contribution of $\C_n^{(1)}$ in the
integral~\eqref{cauchy}  is bounded by 
\beq
\label{C1-estimate}
 4c\, 8^n n^{1/2} \,r_n^{-n}%=O(8^n\, e^{-\log^2n})
=O(8^n n^{1/2-\log n})=o(8^n/n^m)
\eeq
for any $m>0$.

Let us now study the contribution of the contour ${\mathcal C}_n^{(2)}$.
As $t$ varies along ${\mathcal C}_n^{(2)}$, the variable $z$ defined by
$t=(1+z/n)/8$ varies along the contour 
$\Ha_n$.
As $n$ goes to infinity, this contour converges to the contour
$\Ha$ shown on the left side of Figure~\ref{fig-hankel}. Let 
$z\in\Ha$. Then $z\in\Ha_n$ for $n $ large enough, 
$|z|\le |z_n|\sim \log^2n$, and,
as $n$ goes to infinity,
the following approximations hold 
with error terms independent of $z$:
\beq
\label{estimates}\left\{
\begin{array}{rcll}
\displaystyle 
Z(t)&=&\displaystyle 1-2 (-z)^{1/4} n^{-1/4} + O\left( n^{-1/2} \log n \right)\\
1-Z(t)&=&\displaystyle  2 (-z)^{1/4} n^{-1/4}\left( 1 +  O( n^{-1/4} \sqrt{\log n})\right)\\
\\
Z(t)^j& =& \displaystyle \exp(-2\lambda (-z)^{1/4}) \left(1+O(n^{-1/4}{\log
  n})\right)
&\displaystyle (\hbox{recall } j=\lfloor \lambda n ^{1/4}\rfloor)\\
\\
\displaystyle t^{-n-1} &=&\displaystyle \displaystyle 8^{n+1} e^{-z} \left( 1+ O(\log^4 n /n)\right).
\end{array}
\right.
\eeq
Observe  that, for $z \in \Ha$, the real part of $(-z)^{1/4}$
is bounded from below by a positive constant $\alpha$. 
Hence
$$
|\exp(-2\lambda (-z)^{1/4})|=  \exp(-2\lambda \Re (-z)^{1/4}) \le 
\exp(-2\lambda \alpha),
$$
so that $\exp(-2\lambda (-z)^{1/4}) $ does not approach $1$. 
This allows us to write
$$
\frac 1 {1-Z^{j+2}} = \frac 1 { 1 -\exp(-2\lambda (-z)^{1/4})}  \left(1+O(n^{-1/4}{\log
  n})\right).
$$
%so that $\exp(-2\lambda (-z)^{1/4}) $ does not approach $1$. 
Hence, uniformly in $t \in \C_n^ {(2)}$, we have
\begin{eqnarray*}
U_j(t) t^{-n-1}&=&
\frac { \left( 1+Z \right) ^{2}{Z}^{j+1} \left( 1+Z+{Z}^{2} \right)  
\left(1- Z \right) ^{2}}
{ \left( 1+{Z}^{2} \right) 
 \left( 1-{Z}^{j+2} \right)  \left( 1-{Z}^{j+4} \right) }
t^{-n-1}\\
&=& \frac{6. 8^{n+1}}{n^{1/2}} 
\frac{\sqrt{-z}e^{-z}}{\sinh^2(\lambda (-z)^{1/4})} (1+ O( n^{-1/4} \log n))
\end{eqnarray*}
with $8t=1+z/n$. 
Let us now integrate this over ${\mathcal C}_n^{(2)}$:
\begin{eqnarray*}
\displaystyle \int_{{\mathcal C}_n^{(2)}}
 U_j(t)\frac{dt}{t^{n+1}}
&=&\frac{6. 8^n}{n^{3/2}} \int_{\Ha_n}
\frac{\sqrt{-z}e^{-z}(1+ O( n^{-1/4} \log n))}{\sinh^2(\lambda (-z)^{1/4})}
dz\\
&=& \frac{6. 8^{n}}{ n^{3/2}} \left(
\int_{\Ha}
\frac{\sqrt{-z}e^{-z}}{\sinh^2(\lambda (-z)^{1/4})} dz
+o(1)\right).
\end{eqnarray*}
We now put together our estimates of the integrals on $\C_n^{(1)}$
(Eq.~\eqref{C1-estimate}) and $\C_n^{(2)}$ and obtain
$$
[t^n] U_j(t) =\frac{6. 8^{n} n^{-3/2}}{2i \pi} \left(
\int_{\Ha  }
\frac{\sqrt{-z}e^{-z}}{\sinh^2(\lambda (-z)^{1/4})}dz
+o(1)\right).
$$
Using~\eqref{GU} and the estimation $C_n\sim 4^n n
^{-3/2}/\sqrt \pi$, this gives
$$
\PP(N_n > \lambda )\rightarrow
\frac{3}{i \sqrt \pi} \int_{\Ha}
\frac{\sqrt{-z} e^{-z}}{\sinh^2(\lambda (-z)^{1/4})} dz.
$$

\begin{figure}[pbth]
  \begin{center}
    \includegraphics[scale=0.3]{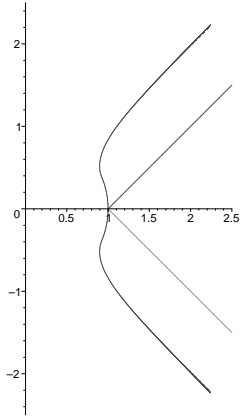}
    \caption{The contours $\Gamma$ (two half lines) and $\J$.}
    \label{figure-contours}
  \end{center}
\end{figure}

The next step in our proof of Theorem~\ref{thm-max} is to set 
$v=(-z)^{1/4}$ in the above integral. As $z$ runs on $\Ha$, the
variable $v$ runs on the contour $\J$ of Figure~\ref{figure-contours}, and the corresponding 
integral is easily seen to coincide with the integral on the contour
$\Gamma$ defined in the statement of the theorem. This gives the first
expression of $G(\lambda)$.

We now want to express $G(\lambda)$ as a real integral. We first observe
that the integration contour $\Gamma$ can be replaced by its
translated version
$$
\Gamma_0= \{-re^{-i\pi/4}, r\in (\infty,0]\} \cup \{re^{i\pi/4},
r\in [0,\infty)\}. 
$$
This parametrization of $\Gamma_0$ by $r$ splits the integral into two
real integrals, and one finds:
\begin{eqnarray*}
G(\lambda) & =&
-  \frac {12}{\sqrt \pi}\int_ 0^\infty \left( \frac 1{\sinh^2(\lambda
  r e^{i\pi/4})} + \frac 1{\sinh^2(\lambda
  r e^{-i\pi/4})} \right) r^5 e^{-r^4} dr \\
&=&
 \frac {48}{\sqrt \pi}\int_ 0^\infty  
\frac{1-\cos (\sqrt 2 \lambda r) \cosh
  (\sqrt 2 \lambda r)}{(\cosh (\sqrt 2 \lambda r) -\cos (\sqrt 2 \lambda r))^2}
r^5 e^{-r^4} dr .
\end{eqnarray*}
The expected expression of $G(\lambda)$ % the limit of the
				% distribution function 
follows, upon setting $u=\sqrt 2 \lambda r$.

\subsection{The limit law and its density}
We now want to prove that
$G(\lambda)$ is the tail distribution function of a random variable.
Since it is the limit of
non-increasing functions, it is non-increasing. Its integral
expressions show that it is a continuous, and even a differentiable
function of $\lambda$ on $(0, +\infty)$. In order to conclude, we
still need to prove that~\cite[Thm.~14.1]{billingsley}
$$
\lim_{\lambda \rightarrow \infty}G(\lambda)=0 \hbox{ and } 
\lim_{\lambda \rightarrow 0}G(\lambda)=1.
$$
In order to prove the first statement, we use the second expression of
$G(\lambda)$ given in the theorem. We note that the function
$$
u\mapsto \frac{1-\cos u \cosh   u}{(\cosh u -\cos u)^2}
$$
is well-defined, bounded  and continuous on $[0, +\infty)$. Moreover,
  as $u $ goes to infinity,
$$
\left|\frac{1-\cos u \cosh   u}{(\cosh u -\cos u)^2}\right|
=O(e^{-u}),
$$
so that the integral
$$
\int_0^\infty\left| \frac{1-\cos u \cosh
  u}{(\cosh u -\cos u)^2}\right| u^5 du
$$
is convergent. The term $1/\lambda ^6$ in the expression of
$G(\lambda)$ then implies the convergence of $G(\lambda)$ to $0$ as
$\lambda \rightarrow \infty$.

In order to study the limit of  $G(\lambda)$ as $\lambda \rightarrow
0^+$, we consider instead the first expression of $G(\lambda)$. Since
$x^2/\sinh^2(x)$ is analytic in the disk of radius $\pi$, with expansion
 $1-x^2/3+ O(x^4)$, there exists a constant $c$ such that for $|v| \le
\pi /(2\lambda)$, 
\beq
\label{int-expansion}
\left| \frac 1 {\sinh ^2(\lambda v)}-\frac 1 {\lambda^2 v^2} +\frac 1
3\right|
\le c \lambda^2 |v|^2.
\eeq
Let us write
$$
\int_ \Gamma\frac{v^5 e^{v^4}}{\sinh^2(\lambda v) } dv
=\int_ \Gamma\left( \frac{1}{\sinh^2(\lambda v)}  -\frac 1 {\lambda^2
    v^2} +\frac 1 3\right) v^5 e^{v^4}dv 
+ \int_ \Gamma \left( \frac {v^3}{\lambda ^2} - \frac {v^5} 3\right)
e^{v^4}dv .
$$
Recall the Hankel expression of the reciprocal of the Gamma function,
valid for any $s\in \cs$:
\beq
\label{gamma-inverse}
\frac 1 {\Gamma (s)} =\frac 1{2i \pi}
\int_{\Ha} (-z)^{-s} e^{-z}dz
= \frac 2 {i \pi}  \int_\Gamma v^{3-4s} e^{v^4} dv.
\eeq
Consequently,
$$
\int_ \Gamma   {v^3} e^{v^4}dv = \frac{i\pi}{2\Gamma(0)}=0 ,
\quad  \quad
\int_ \Gamma   {v^5} e^{v^4}dv =
\frac{i\pi}{2\Gamma(-1/2)}=-\frac{i\sqrt {\pi}} 4,
$$
and we can rewrite
$$
G(\lambda)= \frac{12}{i\sqrt {\pi}}
\int_ \Gamma\frac{v^5 e^{v^4}}{\sinh^2(\lambda v) } dv
=1 + \frac{12}{i\sqrt {\pi}}\int_ \Gamma\left( \frac{1}{\sinh^2(\lambda v)}  -\frac 1 {\lambda^2
    v^2} +\frac 1 3\right) v^5 e^{v^4}dv .
$$
Let us  cut the above integral into two parts,  $|v|\le
\pi/(2\lambda)$ and  $|v|>\pi/(2\lambda)$. The first part is easily
seen to tend to $0$ as $\lambda$ does, thanks
to~\eqref{int-expansion}. For the second part, we observe that for
$\lambda |v| >\pi/2$,
$$
\left| \frac 1 {\sinh ^2(\lambda v)}-\frac 1 {\lambda^2 v^2} +\frac 1
3\right|
$$
is bounded (by a constant independent of $\lambda$ and $v$), that the
integral of $v^5e^{v^4}$ on $\Gamma$ is absolutely convergent, and
that the contour $\{v\in \Gamma : |v|> \pi/(2\lambda)\}$ ``shrinks to
$\infty$'' as $\lambda \rightarrow 0$. We finally conclude that $
G(\lambda)$ tends to $1$ as $\lambda \rightarrow 0$.

Consequently, there exists a random variable $N$ having
distribution function $1-G(\lambda)$, and $N_n$ converges in law to
$N$. Since $G$ is differentiable, 
$N$ has a density with respect to the Lebesgue measure on $\rs_+$,
which is $f(\lambda)=-G'(\lambda)$. The two expressions of $G$ given
in the theorem provide the two expressions of $f$.

\subsection{The moments of $N$}
Let us first prove that for all $k\ge 0$, the tail distribution function
of $N$ satisfies
\beq
\label{G-infty}
G(\lambda)= o(\lambda^{-k}) \hbox{ as } \lambda \rightarrow \infty.
\eeq
This is easily seen to imply the existence of moments of $N$ of all
orders. In order to prove the above bound, we write
$$
G(\lambda)= \frac{24}{i\lambda \sqrt \pi} \int_\Gamma \frac
{v^4(5+4v^4)e^{v^4}}
{e^{2\lambda v}-1} dv.
$$
This is obtained from the first expression of $G( \lambda)$ using an
integration by parts. Now, for $\lambda>0$ and $v\in \Gamma$,
$$
|e^{2\lambda v}-1|\ge |e^{2\lambda v}|-1= e^{2\lambda \Re(v)}-1
\ge  e^{2\lambda }-1.
$$
From this, and from the term $e^{v^4}$ in the integral, we conclude
that there exists a constant $c$ such that
$$
G(\lambda)\le \frac{c}{ e^{2\lambda }-1}.
$$
The bound~\eqref{G-infty}
follows. This bounds also guarantees that for $k\ge 1$,
\beq
\label{moments-tail}
\EE(N^k)= k\int_0^{\infty} \lambda^{k-1} G(\lambda) d\lambda.
\eeq

\noindent
{\bf The generic case: $k\ge 3$}.
Recall the following integral representations of the Riemann zeta function: for
$\Re(s)>1$,
$$
\zeta(s) =\frac 1 {\Gamma(s)} \int_0^\infty \frac{w^{s-1}}{e^w-1}dw
=
\frac 1 {4\Gamma(s+1)} \int_0^\infty \frac{w^s}{\sinh^2(w/2)} dw
=
\frac {2^{s-1}} {\Gamma(s+1)} \int_0^\infty \frac{y^s}{\sinh^2(y)} dy.
$$
The second expression follows from the first one after an integration
by parts.

Let us now combine~\eqref{moments-tail} with the first expression of
$G( \lambda)$:
\beq
\label{ENk-first}
\EE(N^k)= \frac {12k}{i\sqrt\pi} \int_0^\infty \lambda^{k-1} d\lambda
\int_\Gamma \frac{v^5e^{v^4 }}{\sinh^2(\lambda v)} dv.
\eeq
Assume for the moment that we can exchange the order of integration
(this will be justified later). Exchange the integrals, and replace the
variable $\lambda$ by $y/v$, where $y$ is a new variable:
$$
\EE(N^k)= \frac {12k}{i\sqrt\pi}
\int_\Gamma v^{5-k}e^{v^4 }dv
 \int_{v\rs_+}\frac{y^{k-1}}{\sinh^2(y)}dy.
$$
For $k\ge 3$, the function $y\mapsto {y^{k-1}}/{\sinh^2(y)}$ is
meromorphic on $\cs$, with poles at $ik\pi$ for $k \in \zs$ and
$k\not = 0$. From this, and from the strong decay of this function as
$\Re(y) \rightarrow  \infty$, it follows that the integral on $y$ is
actually independent of the choice of $v\in \Gamma$. In particular, it
is equal to its value at $v=1$, which is
$$
 \int_0^\infty\frac{y^{k-1}}{\sinh^2(y)}dy = \frac{4 \Gamma(k)
 \zeta(k-1)}{2^k},
$$
as recalled above. The integral on $v$ is then evaluated in terms of
the Gamma function using~\eqref{gamma-inverse}, and the expected
expression of $\EE(N^k)$ follows.

It remains to justify the exchange of integrals in~
\eqref{ENk-first}. Observe that 
$$
|\sinh(y)|= |e^y -e^{-y}|/2 \ge \left( |e^y| -|e^{-y}|\right)/2
=\sinh(\Re(y)),
$$
so that for $v \in \Gamma$,
$$
\frac 1{|\sinh^2 (\lambda v)|} \le \frac 1{\sinh^2 (\lambda )}.
$$
Moreover, the integral of $v^5e^{v^4 }$ along $\Gamma$ is absolutely
convergent, and so is the integral of $\lambda^{k-1}/\sinh^2 (\lambda
)$
over $R_+$. It follows that the integral~\eqref{ENk-first}, once
converted into two real integrals,  is
absolutely convergent, so that the  integrals can be exchanged.
%in~\eqref{ENk-first}. 

\medskip
\noindent
{\bf The case $k=1$.} We cannot apply exactly the same procedure as
above, because the integral of  $1/\sinh^2 (\lambda)$ over $\rs_+$ is
divergent.
However, in view of~\eqref{gamma-inverse}, we can write
$$
G(\lambda) = \frac {12}{i\sqrt \pi} \int_\Gamma v^5e^{v^4} 
\left( \frac 1 {\sinh^2(\lambda v)} - \frac 1 {\lambda^2 v^2}\right).
$$
Also, replacing $\Gamma$ by $\Gamma_0$ in the latter integral does not
change its value.
The technique is then the same as above:
\begin{eqnarray}
\EE(N)&= &\frac {12}{i\sqrt\pi} \int_0^\infty d\lambda
\int_{\Gamma_0} {v^5e^{v^4 }}
\left( \frac 1 {\sinh^2(\lambda v)} - \frac 1 {\lambda^2 v^2}\right)
dv
\label{EN-first}\\
&=& \frac {12}{i\sqrt\pi}
\int_{\Gamma_0} v^{4}e^{v^4 }dv
 \int_{v\rs_+}\left( \frac 1 {\sinh^2(y)} - \frac 1
     {y^2}\right)dy \nonumber
\end{eqnarray}
(assuming we can change the order of integration). Again, the integral
on $y$ is independent of $v$, and equal to
$$
\int_0^\infty \left( \frac 1 {\sinh^2(y)} - \frac 1
     {y^2}\right)dy = \left[ \frac 1 y - \frac 2
       {e^{2y}-1}\right]_0^\infty
=-1.
$$
Using again~\eqref{gamma-inverse} to evaluate the integral on $v$, one
finds
$$
\EE(N)= -\frac {6\sqrt \pi}{\Gamma(-1/4)}=
\frac{3\sqrt\pi}{2\Gamma(3/4)}.
$$
In order to justify the exchange of integrals in~\eqref{EN-first}, we
wish to prove that~\eqref{EN-first}
is absolutely convergent. In order to do so, we split the integral
over $\Gamma_0$ into two real integrals, corresponding respectively to
$v=re^{i\pi/4}$ and $v=re^{-i\pi/4}$. We are thus led to prove that
$$
\int_0^\infty d\lambda
\int_0^\infty {r^5e^{-r^4 }}
 \left|\frac 1 {\sinh^2(\lambda re^{i\pi/4})} - \frac 1 {i\lambda^2 r^2}\right|
dr
$$
is absolutely convergent (and a similar result when $i$ is
replaced by $-i$). But we can exchange the order of
integration in this integral of \emm positive, functions. Doing so,
and setting $\lambda= y/r$ as above, proves that this integral is
finite.

\smallskip
\noindent{\bf The case $k=2$.}
Let us start from another expression of $G(\lambda)$, obtained by
writing $v= w/\lambda$:
$$
G(\lambda)= \frac{12}{i\sqrt \pi \lambda ^6 }\int_{\lambda \Gamma}
\frac{w^5}{\sinh^2(w)} e^{w^4/\lambda^4} dw
= \frac{12}{i\sqrt \pi \lambda ^6 }\int_{ \Gamma}
\frac{w^5}{\sinh^2(w)} e^{w^4/\lambda^4} dw.
$$
The second expression follows from the analyticity properties of the
integrand. Now, take $\epsilon >0$, and let us evaluate
\begin{eqnarray}
\int_\epsilon ^\infty \lambda G(\lambda)d\lambda
&= &\frac{12}{i\sqrt \pi}\int_\epsilon ^\infty \frac 1 {\lambda^5} d\lambda
\int_{ \Gamma}
\frac{w^5}{\sinh^2(w)} e^{w^4/\lambda^4} dw
\label{N2-first}\\
&= &\frac{12}{i\sqrt \pi}\int_{ \Gamma}
\frac{w^5}{\sinh^2(w)}dw
\int_\epsilon ^\infty \frac {e^{w^4/\lambda^4}} {\lambda^5}d\lambda
\nonumber\\
&= &\frac{12}{i\sqrt \pi}\int_{ \Gamma}
\frac{w^5}{\sinh^2(w)}
\left[ -\frac{e^{w^4/\lambda^4}}{4w^4} \right]_\epsilon ^\infty dw\nonumber\\
&= &\frac{3}{i\sqrt \pi}\int_{ \Gamma} \frac{w}{\sinh^2(w)}
\left( e^{w^4/\epsilon^4}-1\right) dw.\nonumber
\end{eqnarray}
The absolute convergence of integrals that legitimates  the exchange
of integrals in~\eqref{N2-first} is, this time, obvious (thanks to the
fact that $\lambda >\epsilon$). 
Now, the analyticity of the function $w\mapsto {w}/{\sinh^2(w)}$ for
$\Re(w)>0$, and its strong decay as $\Re(w)\rightarrow \infty$, imply
that
$$
\int_{ \Gamma} \frac{w}{\sinh^2(w)}dw=0.
$$
Hence
$$
\begin{array}{llllllll}
\displaystyle \int_\epsilon ^\infty \lambda G(\lambda)d\lambda
&=&
\displaystyle \frac{3}{i\sqrt \pi}\int_{ \Gamma} \frac{w
  e^{w^4/\epsilon^4}}{\sinh^2(w)}dw
&=&\displaystyle \frac{3\epsilon^2}{i\sqrt \pi}\int_{ \Gamma} \frac{v
  e^{v^4}}{\sinh^2(\epsilon v)}dv\\
&=&\displaystyle  \frac{3}{i\sqrt \pi}\int_{ \Gamma} \frac{
  e^{v^4}}{v}dv + o(1)
&=&\displaystyle \frac{3\sqrt\pi} 2 & \hbox{ (by \eqref{gamma-inverse})} .
\end{array}
$$
Now, observe that
$$
2\int_\epsilon ^\infty \lambda G(\lambda)d\lambda = \EE(N^2 \1_{N>
  \epsilon}) - \epsilon^2 G(\epsilon).
$$
 The announced expression of the second moment of $N$ follows.

%%%%%%%%%%%%%%%%%%%%%%%%%%%%%%%%%%%%%%%%%%%%%%
\subsection{Convergence of the moments of $N_n$}
\label{section-moments-Nn}
%%%%%%%%%%%%%%%%%%%%%%%%%%%%%%%%%%%%%%%%%%%%%
In this section, we prove that the moments of $N_n=M_n/n^{1/4}$
converge to the corresponding moments of $N$. In order to do so, we first
express $\EE(M_n^k)$ as the coefficient of $t^n$ in a certain
series. Then, we apply the general consequences of the analysis of
singularities: if this series is regular enough (with a precise
meaning of regular), one can derive the asymptotic behaviour of its
coefficients from the singular behaviour of the series near its
dominant singularities~\cite{flajolet-odlyzko}.

Recall that the series $U_j$, given by~\eqref{Uj-expr}, counts the trees that
contain at least one label larger than $j$. Hence $U_{j-1}-U_j$ counts
the trees having maximal label $j$. Also, note that
\beq
 \label{UV}
U_j= V(Z^j)-V(Z^{j+2}),
\eeq
where
$$
V(x)= \frac{xZ(1+Z)(1-Z^3)}{(1+Z^2)(1-xZ^2)}.
$$
Consequently, for $k\ge 1$,
\beq
\label{EMnk}
\EE(M_n^k)= \frac 1 {2^n C_n} \sum_{j\ge 1} j^k [t^n] (U_{j-1}-U_j)
= \frac 1 {2^n C_n}   [t^n]\sum_{j\ge 0} \left((j+1)^k-j^k\right)U_j.
 \eeq
For $k=1$, this gives
$$
 {2^n C_n}\EE(M_n)= [t^n]  \sum_{j\ge 0} \left(V(Z^j)-V(Z^{j+2})\right)= 
 [t^n]\left( V(1)+V(Z)\right)
= [t^n] \frac{Z(1+2Z+2Z^2)}{1+Z^2}.
$$
By Lemma~\ref{lemma-Z}, the latter series is analytic in
$\cs\setminus[1/8, \infty)$.  
The generic consequences of the
  analysis of singularities apply: one can derive the asymptotic
  behaviour of the coefficients from the singular behaviour of the
  series~\cite{flajolet-odlyzko}. Given that, when $t\rightarrow 1/8$,
$$
\frac{Z(1+2Z+2Z^2)}{1+Z^2}
= \frac 5 6 -6 (1-8t)^{1/4} + O(\sqrt{1-8t}),
$$
the behaviour of the $n$th coefficient of this series is
$$
 [t^n] \frac{Z(1+2Z+2Z^2)}{1+Z^2} = -6 \frac {8^n
 n^{-5/4}}{\Gamma(-1/4)} (1+o(1))
=
\frac 3 2  \frac {8^n
 n^{-5/4}}{\Gamma(3/4)}(1+o(1)).
$$
It remains to divide by $2^n C_n \sim 8^n n^{-3/2}/ \sqrt \pi$ to
conclude that
$$
\EE(M_n n^{-1/4}) \rightarrow \frac{3 \sqrt \pi}{2\Gamma (3/4)},
$$
which is also the first moment of $N$.

Now, by combining the expression~\eqref{UV} of $U_j$ and~\eqref{EMnk},
one obtains, for $k\ge 2$,
\beq
\label{EMk}
2^nC_n \EE(M_n^k)=[t^n]\left( V(1)+(2^k-1)V(Z) +\sum_{j\ge 2}
\left((j+1)^k-j^k-(j-1)^k+(j-2)^k\right) V(Z^j)\right).
\eeq
Observe that $(j+1)^k-j^k-(j-1)^k+(j-2)^k$ is a polynomial in $j$ of
degree $k-2$ and leading coefficient $2k(k-1)$. Let
$$
A_\ell(t) = \sum_{j\ge -1}(j+2)^\ell V(Z^j).
$$
We are going to prove
that, for $\ell\in\ns$,
\beq
\label{a-asympt}
\displaystyle a_n(\ell):= [t^n] A_\ell(t)=\left\{
\begin{array}{lll}
\displaystyle \frac 3 4 \frac {8^n} n  &
\hbox{ if }  \ell=0,\\
\displaystyle \frac{3. 8^n \ell ! \zeta(\ell+1) n^{\ell/4 -1}}
{2^\ell \Gamma(\ell/4)} & \hbox{ if } \ell \ge 1.
\end{array}\right.
\eeq
Assume for the moment this is proved, and let us conclude about the
limiting moments of $N_n=M_n n^{-1/4}$. First, we observe that for $j \ge 0$,
$V(Z^j)$ has (only) a fourth root singularity, so
that the coefficient of $t^n$ in $V(Z^j)$ grows like $8^n n^{-5/4}$,
up to a multiplicative constant. This observation, combined
with~\eqref{EMk} and the
above asymptotics of $ a_n(\ell)$, implies that the dominant term in the
asymptotic  behaviour  of $2^nC_n\EE(M_n^k)$ is that of $2k(k-1)
a_n(k-2)$. After normalizing by $2^n C_n n^{k/4}$, this gives
$$
\EE(M_n^kn^{-k/4}) \rightarrow
\left\{
\begin{array}{lll}
\displaystyle 3\sqrt \pi&
\hbox{ if }  k=2,\\
\displaystyle \frac{24 \sqrt \pi k!\zeta(k-1)}
{2^k \Gamma((k-2)/4)} & \hbox{ if }k \ge 3.
\end{array}\right.
$$
These limiting moments are exactly those of $N$.

\medskip
It remains to study the asymptotic behaviour of the numbers
$a_n(\ell)$ (for $\ell$ fixed, and $n$ going to infinity). We have:
$$
\begin{array}{llll}
 \displaystyle A_\ell(t)&=&\displaystyle \frac{(1+Z)(1-Z^3)}{Z(1+Z^2)} 
\sum_{j\ge 1}j^\ell   \frac {Z^{j}}{1-Z^{j}}\\
&=&\displaystyle \frac{(1+Z)(1-Z^3)}{Z(1+Z^2)} \sum_{j\ge 1, m\ge 1}j^\ell Z^{jm}\\
&=&\displaystyle \frac{(1+Z)(1-Z^3)}{Z(1+Z^2)} \sum_{N\ge 1} Z^N \sigma_\ell(N)
\end{array}
$$
where
$$ 
\sigma_\ell(N)=\sum_{j|N} j^\ell.
$$
The function
$$
D_ \ell(z)=  \sum_{N\ge 1} z^N \sigma_\ell(N)
$$
is easily seen to have radius of convergence $1$. Moreover, as $z$
tends to $1$ in such a way $|\arg(1-z) |<\phi < \pi/2$,
$$
D_ \ell(z)\sim \left\{
\begin{array}{lll}
\displaystyle \frac 1{1-z} \log \left(\frac 1 {1-z}\right) 
& \hbox{ if } \ell=0,\\
\displaystyle \frac{\ell !\zeta(\ell+1)}{(1-z)^{\ell+1}} &
    \hbox { if } \ell \ge 1
\end{array}
\right.
$$
(this can be obtained using a Mellin
transform~\cite{flajolet-binary-trees,flajolet-gourdon-dumas}).  
The above expression of $A_\ell(t)$, combined with Lemma~\ref{lemma-Z}
and these properties of $D_\ell(z)$, shows that $A_\ell(t)$ is analytic
in the domain $\D=\cs\setminus[1/8, \infty)$.
%, and has no other singularity that  $1/8$. 
Moreover, since $|\arg(1-Z)| \le \pi /4 +o(1)$ as $t
  \rightarrow 1/8$ in $\D$, we can use the above estimates of
  $D_\ell(z)$. This gives
$$
A_\ell(t) \sim
\left\{
\begin{array}{lll}
\displaystyle -3 \log 2 + \frac 3 4  \log \left(\frac 1 {1-8t}\right) 
& \hbox{ if } \ell=0,\\
\displaystyle \frac{3 \ell !\zeta(\ell+1)}{2^\ell (1-8t)^{\ell/4}} &
    \hbox { if } \ell \ge 1.
\end{array}
\right.
$$
 The generic results derived from the
  analysis of singularities apply, and give the asymptotic
 behaviour~\eqref{a-asympt} of the numbers $a_n(\ell)$. 
This concludes the proof of Theorem~\ref{thm-max}.

 \subsection{The supremum of the support of the  ise}
\label{section-supremum}
Let us finally prove Proposition~\ref{coro-support}. 
The following  argument 
%is due to Jean-Fran\c cois Marckert. It 
requires a
 detour via discrete snakes and Brownian snakes. We refer
 to~\cite{legall,marckert-mokka-snake,jf-janson} for definitions and
 notation\footnote{We warn the reader that normalizations change from
 one paper to another.}. In
 particular, we  use the following integral representation of the
 random measure $\mu_{\hbox{\sc ise}}$: %~\cite{legall,marckert-mokka-snake}: 
 for any  continuous  bounded function $g$ on $\rs$,
\beq
\label{ise-rep}
\int_\rs g(y) d\mu_{\hbox{\sc \small ise}}(y)= \int_0^1 g( r(t))dt
\eeq
 where $r(.)$ is a random process, continuous on $ [0,1]$, called the
 head of the Brownian  snake. In other words, $\mu_{\hbox{\sc \small
     ise}}$ is the \emm occupation measure, of the process $r$. (Again, the
 definition of $r$ varies from 
 one paper to the other. The above formula fixes our normalization of $r$.)

 The random variable $N_n= M_n n^{-1/4}$ coincides with
 $\max(r_n)$, where $r_n$ is the (normalized) head of the discrete
 snake associated with our tree family. The random process $\sqrt 2 r_n$
 converges weakly to $r$, the head of the Brownian
 snake~\cite{marckert-mokka-snake}.  Since $\max$ 
 is a continuous functional on $\C[0,1]$, this implies that $
\sqrt 2  N_n= \sqrt 2 \max (r_n)$  converges in distribution to $\max
 (r)$. Thus $\max( r)$ has  density $f(\lambda/\sqrt 2)/\sqrt
 2$, where $f$ is defined in Theorem~\ref{thm-max}. 
% In particular,  $\max(r)$ does not charge points.

\begin{figure}[htb]
\begin{center}
\input{fepsilon.pstex_t}
\end{center}
\caption{The functions $f_{\lambda,\epsilon}$, $g_{\lambda,\epsilon}$ and  $h_{\lambda,\epsilon}$.}
\label{fig-epsilon}
\end{figure}
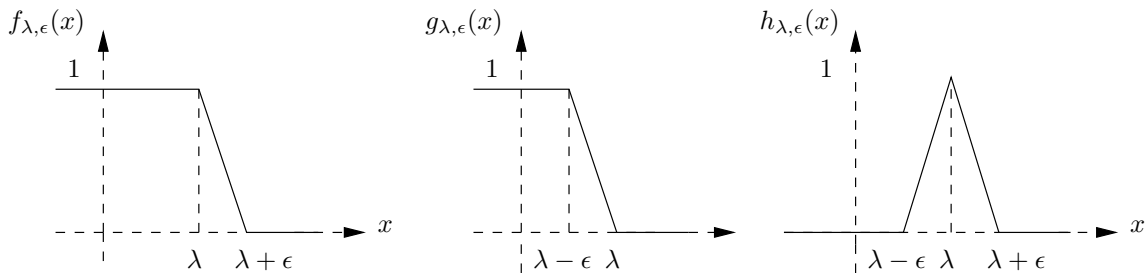

It remains to prove that $\max (r)$ is equal (in distribution) to 
$N_{\hbox{\sc \small ise}}$, the supremum of the support of the  ISE.
Let $\lambda\in \rs$ and $\epsilon >0$. Let
$f_{\lambda, \epsilon} $ % and $g_{\lambda,\epsilon}$ 
be the function plotted on the left-hand side of
Figure~\ref{fig-epsilon}. We have  
$$
N_{\hbox{\sc \small ise}} \le \lambda \Longleftrightarrow 
\mu_{\hbox{\sc \small ise}}(-\infty, \lambda]=1 \Longleftrightarrow\int_\rs
f_{\lambda,1} (y) d\mu_{\hbox{\sc \small ise}}(y)=1.
$$
Thanks to~\eqref{ise-rep}, this gives
$$
N_{\hbox{\sc \small ise}} \le \lambda \Longleftrightarrow 
 \int_0^1 f_{\lambda,1} ( r(t)) dt =1
$$
Taking probabilities yields to
$$
\PP( N_{\hbox{\sc \small ise}} \le \lambda) = 
 \PP\left( 
\int_0^1 f_{\lambda,1} (r(t)) dt =1
\right)
=\PP(\max(r) \le \lambda),
$$
since $r$ is almost surely continuous.

%%%%%%%%%%%%%%%%%%%%%%%%%%%%%%%%%%%%%%%%%%%%%%%%%
\section{A local limit law}
\label{section-limit-local}
%%%%%%%%%%%%%%%%%%%%%%%%%%%%%%%%%%%%%%%%%%%%%%%%%%

For $j \in \zs$, let $X_n(j)$ denote the random variable equal to the
number of nodes having label $j$ in a random tree of  $\mathcal T
_{0,n}$. This quantity is related to the series $S_j(t,u)$ studied in
Proposition~\ref{propo-snake-gf}. In particular,
$$
%P\left(X_n(j)=m\right)=  \frac{[u^mt^n] S_j}{2^n C_n}.
\EE\left(e^{aX_n(j)}\right)=  \frac{[t^n] S_j(t,e^a)}{2^n C_n}.
$$
Also, observe that
$$
X_n(j)=0 \Longleftrightarrow M_n <j,
$$
where $M_n$ is the largest label, studied in the previous section.
Let us define a normalized version of $X_n(j)$ by
$$
Y_n(j)= \frac{X_n(j)}{n^{3/4}}.
$$
Let $\lambda \in \rs$. The aim of this section is to prove that
$Y_n(\lfloor \lambda n^ {1/4}\rfloor)$ converges in distribution, as
$n$ goes to infinity,  to a random variable $Y(\lambda)$ that we
describe by its Laplace transform.  This is achieved in
Theorem~\ref{thm-conv} below, but we first want to present two
consequences of this theorem, which have a simpler formulation. The
first consequence deals with the case $\lambda=0$. Recall
that, up to a normalization by $n^{3/4}$, the random variable
$Y_n(0)$ gives the number of nodes labelled $0$ in a tree rooted at $0$.
\begin{Propo}[{\bf The number of nodes labelled $0$}]
\label{limit-0}
As $n$ goes to infinity, the random variable $3Y_n(0)/\sqrt 2$
converges in distribution to $T^{-1/2}$, where 
 $T$ follows a unilateral stable law of parameter $2/3$.
The convergence of the moments holds as well: for $k\ge 0$,
$$
\EE\left( Y_n(0)^k\right) \rightarrow
\left(\frac{\sqrt 2} 3\right)^k \frac{\Gamma(1+3k/4)}{\Gamma(1+k/2)}=
\left(\frac{\sqrt 2} 3\right)^k \EE(T^{-k/2}). 
$$
\end{Propo}
\noindent  
This proposition will be proved in Section~\ref{moments-proof}.
I am indebted to Alain Rouault, who recognized that the above moments
were related to $T$.
 Recall that $T$ is given by its Laplace transform:
$$
\EE(e^{-aT})= e^{-a^{2/3}} \hbox{ for } a\ge 0.
$$

\smallskip

\begin{figure}[pbth]
  \begin{center}
    \includegraphics[height=5cm,width=5cm]{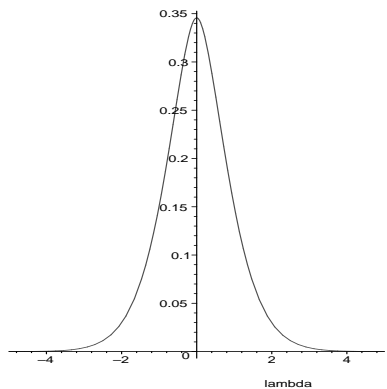}
    \caption{The average number of nodes labelled $\lfloor
    \lambda\sqrt n\rfloor $ in
    a tree of size $n$, when $n\rightarrow \infty$.}
    \label{fig-profil}
  \end{center}
\end{figure}

The second consequence of Theorem~\ref{thm-conv} is an explicit
expansion in $\lambda$ of the limiting first moment of $Y_n(j)$.
\begin{Propo}[{\bf The first moment}]
\label{propo-first}
 Let $\lambda \in \rs$.  Denote $j
=\lfloor \lambda n^{1/4}\rfloor$. Then, as $n$ goes to infinity,
$$
\EE\left( Y_n( j) \right)
\rightarrow
\frac 1 {\sqrt\pi}\sum_{m\ge 0} \frac{(-2|\lambda|)^m}{m!} \cos %\left(
\frac{(m+1)\pi}4%\right)
\Gamma\left(\frac{m+3}4\right).
$$
This function of $\lambda$ is plotted on Figure~{\em\ref{fig-profil}}.
\end{Propo}
\noindent
Similar, but more and more complicated expressions may be written for
the next moments of $Y(\lambda)$. This proposition will be proved in
Section~\ref{first-proof}.
Let us, finally, state our main theorem, from which the two above
propositions derive.

\begin{Theorem} [{\bf A local limit law}]
\label{thm-conv}
Let $\lambda\ge 0$.
The sequence $Y_n(\lfloor \lambda n^{1/4}\rfloor)$ converges in
distribution to a non-negative 
random variable $Y(\lambda)$ whose Laplace transform is given, for
$|a| < 4/\sqrt 3$,  by
$$
\EE\left( e^{aY(\lambda)}\right)=L(\lambda,a)
$$
where
$$
L(\lambda, a)=
1+ \frac {48}{i\sqrt \pi }\int_{\Gamma}
\frac{A(a/v^3)e^{-2\lambda v}}
{(1+A(a/v^3)e^{-2\lambda v})^2}v^{5} e^{v^4} dv,
$$
$A(x)\equiv A$ is the unique solution of
\beq
\label{A-eq}
A=\frac x{24} \frac{(1+A)^3}{1-A}
\eeq
satisfying $A(0)=0$,
and the integral is taken over
$$
\Gamma= \{1-te^{-i\pi/4}, t\in (\infty,0]\} \cup \{1+te^{-i\pi/4},
t\in [0,\infty)\}. 
$$
More precisely, the Laplace transform of $Y_n(\lfloor \lambda n^{1/4}\rfloor)$
converges pointwise to $L(\lambda, \cdot)$ on the interval  $(-
4/\sqrt 3, 4/\sqrt 3)$. The convergence of moments holds as well.
\end{Theorem}
It is believed (or known?) that the random measure $\mu_{{\hbox{\sc
      \small ise}}}$ is almost surely absolutely continuous with
      respect to the Lebesgue measure on
      $\rs$. Eq.~\eqref{limit-density} leads us to the following conjecture.
\begin{Conj}[{\bf The density of the ISE}]
\label{conj-density}
 There exists a random continuous process
$f_{{\hbox{\sc \small ise}}}(\lambda)$, defined for $\lambda\in \rs$,
 such that 
$\mu_{{\hbox{\sc \small ise}}}=f_{{\hbox{\sc \small ise}}}\hbox{Leb}$, where $\hbox{Leb}$ denotes the Lebesgue measure
 on  $\rs$. Moreover,  $f_{{\hbox{\sc \small ise}}}(\lambda)$ satisfies
$$
f_{{\hbox{\sc \small ise}}}(\lambda)\deq \frac 1 {\sqrt 2} \,
Y\left(\frac{|\lambda|}{\sqrt 2}\right), 
$$
where the law of  $Y(\lambda)$ is given in
Theorem~{\em\ref{thm-conv}.}
\end{Conj}
\noindent{\bf Comments}\\
{\bf 1.} The limit random variable $Y(\lambda)$ equals $0$
with a \emm positive probability,  as soon as $\lambda>0$. Indeed, by
the  portmanteau Theorem~\cite[Thm.~11.1.1]{dudley}, 
$$
\PP(Y(\lambda)=0) \ge \limsup \PP(Y_n(\lfloor \lambda n^{1/4}\rfloor
)=0) 
= \limsup \PP(M_n< \lfloor \lambda n^{1/4}\rfloor).
$$
But, by Theorem~\ref{thm-max}, 
$$
\PP(M_n< \lfloor \lambda n^{1/4}\rfloor) \rightarrow 1-G(\lambda) >0.
$$
{\bf 2.} Let us add a few words on the series $A(x)$ defined
by~\eqref{A-eq}, in order to convince ourselves that the integral
giving $L(\lambda, a)$ is well-defined.
 Clearly, the expansion of $A(x)$ at $x=0$ has
non-negative coefficients. Looking at the discriminant of the equation that
defines $A$ shows that $A$ has radius of convergence at least $4/\sqrt
3$. Moreover,  it is
easy to prove that $A(4/\sqrt 3)=2-\sqrt 3=0.26\ldots$ Consequently, $|A(x)|$ is
bounded by $2-\sqrt 3$ for $|x|\le 4/\sqrt 3$.\\
 Since $|v|\ge 1$ for $v \in
\Gamma$, the modulus of $A(a/v^3)$ is bounded from above by $2-\sqrt
3$. Moreover, $\Re(v)\ge 1$, so that $|e^{-2\lambda v}|\le
e^{-2\lambda}<1$. Hence 
$$
\frac{A(a/v^3)e^{-2\lambda v}}
{(1+A(a/v^3)e^{-2\lambda v})^2}
$$
is uniformly bounded on $\Gamma$, and $L(\lambda, a)$ is well-defined.

Note that the series $A(x)$ admits the following closed form expression:
\beq
\label{A-explicit}
A(x)= \frac 2 {1+\frac 2 {\sqrt 3} \cos(\frac{\arccos(-x\sqrt 3 /4)} 3)}-1.
\eeq
This can be checked by proving that this expression
satisfies~\eqref{A-eq} and the initial condition $A(0)=0$. 

%%%%%%%%%%%%%%%%%%%%%%%%%%%%%%%%%%%%%%%%%%%%%%%%%%%%%%%%%%%%%%%%%%
\subsection{Proof of Theorem~\ref{thm-conv}}
\label{local-proof}
Let $\lambda \ge 0$ and $j=\lfloor \lambda n^{1/4}\rfloor$. Let us
first express the Laplace transform of $Y_n(j)$  in terms of the 
\gfs\ $S_j(t,u)$ of Proposition~\ref{propo-snake-gf}:
\beq
\label{moment-series}
\EE\left( e^{aY_n(j)}\right)= \EE\left( e^{an^{-3/4}X_n(j)}\right)
=\frac{[t^n] S_j(t, e^{an^{-3/4}})}{2^nC_n}.
\eeq

Again, we will evaluate this Laplace transform thanks to the analysis
of singularities~\cite{flajolet-odlyzko}.
% The expression of $S_j$ is given in Proposition~\ref{propo-snake-gf}. 
We wish to use again the integration contour $\C_n$ of
Figure~\ref{fig-contour}. This requires to prove that $S_j(t,u)$ is 
analytic in a neigborhood of  this contour (for $n$ large and $u=
e^{an^{-3/4}}$). This is guaranteed by the following lemma. This
lemma naturally includes some properties of the series $\mu$ involved
in the product form~\eqref{product-formSj} of $S_j$.
We denote by  $\I_n$ the part of the complex plane
enclosed by $\C_n$ (including $\C_n$ itself).
\begin{Lemma}[{\bf Analytic properties of $\mu$ and $S_j$}]
\label{lemma-Sj}
Let $a$ be a real number such that $|a|<4/\sqrt 3$.
Then there exists $\epsilon >0$ such that for $n$
large enough, the series $\mu(t,u_n)$, with $u_n=e^{an^{-3/4}}$, is
analytic in the domain
$$
\E_n =\{t : \left|t- 1/ 8 \right| >  1/ ((8+\epsilon)n)\}\setminus[ 1/ 8  , +\infty) .
$$
In particular, $\mu(t,u_n)$ is analytic in a neighborhood of  $\I_n$. 
Its modulus in  $\I_n$ is smaller than $\alpha$, for some $\alpha<1$
independent of $a$ and $n$.
The series $S_j(t,u_n)$ is also analytic in a neighborhood of  $\I_n$. 
\end{Lemma}
\noindent
{\bf Proof.}
The lemma is clear if $a=0$: in this case, $u_n=1$, the series
$\mu(t,u_n)$ vanishes, and the series $S_j$ reduces to the size \gf\
of labelled trees, namely $T$, which is analytic in
$\cs\setminus[1/8,\infty)$.
We now assume that $a\not = 0$ and $|a|< 4/\sqrt 3$. This guarantees
that $A(a)$ is well-defined, where the series $A$ is defined in
Theorem~\ref{thm-conv}.

Let us first study the singularities of the series $\bar \mu\equiv
\bar \mu(z,u)$ defined as the unique \fps \ in $z$ satisfying
$$
\bar \mu= (u-1) \frac{(1+z^2)(1+\bar \mu z)(1+\bar \mu z^2)(1+\bar \mu z^3)}
{(1+z)(1+z+z^2)(1-z)^3(1-\bar \mu z^2)}.
$$
Note that $\bar \mu$ has polynomial coefficients in $u$, and vanishes
when $u=1$.   
Assume that $u$ is a fixed real
number close to, but different from, $1$. Recall that, as all
algebraic \fps , $\bar\mu(t,u)$ has a positive radius of convergence.
Let us perform a classical analysis to detect its possible
singularities. These singularities are found in the
union of two sets $\S_1$ and $\S_2$:

$\bullet$   $\S_1$ is the set of  non-zero roots of the dominant
coefficient of the 
equation defining $\bar \mu$. That is,  $\S_1=\{\pm i\}$ ,

$\bullet$   $\S_2$ is the set of the roots of the discriminant  of the
equation defining $\bar \mu$. For $u=1+x$ and $x$ small, these roots
are found to be
$$
z=\pm 1,\quad  z=-1+O(x),\quad z=e^{\pm 2i\pi/3} +O(x),\quad
z=1+ \omega 12 ^{1/6}x^{1/3} + O(|x|^{2/3}),
$$
where $\omega$ satisfies $\omega ^6=1$. (The term $\omega$ allows us
to write loosely 
$x^{1/3}$ without saying which determination of the cubic root we take.)\\
Observe that the moduli of all these ``candidates for singularities''
go to $1$ as $x$ goes to $0$. 

Now the series $\mu=\mu(t,u)$ involved
in the expression~\eqref{product-formSj} of $S_j(t,u)$ satisfies
$$
\mu(t,u)= \bar \mu(Z(t), u)
$$
where $Z(t)$ is defined by~\eqref{eqZt}. In other words, we could have
defined the series $\bar \mu$ by
$$
\bar \mu(z,u) = \mu\left(\frac{z(1+z^2)}{(1+z)^4}, u\right).
$$
Recall that $Z$ is analytic
in the domain
$\D=\cs\setminus [1/8, \infty)$.  Take $u=u_n=
e^{an^{-3/4}}=1+x$, with $x=an^{-3/4}(1+o(1))$. By Lemma~\ref{lemma-Z}, for $n$
large, the only values of $\S_1 \cup \S_2$ that may be reached by
$Z(t)$, for $t \in\D $, are of the  form
$$
z=1+ \omega 12^{1/6}a^{1/3}n^{-1/4} + O(n^{-1/2}).
$$
In view of~\eqref{eqZt}, these values of $Z(t)$ are reached for
$$
t= \frac 1 8 - \frac{\omega ^4 (12)^{2/3}}{128} \frac{a^{4/3} } n + O(n^{-5/4}).
$$
Since $|a|< 4/\sqrt 3$, these values of $t$ are at distance less than
$1/((8+\epsilon)n)$ of $1/8$, for some $\epsilon >0$, and hence
outside the domain $\E_n$. Consequently, 
$\mu(t,u_n)$ is analytic inside $\E_n$.

\medskip
We now want to bound $\mu(t,u_n)$ inside $\I_n$. Let $t_n\in \I_n$ be such
that
$$
|\mu(t_n,u_n)|= \max _{t\in \I_n} |\mu(t,u_n)|.
$$
In particular, $|\mu(t_n,u_n)| \ge |\mu((1-1/n)/8,u_n)|$. In order to
evaluate the latter quantity, note that $Z((1-1/n)/8)=1-2n^{-1/4}+O(n
^{-1/2})$. Thanks to the closed form expression of $\mu$ given in
Proposition~\ref{propo-mu-closed-form}, and to the
expression~\eqref{A-explicit} of the series $A$, we see that
$\mu((1-1/n)/8,u_n)\rightarrow A(a)$. Since $a\not = 0$, $A(a) \not
=0$, and for $n$ large enough,
\beq
\label{mu-lower}
|\mu(t_n,u_n)| \ge |\mu((1-1/n)/8,u_n)|=|A(a)|+o(1)>0.
\eeq
Recall that all the sets $\I_n$ are included in a ball of finite
radius centered at the origin. Let $\alpha$ be an accumulation point
of the sequence $t_n$. Then  $|\alpha|\le 1/8$. 

Assume first that  $\alpha \not = 1/8$. Then there exists $N$ such
that $\alpha$ is in $\E_n$  for all $n\ge N$, that is, in the
analyticity domain of $\mu(\cdot,u_n)$. Let $t_{n_1}, t_{n_2}, \ldots$
converge to $\alpha$.  By continuity of $\mu$ in
$t$ and $u$, we have
$$
\mu(t_{n_i},u_{n_i}) \rightarrow  \mu(\alpha,1)=0.
$$
This contradicts~\eqref{mu-lower}. Hence the only accumulation point
of $t_n$ is $1/8$, and $t_n$ converges to $1/8$.
Let us thus write
$$
t_n= \frac 1 8 \left( 1 -\frac {x_n} n \right).
$$
We have $x_n=o(n)$, but also $|x_n|>1$ since $t_n$ belongs to $\I_n$.
We wish to estimate $\mu(t_n, u_n)$.
%, using~\eqref{eqMZ}. 
From the singular behaviour of $Z$ (Lemma~\ref{lemma-Z}), we derive
$$
Z(t_n)= 1-2  \left( \frac {x_n} n \right)^{1/4} + O\left(  \left(
\frac {x_n} n \right)^{1/2}\right).
$$
Moreover,
$$
u_n-1=an^{-3/4} \left( 1 + O(n^{-3/4})\right).
$$
This gives
$$
\frac{u_n-1}{(1-Z)^3}= \frac{a}{8x_n^{3/4}}\left(1+ O\left(  \left(
\frac {x_n} n \right)^{1/4}\right)\right).
$$
 If the sequence $x_n$ was
unbounded, then there would exist a subsequence $x_{n_i}$ converging to
infinity. Then $({u_{n_i}-1})/{(1-Z)^3}$ would tend to $0$.   The closed
form expression of $\mu$ given in 
Proposition~\ref{propo-mu-closed-form} implies that $\mu(t_{n_i},u_{n_i})$
would tend to $0$, contradicting~\eqref{mu-lower}. Hence the sequence
$x_n$ is bounded, and one derives from the explicit expressions of
$\mu$ and $A$ that
$$
\mu(t_n,u_n)= A(ax_n^{-3/4}) + o(1).
$$
Since $A$ is bounded by  $2-\sqrt 3$ inside its disk of convergence,
$|\mu(t_n,u_n)|$ is certainly smaller than some $\alpha$ for
$\alpha<1$ and $n$ large enough.
This concludes the proof of the second statement of Lemma~\ref{lemma-Sj}.

\medskip
By continuity of $\mu(t,u_n)$, this function of $t$ is still bounded by $1$
(in modulus) is a neighborhood of $\I_n$. Recall also that the modulus
of $Z(t)$ never reaches $1$ for $t \in \cs \setminus [1/8,
  \infty)$. The form~\eqref{product-formSj} then implies that
  $S_j(t,u_n)$ is an analytic function of $t$ in a neigbourhood of $\I_n$.
\cqfd

Let us now go back to the expression~\eqref{moment-series} of the
Laplace transform of $Y_n(j)$. Thanks to the lemma we have just
proved, we can use the Cauchy formula to extract the coefficient of
$t^n$ in $S_j(t,u_n)$. We use the following expression of $S_j$:
$$
S_j=T+ T\frac{(1-Z)^2(1+Z+Z^2)\mu Z^j}{(1+\mu Z^{j+1})(1+\mu Z^{j+3})},
$$
which is easily derived from~\eqref{product-formSj}. Thus
$$
[t^n] S_j(t,u_n)= 2^nC_n + \frac 1{2i\pi} \int_{\C_n}
 T\frac{(1-Z)^2(1+Z+Z^2)\mu Z^j}{(1+\mu Z^{j+1})(1+\mu Z^{j+3})}
\frac{dt}{t^{n+1}}.
$$
Again, we split the contour $\C_n$ into two parts $\C_n^{(1)}$ and
$\C_n^{(2)}$, shown in Figure~\ref{fig-contour}. As in the proof of
Theorem~\ref{thm-max}, the contribution of
$\C_n^{(1)}$ is easily seen to be $o(8^n/n^m)$ for all $m>0$, thanks
to the results of Lemmas~\ref{lemma-Z} and~\ref{lemma-Sj}.  On
$\C_n^{(2)}$, one has 
$$
t= \frac 1 8 \left( 1 + \frac z n\right)
$$
where $z$ lies in the truncated Hankel contour $\Ha_n$. Conversely,
let $z \in\Ha$. Then $z \in\Ha_n$ for $n$ large enough, and, in
addition to the estimations~\eqref{estimates} already used in the proof of
Theorem~\ref{thm-max}, one finds
\beq
\label{muA}
\mu(t,u_n)= A(a(-z)^{-3/4})(1+o(1)),
\eeq
where $A(x)$ is the series defined by~\eqref{A-eq}. After a few
reductions, one finally obtains
$$
[t^n] S_j(t,u_n)= 2^nC_n +
\frac {12. 8^n n^{-3/2} }{i  \pi}\int_{\Ha }
\frac{A(a(-z)^{-3/4})\exp(-2\lambda (-z)^{1/4})\sqrt{-z}e^{-z}}
{(1+A(a(-z)^{-3/4})\exp(-2\lambda (-z)^{1/4}))^2}dz + o(8^n n^{-3/2}).
$$
It remains to normalize by $2^nC_n=8^n n^{-3/2}/\sqrt \pi$, and then to
set $v=(-z)^{1/4}$ to obtain the expected
expression for the limit of the Laplace transform of $Y_n(j)$, with
$j=\lfloor \lambda n ^{1/4}\rfloor$.

\medskip
The limit Laplace transform $L(\lambda,a)$ is clearly continuous at
$a=0$, and equals $1$ at this point. 
A version of L\'evy's continuity theorem~\cite[Thm.~9.8.2]{dudley}
adapted to Laplace transforms  implies that the sequence 
$Y_n(j)$ converges in distribution to a limit random variable
$Y(\lambda)$ having Laplace transform $L(\lambda, \cdot)$.

\medskip
From the convergence of the Laplace transform in a neighbourhood of
the origin, it is easy to derive that for every $k$,  the sequence of
random variables $Y_n(j)^k$ is \emm uniformly integrable,. 
But then the convergence in distribution implies the convergence of the
moments~\cite[Thm.~5.4]{billingsley2}. This concludes the proof of
Theorem~\ref{thm-conv}.

%%%%%%%%%%%%%%%%%%%%%%%%%%%%%%%%%%%%%%%%%%%%%%%%%
\subsection{Proof of Proposition~\ref{limit-0}}
\label{moments-proof}
When $\lambda=0$, the limiting Laplace transform reduces to
$$
L(0, a)=
1+ \frac {48}{i\sqrt \pi }\int_{\Gamma}
\frac{A(a/v^3)}
{(1+A(a/v^3))^2}v^{5} e^{v^4} dv
=
1+ \frac {12}{i\sqrt \pi }\int_{\Gamma}
\frac{\chi(a/v^3)}
{1+\chi(a/v^3)}v^{5} e^{v^4} dv
$$
where $\chi(x)$ is the unique series in $x$ satisfying
$$
\chi=\frac x 6(1+\chi)^{3/2}.
$$
The Lagrange inversion formula~\cite[p.~38]{stanley-vol2} gives, for $k \ge 1$,
$$
[x^k]  \frac{\chi(x)} {1+\chi(x)}
= \frac 1 {6^k} \frac{\Gamma(3k/2-1)}{k!\Gamma(k/2)}.
$$
Consequently,
$$
L(0, a)=1+ \frac {12}{i\sqrt \pi }
\int_\Gamma
\sum_{k\ge 1} \frac 1 {6^k} \frac{\Gamma(3k/2-1)}{k!\Gamma(k/2)} a^k
v^{5-3k} e^{v^4} dv.
$$
The convergence is absolute, so that we can exchange the sum and the
integral:
$$
L(0, a)=1+ \frac {12}{i\sqrt \pi }
\sum_{k\ge 1} \frac 1 {6^k}
\frac{\Gamma(3k/2-1)}{k!\Gamma(k/2)} a^k
\int_\Gamma
v^{5-3k} e^{v^4} dv.
$$
Using~\eqref{gamma-inverse}, and picking the coefficient of $a^k$, we
find that the 
$k$th moment of the random variable $Y(0)$ is
$$
\EE(Y(0)^k)=\frac {\sqrt \pi}{ 6^{k-1}}
 \frac
{\Gamma(3k/2-1)}{\,\Gamma(k/2)\Gamma((3k-2)/4)}.
$$
The duplication formula,
\[2^{2s-1} \Gamma(s) \Gamma(s+1/2) = \sqrt \pi \, \Gamma(2s),
\]
applied to  $s=(3k-2)/4$, finally gives
$$
 \EE(Y(0)^k)=
 \left( \frac{\sqrt 2} 3\right)^k
\frac {\Gamma(1+3k/4)}{\,\Gamma(1+k/2)} =\lim_{n\rightarrow \infty} 
\EE\left( Y_n( 0) ^k\right).
$$
Since $Y(0)$ has a Laplace transform, it is uniquely determined by its
moments~\cite[Thm.~30.1]{billingsley}. But
$$m_k= \frac
{\Gamma(1+3k/4)}{\,\Gamma(1+k/2)}
$$
is known to be the  $k$th
moment  of $T^{-1/2}$, where $T$ follows a unilateral stable law
of parameter $2/3$
(see~\cite[p.~111]{chaumont-yor}). Proposition~\ref{limit-0} follows.

%%%%%%%%%%%%%%%%%%%%%%%%%%%%%%%%%%%%%%%%%%%%%%%%%
\subsection{Proof of Proposition~\ref{propo-first}}
\label{first-proof}
We have derived above the moments of $Y(0)$ from the expression of its
Laplace transform.  This extends to the moments of $Y(\lambda)$,
for $\lambda >0$: for $k\ge 1$,
$$
\EE(Y(\lambda)^k)= \frac {48. k!}{i\sqrt \pi }\int_{\Gamma}
[a^k]\frac{A(a/v^3)e^{-2\lambda v}}
{(1+A(a/v^3)e^{-2\lambda v})^2}v^{5} e^{v^4} dv.
$$
Since $A(x)=x/24 + O(x^2)$, the case $k=1$ of the above identity reads
$$
\EE(Y(\lambda))= \frac {2}{i\sqrt \pi }\int_{\Gamma}
{e^{-2\lambda v}} v^{2} e^{v^4} dv.
$$
In the above expression, expand the exponential as a series. The
convergence of the sum and integral is absolute, so that one can
exchange them. This gives:
$$
\EE(Y(\lambda))= \frac {2}{i\sqrt \pi }\sum_{m\ge 0} 
\frac{(-2\lambda)^m}{m!}\int_{\Gamma}
 v^{m+2} e^{v^4} dv.
$$
Using~\eqref{gamma-inverse} (which is valid for any $s$ with the
convention $1/\Gamma(-n)=0$ for $n\in \ns$), this can be rewritten as
$$
\EE(Y(\lambda))
={\sqrt \pi }
\sum_{m\ge 0} \frac{(-2\lambda)^m}{m!\Gamma((1-m)/4)} 
=\frac 1 {\sqrt \pi }
\sum_{m\ge 0} \frac{(-2\lambda)^m}{m!}\Gamma\left(\frac{m+3} 4\right) \cos\left(
\frac {(m+1)\pi}4\right).
$$
The last equality follows from the complement formula,
\beq
\label{complement}
\Gamma(s) \Gamma(1-s)= \frac \pi{\sin (\pi s)}.
\eeq

%%%%%%%%%%%%%%%%%%%%%%%%%%%%%%%%%%%%%%%%%%%%%%%%%%%%%%%%%%%%%%
\section{A global limit law, and the distribution function of the ISE}
\label{section-limit-global}
%%%%%%%%%%%%%%%%%%%%%%%%%%%%%%%%%%%%%%%%%%%%%%%%%%
In Section~\ref{section-limit-local}, we have derived from
Proposition~\ref{propo-snake-gf} some \emm 
local, limit results; for instance, a limit law for $X_n(0)/n^{3/4}$,
the (normalized) number of nodes labelled $0$. In this section, we
proceed with a similar study, which aims at deriving from
Proposition~\ref{propo-snake-gf-rep} a \emm global, limit result
--- in particular, 
the limit law of $X_n^+(0)/n$, the normalized number of nodes having a
non-negative label. The technique is copied on
Section~\ref{section-limit-local}, and we do not give all the details.

For $j \in \zs$, let $X^+_n(j)$ denote the random variable equal to the
number of nodes having label at least $j$ in a random tree of
$\mathcal T _{0,n}$. 
Let us define a normalized version of $X^+_n(j)$ by
$$
Y^+_n(j)= \frac{X^+_n(j)}{n}.
$$
These quantities are related to the series $R_j(t,u)$ studied in
Proposition~\ref{propo-snake-gf-rep}. In particular,
$$
\EE\left(e^{aY_n^+(j)}\right)= \EE\left(e^{an^{-1}X_n^+(j)}\right)=
\frac{[t^n] R_j(t,e^{a/n})}{2^n C_n}. 
$$
We extend the definition of $X^+_n$ and $Y^+_n$ to real values
in a natural way by setting $X^+_n(x)=X^+_n(\lceil x\rceil)$ and
$Y^+_n(x)=Y^+_n(\lceil x\rceil)$. 
Let $\lambda \ge 0$.  The aim of this section is to prove that
$Y_n^+( \lambda n^ {1/4})$ converges in distribution, as
$n$ goes to infinity,  to a random variable $Y^+(\lambda)$ that we
describe by its Laplace transform.  This is achieved in
Theorem~\ref{thm-conv-rep} below,  but we first want to present two
consequences of this theorem, which have a simpler formulation. The
first consequence is a striking limit law for $Y^+_n(0)$. Recall
that, up to a normalization by $n$, this random variable gives
the number of nodes having a non-negative label in a tree rooted at $0$.
\begin{Propo}[{\bf The number of non-negative nodes}]
\label{limit-0-G}
 As $n$ goes to infinity, the random variable $Y^+_n(0)$ converges in law to
the uniform distribution on $[0,1]$.
\end{Propo}
\noindent  
This proposition will be proved in Section~\ref{moments-proof-G}.
The second consequence of Theorem~\ref{thm-conv-rep} is an explicit
expansion in $\lambda$ of the limiting first moment of $Y_n^+( \lambda n^{1/4})$.
\begin{Propo}[{\bf The first moment}]
\label{propo-first-G}
 Let $\lambda \ge 0$. Then, as $n$ goes to infinity,
$$
\EE\left( Y_n^+( \lambda n^{1/4}) \right)
\rightarrow
\frac 1 {2 \sqrt\pi}\sum_{m\ge 0} \frac{(-2\lambda)^m}{m!} \cos \left(
\frac{m\pi}4\right)
\Gamma\left(\frac{m+2}4\right).
$$
\end{Propo}
\noindent
This proposition will be proved in Section~\ref{first-proof-G}.

Let us, finally, state our main theorem, from which the two above
propositions derive. 
\begin{Theorem} [{\bf A global limit law}]
\label{thm-conv-rep}
Let $\lambda\ge 0$.% and $j=\lfloor \lambda n^{1/4}\rfloor$. 
The sequence $Y_n^+(\lambda n^{1/4})$ converges in distribution to a
random variable $Y^+(\lambda)$ whose Laplace transform is given, for
$|a|<1$,  by
$$
\EE\left( e^{aY^+(\lambda)}\right)=G(\lambda,a),
$$
where
$$
G(\lambda, a)=
1+ \frac {48}{i\sqrt \pi }\int_{\Gamma}
\frac{B(a/v^4)e^{-2\lambda v}}
{(1+B(a/v^4)e^{-2\lambda v})^2}v^{5} e^{v^4} dv,
$$
\beq
\label{B-eq}
B(x)=-\frac{(1-D)(1-2D)}{(1+D)(1+2D)}, \quad  \quad 
D=\sqrt{\frac{1+\sqrt{1-x}}2},
\eeq
and the integral is taken over
$$
\Gamma= \{1-te^{-i\pi/4}, t\in (\infty,0]\} \cup \{1+te^{-i\pi/4},
t\in [0,\infty)\}. 
$$
Moreover, the Laplace transform of  $Y_n^+(\lambda n^{1/4})$ converges
pointwise to $G(\lambda,\cdot)$ on the interval $(-1,1)$. The
convergence of moments holds as well.

The sequence $Y_n^+(-\lambda n^{1/4})$ converges in distribution to
the random variable $1- Y^+(\lambda)$.
\end{Theorem}
This theorem will be proved in the next subsection. In view of the
following proposition, it tells us about the law of the distribution
function of the ISE.
\begin{Propo}[{\bf The tail distribution function of the ISE}]
\label{coro-tail}
Let $g_{{\hbox{\sc \small ise}}}(\lambda)= \mu_{{\hbox{\sc \small
      ise}}}(\lambda, +\infty)$ denote the tail
distribution function  of the 
ISE. Then for $\lambda \ge 0$,
$$
g_{{\hbox{\sc \small ise}}}(\lambda)
% = \mu_{{\hbox{\sc \small      ise}}}(\lambda, +\infty)
\deq  Y^+(\lambda /\sqrt 2), 
$$
where the law of the variable $Y^+(\lambda )$ is given in
Theorem~{\em\ref{thm-conv-rep}}. In particular,  $g_{{\hbox{\sc \small
      ise}}}(0)$ is uniformly distributed on $[0,1]$. The random
variable $g_{{\hbox{\sc \small ise}}}(-\lambda)$ 
has the same distribution as $1-Y^+(\lambda /\sqrt 2)$.
 \end{Propo}
\noindent{\bf Comments}\\
{\bf 1.} The law of $g_{{\hbox{\sc \small ise}}}(0)$ was already given by
Aldous~\cite[Eq.~(12)]{aldous}.
 \\
{\bf 2.} 
Let us add a few words on the series $B$ and $D$ to convince
ourselves that the integral giving $G(\lambda, a)$ is well-defined as
long as $|a|<1$. Let $E(x)=1-D(x)$. Then $E$ admits the following expansion:
$$
E(x)= 1-\sqrt{1-\frac{1-\sqrt{1-x}} 2} =
2 \sum_{n\ge 1} C_{n-1} \left( \frac
{1-\sqrt{1-x}} 8\right)^n
$$
where $C_n$ is the $n$th Catalan number. Similarly,
$$
{1-\sqrt{1-x}}= 2 \sum_{n\ge 1} C_{n-1} 4^{-n} x^n,
$$
and these two identities imply that $E(x)$ has non-negative
coefficients. Moreover, its radius of convergence is easily seen to
be $1$, so that $|E(x)|\le E(1)=1-1/\sqrt 2 $ for $|x|\le 1$.
Moreover, expressing $B$ in terms of  $E$ gives:
$$
B= \frac{E(1-2E)}{(2-E)(3-2E)},
$$
which shows that $B(x)$ is also analytic for $|x|<1$ and satisfies
in this domain
$$
|B(x)|\le \frac{E(1)(1+2E(1))}{(2-E(1))(3-2E(1))}=22\sqrt 2-31
=0.11...
$$
For $v \in \Gamma$, $|v|\ge 1$ and $\Re(v)\ge 1$. This implies that 
$$
\frac{B(a/v^4)e^{-2\lambda v}}
{(1+B(a/v^4)e^{-2\lambda v})^2}
$$
is uniformly bounded on $\Gamma$, and $G(\lambda, a)$ is well-defined.

\subsection{Proof of Theorem~\ref{thm-conv-rep}}
Let $j= \lceil \lambda n^{1/4}\rceil$. 
Given that the product forms for the series $S_j$ and $R_j$ are very
similar, it is not surprising that we use an approach copied on
that of the previous section. We start from
$$
\EE(e^{aY_n^+(j)})=\EE(u_n^{X_n^+(j)})= \frac{[t^n] R_j(t,u_n)}{2^n C_n},
$$
with $u_n= e^{a/n}$. For technical reasons, we choose to modify
slightly the integration contour of
Figure~\ref{fig-contour}. The Hankel part of this contour, which was
lying at distance $1/8$ of the real axis, is now moved a bit further,
at distance $1/6$ of the real axis. More precisely, the new contour
$\overline \C_n$ consists of
two parts $\overline \C_n^{(1)}$ and  $\overline \C_n^{(2)}$ such that
\begin{itemize}
\item[$\bullet$ ] $\overline \C_n^{(1)}$ is an arc  of radius
$(1+ \log^2 n /n)/8$, centered at the origin;
\item[$\bullet$ ] $\overline \C_n^{(2)}$ is a Hankel contour around $1/8$, at distance
$1/(6n)$ of the real axis, 
which meets  $\overline \C_n^{(1)}$ at both ends.
\end{itemize}
We first need to prove that the series
$R_j(t,u_n)$ is analytic in a neighborhood of   $\overline{\I}_n$, the region
lying inside the integration contour $\overline\C_n$. The following lemma
is the counterpart of Lemma~\ref{lemma-Sj}.
\begin{Lemma}[{\bf Analytic properties of $\nu$ and $R_j$}]
\label{lemma-Rj}
Let $a$ be a real number such that $|a|<1$.
Then 
$\nu(t,u_n)$ is analytic in a neighborhood of $\overline{\I}_n$. 
Its modulus in  $\overline{\I}_n$ is smaller than $\alpha$, for some $\alpha<1$
independent of $a$ and $n$.
The series $R_j(t,u_n)$ is also analytic in a neighborhood of
$\overline{\I}_n$. 
\end{Lemma}
\noindent
{\bf Proof.} Again, the lemma is obvious if $a=0$. We thus assume
$a\not =0$ and $|a|<1$.

Let us first study the singularities of the series  $\bar \nu\equiv
\bar \nu(z,u)$ defined by
$$
\bar \nu(z,u)= \nu\left( \frac{z(1+z^2)}{(1+z)^4},u\right).
$$
According to Proposition~\ref{propo-snake-gf-rep}, $\bar \nu$  is a \fps \ in
$z$ with polynomial coefficients in $u$, and by~\eqref{eqZt}, one has:
$$
\nu(t,u)= \bar \nu(Z(t), u).
$$
In the course of the proof of Proposition~\ref{propo-snake-gf-rep}, we
have obtained a polynomial equation $P(\nu,Z,u)=0$, of degree 4 in
$\nu$, relating  $\nu(t,u), Z(t)$ and the variable
$u$. This equation is not written in the paper (it is a bit too big),
but it can be easily obtained using the expression of $\nu$ given in
Proposition~\ref{propo-nu}. By definition of $\bar \nu$, we have
$P(\bar \nu, z, u)=0$.

 Assume that $u$ is a fixed real
number close to $1$. That is, $u=1+x$, with $x$ small. In order to
study the singularities of $\bar \nu$, we look again at the zeroes of
the leading coefficient of $P$, and at the zeroes of its
discriminant. This gives several candidates for the
singularities of $\bar \nu(z,u)$, which we classify in three series
according to their behaviour when $x$ is small. First, some
candidates tend to a limit that is different from 1,
$$
z=- 1, \quad z= \pm i,  \quad 
z= e^{\pm 2i\pi/3},  \quad 
z= e^{\pm 2i\pi/3} + O(x), \quad 
z=-1+O(x).
$$
Then, some candidates tend to 1 and lie at distance at most
$|x|^{1/4}$ of 1 (up to a multiplicative constant):
$$
z= 1+ \omega (cx)^{1/4} + O(\sqrt{|x|}),
$$
%z=1,
%z=1+2 \omega |x|^{1/4} + O(\sqrt{|x|}),
%z=1+ \omega\left( \frac{64 |x|}3\right)^{1/4} + O(\sqrt{|x|}),
%z=1+ e^{i\pi/4} \omega\left( \frac{16 |x|}3\right)^{1/4} + O(\sqrt{|x|}),
%$$
where $\omega$ is a fourth root of unity and  $c$ is in the set $\{ 0,
16,  {64 }/3, -16/3\}$. 
Finally, some candidates tend to 1 but lie further away from 1 (more
precisely, at distance $|x|^{1/6}$):
$$
z=1+ 2 e^{i\pi/6} \omega' x^{1/6} + O({|x|^{1/3}}),
$$
where $\omega'$ is a sixth root of unity.

Let us now consider $\nu(t,u)=\bar \nu(Z(t),u)$ with $u=u_n=
e^{a/n}=1+x$, where $x=a/n(1+o(1))$. Recall that $Z$ is
analytic in $\cs \setminus [1/8, \infty)$. By Lemma~\ref{lemma-Z}, the 
  series $Z(t)$ never approaches any root of unity different from
  1. Hence for $n$ large enough, $Z(t)$ never reaches any of
  the candidates $z$ of the first series.

The candidates of the second series are of the   form
$$
z= 1 + \omega  (ac/n)^{1/4} + O( n^{-1/2})
$$
for some constant $c$, with $|c|\le 64/3$, depending on the candidate.
By~\eqref{eqZt}, $Z(t)$ may only reach these values for
$$
t=\frac 1 8 -  \frac{ac}{128n} + O(n^{-5/4}).
$$
Since $|a|<1$, there exists $\epsilon>0$ such that these  values lie
at distance less that $1/((6+\epsilon)n)$ of $1/8$, that is, outside a
neighborhood of the domain $\overline{\I}_n$.

The candidates of the third series are more worrying: $Z(t)$ may
reach them for
\beq
\label{t-worrying}
t=\frac 1 8 - \frac{\omega ''}{8}  \left(a/n\right)^{2/3} + O(n^{-5/6}), 
\eeq
where $\omega''$ is a cubic root of unity,
and these values may lie inside $\overline{\I}_n$.
If $a>0$ and $\omega''=e^{\pm 2i\pi/3}$, or if $a<0$ and $\omega''=e^{
  2i\pi/3}$, the modulus of the above value of $t$ is 
found to be $1/8(1+ cn^{-2/3} + o(n^{-2/3}))$, for some \emm positive,
constant $c$: this is larger than the radius of the contour $\overline
\C_n$, which implies that $t$ lies outside a neighborhood of
$\overline{\I}_n$. However, if 
 $a>0$ and $\omega''=1$, or  if $a<0$ and  $\omega''=1$ or $e^{-
  2i\pi/3}$,  the above value of $t$  lies definitely  inside
$\overline{\I}_n$.  Its modulus is $1/8(1- cn^{-2/3} +
o(n^{-2/3}))$, for some  positive constant $c$. 

In order to rule out the possibility that $\nu(t,u_n)$ has such a
singularity, we are going to prove, by having a close look at the 
expression of $\nu$ given in Proposition~\ref{propo-nu}, that the
radius of convergence of $\nu(t,u_n)$ is at least $1/8-O(1/n)$.
We use below the notation of Proposition~\ref{propo-nu}.

Clearly, the series $V(t,u_n)$ has radius of convergence $\min
(1/8,1/(8u_n))$. In particular, this radius is at least
$\rho_n:=1/(8(1+|x|))$ (with $u_n=1+x$).
 Moreover,  the series $V$ admits the following
expansion 
$$
V(t,1+x)= \frac 1 4 \left( 1-\sqrt{1-\frac{8tx}{1-8t}}\right)
= \frac 1 2 \sum_{n \ge 1} C_{n-1} \left( \frac{2tx}{1-8t}\right)^n,
$$
where $C_n$ is the $n$th Catalan number. This shows that $V(t,1+ |x|)$
is a series in $t$ with \emm positive, coefficients and that for all
$t$ such that $|t|\le \rho_n$,
$$
\left|V(t,1+x)\right| \le V(|t|, 1+ |x|) \le V\left( \frac 1{
  8(1+|x|)},   1+ |x|\right) = \frac 1 4 .
$$
The next step is to prove that  $\Delta(t,u_n)$ never vanishes for $|t|\le
\rho_n$. Indeed,
$$
\Delta=
%(1-V)^2-\frac{4 ZV^2}{(1+Z)^2} = 
(1-V)^2 - 4WV^2,
$$
where $W\equiv W(t)$ is the formal power series in $t$ defined by~
\eqref{eqW}. This series has radius $1/8$, and non-negative
coefficients.
Hence for all $t$ such that $|t|\le 1/8$, one has $|W(t)|\le W(1/8)=
1/4$. Consequently, for  $|t|\le \rho_n$,
$$
|\Delta(t,1+x)|\ge (1-|V(t,1+x)|)^2 - 4|W(t)||V(t,1+x)|^2 \ge \left( 1 - \frac 1
4\right)^2 - \frac 1 {16} = \frac 1 2.
$$
Hence $\Delta(t,u_n)$ does not vanish in the centered disk of
radius $\rho_n$. It follows that the series $P(t,u_n)$ is analytic
inside this disk. 

According to the expression of $\nu$  given
in Proposition~\ref{propo-nu}, the series $\nu(t,u_n)$ is meromorphic
for   $|t|\le \rho_n$. The final question we need to answer is whether
$\nu$ has poles in this disk, and where. Returning to the polynomial
$P$ such that $P(\nu, Z, u)=0$ shows that this can only happen if the
coefficient of $\nu^4$ in this polynomial vanishes. But this can only
occur if $z=Z(t)$ has one of the following forms:
$$
z=\pm i, \quad z= e^{\pm 2i\pi/3} +O(x), \quad z=-1 + O(x), \quad z= 1+ \omega
(64x/3)^{1/4} + O(x^{1/2}).
$$
 As argued above, only the last value of $z$ is likely to be reached
by $Z(t)$, and this may only occur if
$$
t=\frac 1 8 - \frac{1}{6} (a/n) + O(n^{-5/4}).
$$
Consequently, the radius of $\nu(t,u_n)$ is at least $1/8-O(1/n)$,
and this proves that the values~\eqref{t-worrying} that have been
shown to lie in the centered disk of radius $1/8$, 
are not, after all, singularities of $\nu(t,u_n)$. This completes our
proof that $\nu(t,u_n)$ is analytic in a neighborhood of  $\overline{\I}_n$.

\medskip
We now want to bound $\nu(t,u_n)$ inside $\overline\I_n$. From now on,
we can walk safely in the steps of the proof of Lemma~\ref{lemma-Sj}.
Let $t_n\in \overline\I_n$ be such that
$$
|\nu(t_n,u_n)|= \max _{t\in \overline\I_n} |\nu(t,u_n)|.
$$
We first give a lower bound for this quantity, by estimating
$\nu(t,u_n)$ for $t=1/8-1/(6n)$. This is easily done by combining the
closed form expressions of $\nu$ (Proposition~\ref{propo-nu}) and $B$
(Theorem~\ref{thm-conv-rep}). One obtains:
$$
|\nu(t_n,u_n)| \ge |\mu(1/8-1/(6n), u_n)|=|B(3a/4)|+o(1)>0.
$$
This lower bound is then used to rule out the possibility that the
sequence $t_n$ has an accumulation point different from $1/8$. Thus
$t_n$ converges to $1/8$, and one can write
$$
t_n= \frac 1 8 \left( 1 -\frac {x_n} n \right).
$$
We have $x_n=o(n)$, but also $|x_n|>4/3$ since $t_n$ belongs to
$\overline \I_n$. We want to estimate $\nu(t_n, u_n)$.
Since
$$
Z(t_n)= 1-2  \left( \frac {x_n} n \right)^{1/4} + O\left(  \left(
\frac {x_n} n \right)^{1/2}\right)
$$
and
$$
u_n-1=a/n  \left( 1 + O(1/n)\right),
$$
one has
$$
\frac{u-1}{(1-Z)^4}= \frac{a}{16x_n}\left(1+ O\left(  \left(
\frac {x_n} n \right)^{1/4}\right)\right).
$$
The closed form expressions of $\nu$ and $B$ imply that the sequence
$x_n$ is bounded and
$$
\nu(t_n,u_n)= B(a/x_n) + o(1).
$$
Since $B$ is bounded by  $0.12$ inside its disk of convergence,
$|\nu(t_n,u_n)|$ is certainly smaller than some $\alpha$ for
$\alpha<1$ and $n$ large enough.
This concludes the proof of the second statement of Lemma~\ref{lemma-Rj}.

\medskip
By continuity of $\nu(t,u_n)$, this function of $t$ is still bounded by $1$
(in modulus) is a neighborhood of $\overline \I_n$. 
Recall also that the modulus
of $Z(t)$ never reaches $1$ for $t \in \cs \setminus [1/8,
  \infty)$. The form~\eqref{product-formRj} then implies that
  $R_j(t,u_n)$ is an analytic function of $t$ in a neigborhood of
  $\overline \I_n$. 
\cqfd

Once this rather painful lemma is at last established, the rest of the
proof of Theorem~\ref{thm-conv-rep} copies  the end of the
proof of Theorem~\ref{thm-conv}, with $S_j, \mu$ and $A$ respectively
replaced by $R_j$,  $\nu$ and $B$. The counterpart of~\eqref{muA} is
$$
 \nu(t,u_n)= B(-a/z) (1+o(1)).
$$
Recall that the Hankel part of the contour $\overline \C_n$ is now at
distance $1/(6n)$ of the real axis. Hence, when $n$ goes to infinity,
one finds
 $$
[t^n] R_j(t,u_n)= 2^nC_n +
\frac {12. 8^n n^{-3/2} }{i  \pi}\int_{4/3 \Ha }
\frac{B(-a/z)\exp(-2\lambda (-z)^{1/4})\sqrt{-z}e^{-z}}
{(1+B(-a/z)\exp(-2\lambda (-z)^{1/4}))^2}dz + o(8^n n^{-3/2}).
$$
After normalizing by $2^n C_n$ and setting $v=(-z)^{1/4}$, this gives
$$
\EE(e^{aY_n^+(j)}) \rightarrow
1+ \frac {48}{i\sqrt \pi }\int_{(4/3)^{1/4}\Gamma}
\frac{B(a/v^4)e^{-2\lambda v}}
{(1+B(a/v^4)e^{-2\lambda v})^2}v^{5} e^{v^4} dv,
$$
but the analyticity properties of the integrand allow us to replace
the integration contour by $\Gamma$.

%%%%%%%%%%%%%%%%%%%%%%%%%%%%%%%%%%%%%%%%%%%%%%%%%
\subsection{Proof of Proposition~\ref{limit-0-G}}
\label{moments-proof-G}
When $\lambda=0$, the limiting Laplace transform reduces to
$$
G(0, a)=
1+ \frac {48}{i\sqrt \pi }\int_{\Gamma}
\frac{B(a/v^4)}
{(1+B(a/v^4))^2}v^{5} e^{v^4} dv
=
1+ \frac {4}{3 i\sqrt \pi }\int_{\Gamma}
\frac{\chi(a/v^4)(3-\chi(a/v^4))}
{1+\chi(a/v^4)}v^{5} e^{v^4} dv
$$
where $\chi(x)$ is the unique \fps\ in $x$  satisfying
$$
\chi=\frac x 4(1+\chi)^{2}.
$$
The Lagrange inversion formula gives, for $k \ge 1$,
$$
[x^k]  \frac{\chi(x)(3-\chi(x))} {1+\chi(x)}
= \frac {6} {4^k} \frac{(2k-2)!}{(k-1)!(k+1)!}.
$$
Consequently,
$$
G(0, a)=1+ \frac {8}{i\sqrt \pi }\int_{\Gamma}
\sum_{k\ge 1} \frac 1 {4^k} \frac{(2k-2)!}{(k-1)!(k+1)!} a^k
v^{5-4k} e^{v^4} dv.
$$
The convergence is absolute, so that we can exchange the sum and the
integral:
$$
G(0, a)=1+\frac {8}{i\sqrt \pi } \sum_{k\ge 1} \frac 1 {4^k}
 \frac{(2k-2)!}{(k-1)!(k+1)!} a^k
\int_{\Gamma}v^{5-4k} e^{v^4} dv.
$$
Using~\eqref{gamma-inverse}, and picking the coefficient of $a^k$, we find that the
$k$th moment of the random variable $Y^+(0)$ is
$$
\EE(Y(0)^k)=  \frac {8}{i\sqrt \pi } 
 \frac 1 {4^k}
 \frac{(2k-2)!}{(k-1)!(k+1)}\frac {i \pi}{2\Gamma(k-1/2)} = \frac 1 {k+1}.
$$
The unique distribution having  its $k$th
moment equal to $1/(k+1)$  is the uniform distribution on
$[0,1]$. Proposition~\ref{limit-0-G} follows.

%%%%%%%%%%%%%%%%%%%%%%%%%%%%%%%%%%%%%%%%%%%%%%%%%
\subsection{Proof of Proposition~\ref{propo-first-G}}
\label{first-proof-G}
We have derived above the moments of $Y^+(0)$ from the expression of its
Laplace transform.  This extends to the moments of $Y^+(\lambda)$,
for $\lambda >0$: for $k\ge 1$,
$$
\EE(Y^+(\lambda)^k)= \frac {48. k!}{i\sqrt \pi }\int_{\Gamma}
[a^k]\frac{B(a/v^4)e^{-2\lambda v}}
{(1+B(a/v^4)e^{-2\lambda v})^2}v^{5} e^{v^4} dv.
$$
Since $B(x)=x/48 + O(x^2)$, the case $k=1$ of the above identity gives
$$
\EE(Y^+(\lambda))= \frac {1}{i\sqrt \pi }\int_{\Gamma}
{e^{-2\lambda v}} v e^{v^4} dv.
$$
In the above expression, expand the exponential as a series. The
convergence of the sum and integral is absolute, so that one can
exchange them. This gives:
$$
\EE(Y^+(\lambda))= \frac {1}{i\sqrt \pi }\sum_{m\ge 0} 
\frac{(-2\lambda)^m}{m!}\int_{\Gamma}
 v^{m+1} e^{v^4} dv.
$$
Using~\eqref{gamma-inverse}, this can be rewritten as
$$
\EE(Y^+(\lambda))=\frac {\sqrt \pi } 2
\sum_{m\ge 0} \frac{(-2\lambda)^m}{m!\Gamma((2-m)/4)} 
=\frac 1 {2 \sqrt \pi }
\sum_{m\ge 0} \frac{(-2\lambda)^m}{m!}\Gamma\left(\frac{m+2} 4\right) \cos\left(
\frac {m\pi}4\right).
$$
The last equality follows from the complement formula~\eqref{complement}.

\subsection{The distribution function of the {\sc ISE}}
\label{section-coro-tail}
\noindent Let us finally prove Proposition~\ref{coro-tail}.

Let $\mu_n$ be a sequence of random probability measures on $\rs$ converging
weakly to a probability measure $\mu$. Let $F_n$ denote the (random)
distribution function of $\mu_n$: for $ \lambda \in \rs$, 
$$
F_n(\lambda) = \mu_n(-\infty, \lambda].
$$
Similarly, let $F$ denote the distribution function of $\mu$. It is
not very hard 
to prove that, for all $\lambda\in\rs$ such that $\mu\{\lambda\}=0$,
$F_n(\lambda)$ converges in distribution to $F(\lambda)$.
(We prove this in the appendix of the paper, but it is certainly
written somewhere in the literature.)

Let us now apply this general result to our context. The probability
measure $\mu_n$ is given by~\eqref{random-measure}, with $c=\sqrt
2$. It is known to 
converge to the random measure $\mu_{\hbox{\sc \small ise}}$. 
Assume for the moment  that this measure does not assign a positive weight
to any point. Then, with the above notation, $F_n(\lambda)$ converges in
distribution to $F(\lambda)$, for all $\lambda\in \rs$.  But, given
the definition~\eqref{random-measure} of $\mu_n$,
$$
F_n(\lambda)
%=\mu_n(-\infty, \lambda]
= 1-\mu_n(\lambda, \infty)= 1
-\frac 1{n+1} X_n^+( \lambda n^{1/4}/\sqrt 2)  
+ \frac 1{n+1} X_n( \lambda n^{1/4}/\sqrt 2),
$$
where the definition of $X_n$ is extended to all reals  by
$X_n(x)=0$ if $x\in \rs\setminus \zs$. By Theorems~\ref{thm-conv} and~\ref{thm-conv-rep}, 
the right-hand side converges in distribution to $1-Y^+(\lambda/\sqrt
2)$.  Consequently, the tail distribution function of the ISE (that
is, $\mu_{\hbox{\sc \small ise}}(\lambda, \infty)$) has the
same law as $Y^+(\lambda/\sqrt2)$.

It remains to prove  that $\mu_{\hbox{\sc \small ise}}$ does not weight
  points positively (almost surely). Let $\lambda\in \rs$. Then
\beq
\label{E-P}
\PP(\mu_{\hbox{\sc \small ise}}\{\lambda\}>0)=0 \Longleftrightarrow \EE(\mu_{\hbox{\sc \small ise}}\{\lambda\})=0.
\eeq
Let $\epsilon >0$, and let $h_{\lambda, \epsilon}$ be the function
plotted on the right-hand side of Figure~\ref{fig-epsilon}.
Then
$$
\begin{array}{lllll}
\mu_{\hbox{\sc \small ise}}\{\lambda\}&=&
\displaystyle \lim_{\epsilon \rightarrow 0^+} \int_\rs h_{\lambda, \epsilon}(y) d\mu_{\hbox{\sc \small ise}}(y)\\
&=&\displaystyle\lim_{\epsilon \rightarrow 0^+}\int_0^1 h_{\lambda, \epsilon}(
r(t)) dt & (\hbox{by } \eqref{ise-rep})\\
&=& \displaystyle\int_0^1 \1_{\lambda=  r(t)} dt.
\end{array}
$$
Taking expectations, we obtain
$$
\EE(\mu_{\hbox{\sc \small ise}}\{\lambda\})
= \EE\left(\int_0^1 \1_{\lambda=  r(t)}
dt\right) = \int_0^1 \PP(\lambda=  r(t)) dt.
$$
But $\PP(\lambda=  r(t))=0$ for all $t\in (0,1)$ and $\lambda$, since
$r(t)$ has a density with respect to the Lebesgue measure for all
$t$. By~\eqref{E-P}, we conclude that $\mu_{\hbox{\sc \small
    ise}}$ does not weight points positively.

The last statement of Proposition~\ref{coro-tail} is then easily proven,
using a the symmetry of $\mu_{\hbox{\sc \small ise}}$ and the fact
that it does not assign a positive probability to any point.

%%%%%%%%%%%%%%%%%%%%%%%%%%%%%%%%%%%%%%%%%%%%%%%%%%%%%%%%%%%%
%%%%%%%%%%%%%%%%%%%%%%%%%%%%%%%%%%%%%%%%%%%%%%%%%%%%%%%%%%%%%%
\section{Other tree models and universality}
\label{section-universality}
\subsection{Trees with increments $0, \pm 1$}
We consider in this section a slight variation on the previous family
of trees: the increments of the labels along edges may now be $0, \pm
1$. This family of trees has attracted a lot of interest in relation to
planar
maps~\cite{bdg-statistics,chassaing-schaeffer,cori-vauquelin,jf-mokka-brownianmap}. 
\subsubsection{Enumerative results}
\label{section-enumerative-0pm1}
As above, let $T_j \equiv T_j(t)$ be the \gf\ of labelled trees
in which all labels are at most $j$, counted by their number of
edges. Let $S_j \equiv S_j(t,u)$ be the 
\gf\ of labelled trees, counted by the number of edges (variable $t$)
and the number of nodes labelled $j$ (variable $u$). Finally, let $R_j
\equiv R_j(t,u)$ be the 
\gf\ of labelled trees, counted by the number of edges
and the number of nodes having label $j$ at least. 
As above, it is easy to write an infinite system of equations defining
any of the families $T_j$, $S_j$ or $R_j$. The only difference with
our first family of trees  is that
a third case arises in the decomposition of trees illustrated
by Figure~\ref{fig-decomp}: the leftmost child of the root
may have label $j$. In particular,  the \gf \  $T\equiv T(t)$ counting
plane labelled trees now satisfies
$$
T=1+3tT^2,
$$
while  for $j\ge 0$,
\beq
\label{rec-TjP}
T_j=1+t(T_{j-1}+T_j+T_{j+1})T_j.
\eeq
The equations of Lemmas~\ref{lemma-Sj-rec} and~\ref{lemma-Rj-rec} are
modified in a similar way. 
The three infinite systems of equations thus obtained can be solved
using the same 
techniques as in Section~\ref{section-snake-gf}. The solutions are
expressed in terms of the above
series $T\equiv T(t)$ and of the  unique formal power series
$Z\equiv Z(t)$, with constant term $0$, satisfying
\beq
\label{eqZtP}
Z=t\, \frac{(1+4Z+Z^2)^2}{1+Z+Z^2}.
\eeq
Observe that $T$ and $Z$ are related by:
%\beq
%\label{eqTZP}
$$
T=\frac{1+4Z+Z^2}{1+Z+Z^2}.
$$ %\eeq
We state without proof the counterparts of
Propositions~\ref{propo-snake-bounded},~\ref{propo-snake-gf} 
and~\ref{propo-snake-gf-rep}. 
\begin{Propo}[{\bf Trees with small labels \cite{bdg-geodesic,bdg-statistics}}]
\label{propo-snake-boundedP}
Let $T_j\equiv T_j(t)$ be the \gf\ of trees having no label greater
than $j$. Then $T_j$ is algebraic of degree (at most) $2$. In particular,
$$
T_0=1-16\,t
+18\,t  {T_0}
-27\,{t}^{2}{{ T_0}}^{2}.
$$
Moreover,  for all $j \ge -1$,
%\beq\label{product-formTjP}
$$
T_j= T\,\frac{(1- Z^{j+1})(1- Z^{j+4})}{(1- Z^{j+2})(1-Z^{j+3})},
%\eeq
$$
where $Z\equiv Z(t)$ is given by~\eqref{eqZtP}.
\end{Propo}
\noindent
{\bf Remarks}\\
{\bf 1.} As observed in~\cite[p.~645]{bdg-statistics}, there is an ``invariant'' function
attached to equations of the form~\eqref{rec-TjP}:
 for $j\ge 0$,
$$
I(T_{j-1},T_j)=I(T_{j},T_{j+1})
$$
where  $I$ is now given by
$$
I(x,y)= xy(1-t(x+y))-x-y.
$$
As explained in the remark that follows
Propositions~\ref{propo-snake-gf} and~\ref{propo-mu-closed-form},
this can be used to derive rapidly from~\eqref{rec-TjP} the value of 
$T_0$.\\
{\bf 2.} As was the case for trees with increments $\pm1$,
  the trees counted by $T_0$ (equivalently, the trees having only
non-negative labels) are closely related to planar
maps. More precisely, there is a one-to-one correspondence between
non-negative trees having $n$ edges and planar maps having $n$
edges~\cite{chassaing-schaeffer,cori-vauquelin}. 
The coefficients of $T_0(t)$ are also remarkably simple:
$$
T_0(t)= 
\frac{(1-12t)^{3/2}-1+18t}{54t^2}
=\sum_{n\ge 0} \frac{2. 3^{n}}{(n+1)(n+2)}{{2n}\choose n}t^n.
$$
A combinatorial explanation for the algebraicity of $T_0$ is given
in~\cite{cori-vauquelin}.

\begin{Propo}[{\bf The number of nodes labelled $j$}]
\label{propo-snake-gfP}
For any $j\in \zs$, the \gf \ $S_j\equiv S_j(t,u)$ that counts
labelled trees by the number of edges and the number of nodes labelled $j$ is
algebraic of degree at most $4$ over $\qs(T,u)$ (and hence has degree
at most $8$ over $\qs(t,u)$). 
More precisely,
$$
\frac{9T^4(u-1)^2}{(T-S_0)^2}= 9T^2-2T(T-1)(2T+1)S_0+(T-1)^2S_0^2
%\eeq
$$
and all the $S_j$ belong to $\qs(t,u,S_0)$.
Moreover,  for all $j \ge 0$,
%\beq\label{product-formSjP}
$$
S_j= T\,\frac{(1+\mu Z^j)(1+\mu Z^{j+3})}{(1+\mu Z^{j+1})(1+\mu Z^{j+2})},
%\eeq
$$
where $Z\equiv Z(t)$ is given by~\eqref{eqZtP}
and $\mu\equiv \mu(t,u)$ is the unique \fps\ in $t$
satisfying
%\beq\label{eqMZP}
$$
\mu= (u-1) \frac{(1+Z+Z^2)(1+\mu Z)^2(1+\mu Z^2)^2}
{(1+Z)^2(1-Z)^3(1-\mu^2 Z^3)}.
%\eeq
$$
The series $\mu(t,u)$  has polynomial coefficients in $u$,
and satisfies $\mu(t,1)=0$. It has degree $4$ over $\qs(Z,u)$ and $16$
over $\qs(t,u)$. 
\end{Propo}

\begin{Propo}[{\bf  The number of nodes labelled $j$ or more}]
\label{propo-snake-gf-repP}
Let $j \in \zs$. The \gf\ $R_j(t,u)\equiv R_j$
that counts labelled trees by the number of edges and the number of
nodes labelled $j$ or more
is  algebraic over $\qs(t,u)$, of degree at most $8$. 
 It has degree
at most $2$ over $\qs(T,\tilde T)$, where $T\equiv T(t)$ and $\tilde
T\equiv T(tu)$. More precisely, it belongs to the
extension of $\qs(T,\tilde T)$ generated by 
$$
\sqrt{4(T+\tilde T)^2-T\tilde T (4+3T\tilde T)}.
$$
Moreover,  for all $j \ge 0$,
%\beq\label{product-formRjP}
$$
R_j= T\,\frac{(1+\nu Z^j)(1+\nu Z^{j+3})}{(1+\nu Z^{j+1})(1+\nu Z^{j+2})},
%\eeq
$$
where $Z\equiv Z(t)$ is given by~\eqref{eqZtP}
and $\nu\equiv \nu(t,u)$ is a \fps\ in $t$, with polynomial
coefficients in $u$, which is algebraic
%  of degree $2$ over $\qs(Z(t),Z(tu))$, 
of degree $4$ over $\qs(u,Z)$, and  of degree $16$
over $\qs(t,u)$. 
This series satisfies $\nu(t,1)=0$.
The first terms in its expansion are:
$$
\nu(t,u)=(u-1)\Big(
1+3\,ut+ \left(15\,u+ 14\,{u}^{2} \right) {t}^{2}
+ \left( 104\,u+117\,{u}^{2}+83\,{u}^{3} \right) {t}^{3}+
O(t^4)\Big) .
$$
\end{Propo}

\subsubsection{Limit laws}
We now endow the set of labelled trees having $n$ edges with the
uniform distribution, and consider the same random variables as for
our first family of trees: $M_n$,  the largest label,  $X_n(j)$, the
number of nodes having label $j$, % in a random tree of  rooted at $0$,
and finally  $X_n^+(j)$, the
number of nodes having label $j$ at least.

Again, we can prove that 
 $M_nn^{-1/4}$ converges in law to  $N_{{\hbox{\sc \small
   ise}}}/\sqrt 3 $,  where $N_{{\hbox{\sc \small
   ise}}}$ is the supremum of the support of the ISE, and that
 for all $\lambda\in \rs$, the sequence $X_n^+( \lambda n^{1/4})/ n$
 converges in law to  $g_{{\hbox{\sc \small ise}}}(\sqrt{3 }\lambda 
)$ where $g_{{\hbox{\sc \small ise}}}$ is the tail
distribution function of the ISE. The arguments are the same as for
 our first class of trees (Sections~\ref{section-supremum}
 and~\ref{section-coro-tail}). 

Hence we could  just as well have started from the enumerative results of
Section~\ref{section-enumerative-0pm1}, rather than from those of
Section~\ref{section-snake-gf}, to characterize the laws of
$N_{{\hbox{\sc \small  ise}}} $ and $g_{{\hbox{\sc \small
      ise}}}(\lambda)$ (Propositions~\ref{coro-support}
and~\ref{coro-tail}). More remarkably, we have performed on
$X_n(j)$ an analysis similar to that of
Section~\ref{section-limit-local}, and obtained \emm the same, local
limit law. In other words, for all $\lambda \ge 0$,
the sequence  $X_n(\lfloor \lambda n^{1/4}\rfloor) n^{-3/4}$ converges
in law to  $
\sqrt 3  f_{{\hbox{\sc \small ise}}}(\sqrt 3 \lambda )
$
where $f_{{\hbox{\sc \small ise}}}$ is the conjectured
density of the ISE, given in Conjecture~\ref{conj-density}.

In all three cases, the convergence of the moments holds as well.

%%%%%%%%%%%%%%%%%%%%%%%%%%%%%%%%%%%%%%%%%%%%%%%%%%
\subsection{Naturally embedded binary trees}
 We study in the section the  incomplete binary
trees\footnote{The author has obtained similar, but slightly heavier
  results for embedded  \emm complete,\/  binary trees.}
carrying their natural labelling, as shown on the right of
Figure~\ref{fig-three-families}. Such trees are either
empty, or have a root, to 
which a left and right subtree (both possibly empty) are
attached. A (minor) difference with the two previous  families of trees
is that the main enumeration parameter is now the number of nodes
rather than the number of edges.

\subsubsection{Enumerative results}
\label{section-enumerative-binary}
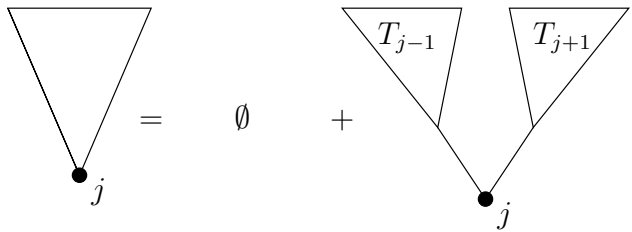
\begin{figure}[hbt]
\begin{center}
\input{tree-decomp-binary.pstex_t}
\end{center}
\caption{The decomposition of naturally labelled binary trees rooted
  at $j$.}
\label{fig-decomp-binary}
%\hrule
\end{figure}

 Let $T_j \equiv T_j(t)$ be the \gf\ of (naturally labelled) binary trees
in which all labels are at most $j$, counted by their number of
 nodes. Let $S_j \equiv S_j(t,u)$ be the 
\gf\ of binary trees, counted by the number of  nodes (variable $t$)
and the number of nodes labelled $j$ (variable $u$). Finally, let $R_j
\equiv R_j(t,u)$ be the 
\gf\ of binary  trees, counted by the number of nodes
and the number of nodes having label $j$ at least. 
It is easy to write an infinite system of equations defining
any of the families $T_j$, $S_j$ or $R_j$. The decomposition of trees
that was crucial in Section~\ref{section-snake-gf} is now replaced by the
 decomposition of Figure~\ref{fig-decomp-binary}. 
% Observe that there is now  no choice for   the labels of the
% children of the root. In particular,  
The \gf \ $T\equiv T(t)$ counting 
naturally labelled binary trees satisfies
$$
T=1+tT^2,
$$
(as it should!) while  for $j\ge 0$,
\beq
\label{rec-TjB}
T_j=1+tT_{j-1}T_{j+1}.
\eeq
Note that the initial condition is now $T_{-1}=1$ (accounting for
the empty tree). 

The equations of Lemmas~\ref{lemma-Sj-rec} and~\ref{lemma-Rj-rec}
respectively become:
\beq
\label{rec-Sj-binary}
S_j=\left\{
\begin{array}{llll}
1+tS_{j-1}S_{j+1} & \hbox{if } j\not = 0, \\
1+tuS_1^2  & \hbox{if } j = 0,
\end{array}
\right.
\eeq
while
\beq
\label{rec-RjB}
R_j= 1+tR_{j-1}R_{j+1}  \quad \hbox{for } j\ge1,
\eeq
and
\beq
\label{eq-symmetryB}
R_{-j}(t,u)=R_{j+1}(tu,1/u)\quad \hbox{for all } j \in \zs.
\eeq

The three infinite systems of equations thus obtained can be solved
using the same 
techniques as in Section~\ref{section-snake-gf}. The solutions are
expressed in terms of the above
series $T\equiv T(t)$ and of the  unique formal power series
$Z\equiv Z(t)$, with constant term $0$, satisfying
\beq
\label{eqZtB}
Z=t\, {\frac{ \left( 1+{Z}^{2} \right) ^{2
}} {  1-Z+{Z}^{2}  }}.
\eeq
Observe that $T$ and $Z$ are related by:
$$
T=\frac{1+Z^2}{1-Z+Z^2}.
$$ %\eeq
We state without proof the counterparts of
Propositions~\ref{propo-snake-bounded},~\ref{propo-snake-gf} 
and~\ref{propo-snake-gf-rep}. Once again, the results below are dying
for combinatorial explanations!
\begin{Propo}
Let $T_j\equiv T_j(t)$ be the \gf\ of binary trees having no label greater
than $j$. Then $T_j$ is algebraic of degree (at most) $2$. In particular,
$$
T_0= \frac {(1-4t)^{3/2}-1+8t-2t^2}{2t(1+t)}.
$$
Moreover,  for all $j \ge -1$,
%\beq\label{product-formTjB}
$$
T_j= T\,\frac{(1- Z^{j+2})(1-Z^{j+7})}{(1- Z^{j+4})(1- Z^{j+5})},
%\eeq
$$
where $Z\equiv Z(t)$ is given by~\eqref{eqZtB}.
\end{Propo}
\noindent
It is easy to check that the above series $T_j$ satisfy the
equations~\eqref{rec-TjB} and the initial condition $T_{-1}=1$. The
method we used to 
\emm discover,  this product form  is again borrowed
from~\cite{bdg-geodesic}.
\smallskip

\noindent{\bf Remark.}
For this family of trees as well, we have found  an ``invariant'' function
attached to equations of the form~\eqref{rec-TjB}:
 for $j\ge 0$,
$$
I(T_{j-1},T_j)=I(T_{j},T_{j+1})
$$
where 
$$
I(x,y)=\left( x+y \right) {t}^{2}+{\frac { \left( {x}^{2}-x-y+{y}^{2} \right) t}{xy}}+{\frac {-1+x+y}{xy}}.
$$
This can be used to derive rapidly from~\eqref{rec-TjB} the value of 
$T_0$.

\begin{Propo}[{\bf The number of nodes labelled $j$}]
\label{propo-snake-gfB}
For any $j\in \zs$, the \gf \ $S_j\equiv S_j(t,u)$ that counts
binary trees by the number of nodes and the number of nodes labelled $j$ is
algebraic of degree at most $4$ over $\qs(T,u)$ (and thus has degree
at most $8$ over $\qs(t,u)$). More precisely,
$$
\frac{T^2(u-1)^2}{u(T-S_0)^2}
=
\frac{(T-1)^4S_0^2-2TS_0 (T-1)^2(3-9T+7T^2)+T^2(T^2+T-1)^2}
{(T-1)(S_0-1)(T^2+TS_0-S_0)^2}
%\eeq
$$
and all the series $S_j$ belong to $\qs(t,u,S_0)$.
Moreover,  for all $j \ge 0$,
%\beq\label{product-formSjB}
$$
S_j= T\,\frac{(1+\mu Z^j)(1+\mu Z^{j+5})}{(1+\mu Z^{j+2})(1+\mu Z^{j+3})},
%\eeq
$$
where $Z\equiv Z(t)$ is given by~\eqref{eqZtB}
and $\mu\equiv \mu(t,u)$ is the unique \fps\ in $t$
satisfying
%\beq\label{eqMZB}
$$
\mu= (u-1) \frac{Z(1+\mu Z)^2(1+\mu Z^2)(1+\mu Z^6)}
{(1+Z)^2(1+Z+Z^2)(1-Z)^3(1-\mu^2 Z^5)}.
%\eeq
$$
The series $\mu(t,u)$  has polynomial coefficients in $u$,
and satisfies $\mu(t,1)=0$. It has degree $4$ over $\qs(Z,u)$ and $16$
over $\qs(t,u)$. 
\end{Propo}
\noindent{\bf Comment on the proof.} The proof is similar to the proof of
Proposition~\ref{propo-snake-gf-rep} until one computes the equation
satisfied by $S_0$. But then, the relation $S_0=1+tS_1^2$ does not
allow us to conclude that $S_1$ belongs to $\qs(t,u,S_0)$. Instead, we
compute the algebraic equation satisfied by $S_1$. It is found to have
degree $4$ over $\qs(u,T)$. The above relation between $S_0$ and $S_1$
shows that $S_0$ belongs to the extension of $\qs(t,u)$ generated by
$S_1$. Comparing the degrees implies finally that
$\qs(t,u,S_0)=\qs(t,u,S_1)$. Then~\eqref{rec-Sj-binary} shows, by
induction on $j$, that all the series $S_j$ belong to this field.

\begin{Propo}[{\bf  The number of nodes labelled $j$ or more}]
\label{propo-snake-gf-repB}
Let $j \in \zs$. The \gf\ $R_j(t,u)\equiv R_j$
that counts binary trees by the number of nodes and the number of
nodes labelled $j$ or more
is  algebraic over $\qs(t,u)$.
% , of degree at most $16$.  It has degree at most $4$ over
% $\qs(T,\tilde T)$, where $T\equiv T(t)$ and $\tilde T\equiv T(tu)$. 
More precisely, $R_0$ has degree $16$ over
$\qs(t,u)$ and degree $4$ over  $\qs(T,\tilde T)$, with $\tilde T=
T(tu)$, and all the series
$R_j$ belong to  $\qs(T,\tilde T, R_0)=\qs(t,u,R_0)$.
Moreover,  for all $j \ge 0$,
%\beq\label{product-formRjB}
$$
R_j= T\,\frac{(1+\nu Z^j)(1+\nu Z^{j+5})}{(1+\nu Z^{j+2})(1+\nu Z^{j+3})},
%\eeq
$$
where $Z\equiv Z(t)$ is given by~\eqref{eqZtP}
and $\nu\equiv \nu(t,u)$ is a \fps\ in $t$, with polynomial
coefficients in $u$, which is algebraic
%  of degree $2$ over $\qs(Z(t),Z(tu))$, 
of degree $8$ over $\qs(u,Z)$ and  $32$ over $\qs(t,u)$. 
This series satisfies $\nu(t,1)=0$. 
The first terms in its expansion are:
$$
\nu(t,u)=(u-1)\Big(
t+(u+1)t^2+(2u^2+3u+3)t^3+O(t^4)\Big) .
$$
\end{Propo}
\noindent{\bf Comment on the proof.} The proof is similar to the proof of
Proposition~\ref{propo-snake-gf-rep} until one computes the equation
satisfied by $R_0$. One finds that $R_0$ has degree $4$ over
$\qs(T,\tilde T)$, and degree 16 over
$\qs(t,u)$. Using~\eqref{eq-symmetryB}, one then derives an equation
satisfied by $R_1$. Strangely enough, it turns out that the minimal
polynomials of  $R_0$ and $R_1$ over $\qs(T,\tilde T)$ (or over
$\qs(t,u)$) \emm are the same,. The two series are 
of course different:
$$
\begin{array}{lll}
R_0(t,u)&=&1+tu+u(1+u)t^2+u(2u^2+2u+1)t^3+u(1+u)(4u^2+u+2)t^4+O(t^5),\\
R_1(t,u)&=&1+t+(1+u)t^2+(u^2+2u+2)t^3+(1+u)(2u^2+u+4)t^4+O(t^5).
\end{array}
$$

Let  $P(x)$ be the minimal polynomial of $R_0$ and $R_1$ over
$\GK\equiv\qs(T,\tilde T)$.
We want to prove that $R_1$ belongs to the extension of
$\GK$ generated by $R_0$. Note that this property does not simply
follow from the fact that $R_0$ and $R_1$ are conjugate roots of
$P$. For instance, for a \emm generic, polynomial $P$ of degree 4 over $\GK$,
with Galois group $\Sn_4$, the four extensions of $\GK$ generated by the
roots of $P$ are  different. We are going to determine the Galois
group of our polynomial $P$, using the general strategy described
in~\cite[p.~141--142]{bastida}. The 
resolvent cubic  of $P$, which we denote $R$ below, is 
found to factor into a linear term and a quadratic one. Hence the
Galois group of $R$ over $\GK$ has order 2. This implies
 that the Galois group $G$ of $P$ over $\GK$ is either the cyclic
 group of order 4 or 
 the dihedral group of order 8. In the former case, the four
 extensions of $\GK$ generated by the roots of $P$ coincide (and are
 equal to the splitting field of $P$) and we are done. In the latter case,
 there exists a labelling of the four roots of $P$, say $X_0,X_1,
X_2, X_3$, such that the group $G$, seen as a subgroup of the
permutations of $\{0,1,2,3\}$, is
$$G=\{ \hbox{id}, (0,1), (2,3), (0,1) (2,3), (0,2) (1,3), (0,3) (1,2),
(0, 2,1,3), (0,3,1,2)\}.$$
Then

$\bullet$ the simple extensions of $\GK$ generated by the $X_i$ satisfy
    $\GK(X_0)=\GK(X_1)$ and     $\GK(X_2)=\GK(X_3)$,

$\bullet$ the root of the resolvent that belongs to $\GK$
is $Y=X_0X_1+X_2X_3$.

The root of $R$ that belongs to $\GK=\qs(T,\tilde T)$ is found to be 
$Y=t^{-3} /(u(1+u)) + O(t^{-2})$. We already know two roots of $P$,
namely $R_0$ and $R_1$, which are equal to $1+O(t)$. The other two
roots are respectively of the form $Q_2=-t^{-2}/u+ O(t^{-1})$ and
$Q_3=-t^{-1}/(1+u)+ O(1)$.  
From the value of $Y$, we conclude that the above properties hold with
$X_0=R_0$, 
$X_1=R_1$, $X_2=Q_2$ and  $X_3=Q_3$. In particular, $R_0$ and $R_1$
belong to the same extension of degree 4 of $\qs(T,\tilde T)$. 

Then, an induction on $j$, based on~\eqref{rec-RjB}, implies that for
all $j\ge 2$, the series $R_j$ belongs to the extension of
$\GK=\qs(T,\tilde T)$ generated by $R_0$. Since $R_1(t,u)=R_0(tu,1/u)$ and
$\GK(R_0)=\GK(R_1)$, the field $\GK(R_0)$ is invariant under the
transformation $A(t,u) \mapsto  A(tu,1/u)$. This property,
combined with~\eqref{eq-symmetryB}, implies that for $j \ge 1$, the
series $R_{-j}$ belongs to $\GK(R_0)$.

\subsubsection{Limit laws}
We now endow the set of binary trees having $n$ nodes with the
uniform distribution, and consider the same random variables as for
above: $M_n$,  the largest label,  $X_n(j)$, the
number of nodes having label $j$, 
and  $X_n^+(j)$, the
number of nodes having label $j$ at least.

Again, we can prove that 
 $M_nn^{-1/4}$ converges in law to  $N_{{\hbox{\sc \small
   ise}}}$,  where $N_{{\hbox{\sc \small
   ise}}}$ is the supremum of the support of the ISE, and that
 for all $\lambda\in \rs$, the sequence $X_n^+( \lambda n^{1/4})/ n$
 converges in law to  $g_{{\hbox{\sc \small ise}}}(\lambda 
)$ where $g_{{\hbox{\sc \small ise}}}$ is the tail
distribution function of the ISE. The arguments are the same as for
 our first class of trees (Sections~\ref{section-supremum}
 and~\ref{section-coro-tail}). 
Hence we could  just as well have started from the enumerative results of
Section~\ref{section-enumerative-binary}, rather than from those of
Section~\ref{section-snake-gf}, to characterize the laws of
$N_{{\hbox{\sc \small  ise}}} $ and $g_{{\hbox{\sc \small
      ise}}}(\lambda)$ (Propositions~\ref{coro-support}
and~\ref{coro-tail}). 

More remarkably, we have performed on
$X_n(j)$ an analysis similar to that of
Section~\ref{section-limit-local}, and obtained \emm the same, local
limit law. In other words, for all $\lambda \ge 0$,
the sequence  $X_n(\lfloor \lambda n^{1/4}\rfloor) n^{-3/4}$ converges
in law to  $
 f_{{\hbox{\sc \small ise}}}(\lambda )
$
where $f_{{\hbox{\sc \small ise}}}$ is the conjectured
density of the ISE, given in Conjecture~\ref{conj-density}.

In all three cases, the convergence of the moments holds as well.

\bigskip

%%%%%%%%%%%%%%%%%%%%%%%%%%%%%%%%%%%%%%%%%%%%%%%%%%

\noindent{\bf Acknowledgements.} I am extremely grateful to
Jean-Fran\c cois Marckert, who not only suggested the topic of this
paper, but also spent a lot of time explaining me the probabilistic
implications of the limit results I had obtained. My thanks also go to
Alain Rouault, who identified the moments occurring in
Proposition~\ref{limit-0}, and to Philippe Flajolet, Jean-François Le
Gall and Guy Louchard 
for providing useful  references.

\bibliographystyle{plain}
\bibliography{ise}

\bigskip
\noindent{\bf Appendix: convergence of random distribution functions}\\
We want to prove the result stated without proof at the beginning of
the proof of Section~\ref{section-coro-tail}. 

Recall that a sequence of real random variables $Z_n$ converges in law to
another random variable $Z$ if and only if for all $x \in \rs$ such
that $\PP(Z=x)=0$,
$$
\lim_n \PP(Z_n\le x) =\PP(Z\le x).
$$
This implies   the so-called \emm
portmanteau, inequality: for \emm all,\/ $x \in \rs$,
\beq
\label{conv-law}
\PP(Z<x)\le \liminf \PP(Z_n\le x)\le \limsup \PP(Z_n\le x) \le\PP(Z\le x).
\eeq

Let us now use the notation of Section~\ref{section-coro-tail}. The
 convergence of $\mu_n$ to $\mu$ implies  that for any 
 bounded  Lipschitz function $f$ on $\rs$~\cite[p.~71--74]{VV}:
$$
\int_\rs f(x) d\mu_n(x) \dconv \int_\rs f(x) d\mu(x).
$$
Let $\lambda\in \rs$ and let  $f_{\lambda,\epsilon}$ and $g_{\lambda,\epsilon}$
be the functions plotted in Figure~\ref{fig-epsilon}. Then
$$
\int_\rs g_{\lambda,\epsilon}(y) d\mu_n(y) \le 
F_n(\lambda) = \mu_n(-\infty, \lambda] 
\le \int_\rs f_{\lambda,\epsilon}(y) d\mu_n(y).
$$
 Hence, for all $x\in \rs$,
$$
\PP\left(\int_\rs f_{\lambda,\epsilon}(y) d\mu_n(y) \le x\right)
\le \PP(F_n(\lambda)\le x)\le \PP\left(\int_\rs
g_{\lambda,\epsilon}(y) d\mu_n(y) \le x\right) .
$$
Since $\mu_n$ converges to $\mu$, and $g_{\lambda,\epsilon}$ is a
bounded Lipschitz function,
$$%Z_n=
 \int_\rs g_{\lambda,\epsilon}(y) d\mu_n(y) \dconv
% converges in law to 
 \int_\rs g_{\lambda,\epsilon}(y) d\mu(y) .
$$
A similar result holds for the integral involving $f_{\lambda,\epsilon}$.
Thus~\eqref{conv-law} implies
$$
\PP\left(\int_\rs f_{\lambda,\epsilon}(y) d\mu(y) < x\right) \le 
\liminf \PP(F_n(\lambda)\le x)\le \limsup \PP(F_n(\lambda)\le x)\le
 \PP\left(\int_\rs g_{\lambda,\epsilon}(y) d\mu(y) \le x\right) .
$$
The integral occurring in the
rightmost expression of this inequality is bounded from below by $\mu(-\infty,
\lambda-\epsilon]$, while the integral involving
$f_{\lambda,\epsilon}$ is bounded from above by $\mu(-\infty,
\lambda+\epsilon]$.  Hence
$$
\PP(\mu(-\infty, \lambda+\epsilon]<x)  \le 
\liminf \PP(F_n(\lambda)\le x)\le \limsup \PP(F_n(\lambda)\le x)\le
\PP(\mu(-\infty, \lambda-\epsilon]\le x).
$$
Taking the limit $\epsilon \rightarrow 0^+$ gives:
$$
\PP(\mu(-\infty, \lambda]<x)  \le 
\liminf \PP(F_n(\lambda)\le x)\le \limsup \PP(F_n(\lambda)\le x)\le
\PP(\mu(-\infty, \lambda)\le x).
$$
If, in addition, the measure $\mu$ does not assign a positive
probability to $\lambda$, the
rightmost expression in the above inequality equals 
$\PP(\mu(-\infty, \lambda] \le x)$. The inequality becomes
$$
\PP(F(\lambda)< x) \le \liminf \PP(F_n(\lambda) \le x)
\le \limsup \PP(F_n(\lambda) \le x)\le \PP(F(\lambda)\le x),
$$
so that for all $x$ such that $\PP(F(\lambda)= x)=0$,
$$
\lim\PP(F_n(\lambda) \le x) = \PP(F(\lambda)\le x).
$$
That is,  $F_n(\lambda)$ converges in law to $F(\lambda)$.

%\tableofcontents
\end{document}

%% file: tree-examples.pstex_t
\begin{picture}(0,0)%
\includegraphics{tree-examples.pstex}%
\end{picture}%
\setlength{\unitlength}{2368sp}%
\begingroup\makeatletter\ifx\SetFigFont\undefined%
\gdef\SetFigFont#1#2#3#4#5{%
  \reset@font\fontsize{#1}{#2pt}%
  \fontfamily{#3}\fontseries{#4}\fontshape{#5}%
  \selectfont}%
\fi\endgroup%
\begin{picture}(10575,2620)(526,-1919)
\put(1351,-1861){\makebox(0,0)[lb]{\smash{\SetFigFont{9}{10.8}{\familydefault}{\mddefault}{\updefault}{\color[rgb]{0,0,0}$0$}%
}}}
\put(2251,-1036){\makebox(0,0)[lb]{\smash{\SetFigFont{9}{10.8}{\familydefault}{\mddefault}{\updefault}{\color[rgb]{0,0,0}$-1$}%
}}}
\put(4426,-1861){\makebox(0,0)[lb]{\smash{\SetFigFont{9}{10.8}{\familydefault}{\mddefault}{\updefault}{\color[rgb]{0,0,0}$0$}%
}}}
\put(4651,-1036){\makebox(0,0)[lb]{\smash{\SetFigFont{9}{10.8}{\familydefault}{\mddefault}{\updefault}{\color[rgb]{0,0,0}$1$}%
}}}
\put(526,-1036){\makebox(0,0)[lb]{\smash{\SetFigFont{9}{10.8}{\familydefault}{\mddefault}{\updefault}{\color[rgb]{0,0,0}$1$}%
}}}
\put(3601,-1036){\makebox(0,0)[lb]{\smash{\SetFigFont{9}{10.8}{\familydefault}{\mddefault}{\updefault}{\color[rgb]{0,0,0}$1$}%
}}}
\put(1576,-1036){\makebox(0,0)[lb]{\smash{\SetFigFont{9}{10.8}{\familydefault}{\mddefault}{\updefault}{\color[rgb]{0,0,0}$1$}%
}}}
\put(5326,-1036){\makebox(0,0)[lb]{\smash{\SetFigFont{9}{10.8}{\familydefault}{\mddefault}{\updefault}{\color[rgb]{0,0,0}$0$}%
}}}
\put(1351,-286){\makebox(0,0)[lb]{\smash{\SetFigFont{9}{10.8}{\familydefault}{\mddefault}{\updefault}{\color[rgb]{0,0,0}$2$}%
}}}
\put(1801,-286){\makebox(0,0)[lb]{\smash{\SetFigFont{9}{10.8}{\familydefault}{\mddefault}{\updefault}{\color[rgb]{0,0,0}$0$}%
}}}
\put(2176,-286){\makebox(0,0)[lb]{\smash{\SetFigFont{9}{10.8}{\familydefault}{\mddefault}{\updefault}{\color[rgb]{0,0,0}$-2$}%
}}}
\put(4426,-286){\makebox(0,0)[lb]{\smash{\SetFigFont{9}{10.8}{\familydefault}{\mddefault}{\updefault}{\color[rgb]{0,0,0}$1$}%
}}}
\put(4876,-286){\makebox(0,0)[lb]{\smash{\SetFigFont{9}{10.8}{\familydefault}{\mddefault}{\updefault}{\color[rgb]{0,0,0}$0$}%
}}}
\put(5251,-286){\makebox(0,0)[lb]{\smash{\SetFigFont{9}{10.8}{\familydefault}{\mddefault}{\updefault}{\color[rgb]{0,0,0}$-1$}%
}}}
\put(9901,-1711){\makebox(0,0)[lb]{\smash{\SetFigFont{9}{10.8}{\familydefault}{\mddefault}{\updefault}{\color[rgb]{0,0,0}$0$}%
}}}
\put(9226,-1036){\makebox(0,0)[lb]{\smash{\SetFigFont{9}{10.8}{\familydefault}{\mddefault}{\updefault}{\color[rgb]{0,0,0}${ -1}$}%
}}}
\put(10651,-1186){\makebox(0,0)[lb]{\smash{\SetFigFont{9}{10.8}{\familydefault}{\mddefault}{\updefault}{\color[rgb]{0,0,0}$1$}%
}}}
\put(11101,-436){\makebox(0,0)[lb]{\smash{\SetFigFont{9}{10.8}{\familydefault}{\mddefault}{\updefault}{\color[rgb]{0,0,0}$ 2$}%
}}}
\put(9826,-361){\makebox(0,0)[lb]{\smash{\SetFigFont{9}{10.8}{\familydefault}{\mddefault}{\updefault}{\color[rgb]{0,0,0}$ 0$}%
}}}
\put(10801,164){\makebox(0,0)[lb]{\smash{\SetFigFont{9}{10.8}{\familydefault}{\mddefault}{\updefault}{\color[rgb]{0,0,0}$ 1$}%
}}}
\put(10201,164){\makebox(0,0)[lb]{\smash{\SetFigFont{9}{10.8}{\familydefault}{\mddefault}{\updefault}{\color[rgb]{0,0,0}$1$}%
}}}
\end{picture}

%% file: profils.pstex_t
\begin{picture}(0,0)%
\includegraphics{profils.pstex}%
\end{picture}%
\setlength{\unitlength}{1973sp}%
\begingroup\makeatletter\ifx\SetFigFont\undefined%
\gdef\SetFigFont#1#2#3#4#5{%
  \reset@font\fontsize{#1}{#2pt}%
  \fontfamily{#3}\fontseries{#4}\fontshape{#5}%
  \selectfont}%
\fi\endgroup%
\begin{picture}(10362,3514)(451,-3119)
\put(451,-2761){\makebox(0,0)[lb]{\smash{{\SetFigFont{8}{9.6}{\familydefault}{\mddefault}{\updefault}{\color[rgb]{0,0,0}$0$}%
}}}}
\put(451,-2161){\makebox(0,0)[lb]{\smash{{\SetFigFont{8}{9.6}{\familydefault}{\mddefault}{\updefault}{\color[rgb]{0,0,0}$1$}%
}}}}
\put(451,-1561){\makebox(0,0)[lb]{\smash{{\SetFigFont{8}{9.6}{\familydefault}{\mddefault}{\updefault}{\color[rgb]{0,0,0}$2$}%
}}}}
\put(451,-961){\makebox(0,0)[lb]{\smash{{\SetFigFont{8}{9.6}{\familydefault}{\mddefault}{\updefault}{\color[rgb]{0,0,0}$3$}%
}}}}
\put(451,-361){\makebox(0,0)[lb]{\smash{{\SetFigFont{8}{9.6}{\familydefault}{\mddefault}{\updefault}{\color[rgb]{0,0,0}$4$}%
}}}}
\put(5401,-361){\makebox(0,0)[lb]{\smash{{\SetFigFont{8}{9.6}{\familydefault}{\mddefault}{\updefault}{\color[rgb]{0,0,0}$2$}%
}}}}
\put(5401,-961){\makebox(0,0)[lb]{\smash{{\SetFigFont{8}{9.6}{\familydefault}{\mddefault}{\updefault}{\color[rgb]{0,0,0}$3$}%
}}}}
\put(5401,-1561){\makebox(0,0)[lb]{\smash{{\SetFigFont{8}{9.6}{\familydefault}{\mddefault}{\updefault}{\color[rgb]{0,0,0}$4$}%
}}}}
\put(5401,-2161){\makebox(0,0)[lb]{\smash{{\SetFigFont{8}{9.6}{\familydefault}{\mddefault}{\updefault}{\color[rgb]{0,0,0}$2$}%
}}}}
\put(5401,-2761){\makebox(0,0)[lb]{\smash{{\SetFigFont{8}{9.6}{\familydefault}{\mddefault}{\updefault}{\color[rgb]{0,0,0}$1$}%
}}}}
\put(7276,239){\makebox(0,0)[lb]{\smash{{\SetFigFont{8}{9.6}{\familydefault}{\mddefault}{\updefault}{\color[rgb]{0,0,0}$2$}%
}}}}
\put(7801,239){\makebox(0,0)[lb]{\smash{{\SetFigFont{8}{9.6}{\familydefault}{\mddefault}{\updefault}{\color[rgb]{0,0,0}$2$}%
}}}}
\put(8476,239){\makebox(0,0)[lb]{\smash{{\SetFigFont{8}{9.6}{\familydefault}{\mddefault}{\updefault}{\color[rgb]{0,0,0}$4$}%
}}}}
\put(9076,239){\makebox(0,0)[lb]{\smash{{\SetFigFont{8}{9.6}{\familydefault}{\mddefault}{\updefault}{\color[rgb]{0,0,0}$2$}%
}}}}
\put(9676,239){\makebox(0,0)[lb]{\smash{{\SetFigFont{8}{9.6}{\familydefault}{\mddefault}{\updefault}{\color[rgb]{0,0,0}$1$}%
}}}}
\put(10276,239){\makebox(0,0)[lb]{\smash{{\SetFigFont{8}{9.6}{\familydefault}{\mddefault}{\updefault}{\color[rgb]{0,0,0}$1$}%
}}}}
\put(10276,-3061){\makebox(0,0)[lb]{\smash{{\SetFigFont{8}{9.6}{\familydefault}{\mddefault}{\updefault}{\color[rgb]{0,0,0}$3$}%
}}}}
\put(9676,-3061){\makebox(0,0)[lb]{\smash{{\SetFigFont{8}{9.6}{\familydefault}{\mddefault}{\updefault}{\color[rgb]{0,0,0}$2$}%
}}}}
\put(9076,-3061){\makebox(0,0)[lb]{\smash{{\SetFigFont{8}{9.6}{\familydefault}{\mddefault}{\updefault}{\color[rgb]{0,0,0}$1$}%
}}}}
\put(8476,-3061){\makebox(0,0)[lb]{\smash{{\SetFigFont{8}{9.6}{\familydefault}{\mddefault}{\updefault}{\color[rgb]{0,0,0}$0$}%
}}}}
\put(7876,-3061){\makebox(0,0)[lb]{\smash{{\SetFigFont{8}{9.6}{\familydefault}{\mddefault}{\updefault}{\color[rgb]{0,0,0}$-1$}%
}}}}
\put(7276,-3061){\makebox(0,0)[lb]{\smash{{\SetFigFont{8}{9.6}{\familydefault}{\mddefault}{\updefault}{\color[rgb]{0,0,0}$-2$}%
}}}}
\end{picture}%

%% file: tree-decomp.pstex_t
\begin{picture}(0,0)%
\includegraphics{tree-decomp.pstex}%
\end{picture}%
\setlength{\unitlength}{3947sp}%
\begingroup\makeatletter\ifx\SetFigFont\undefined%
\gdef\SetFigFont#1#2#3#4#5{%
  \reset@font\fontsize{#1}{#2pt}%
  \fontfamily{#3}\fontseries{#4}\fontshape{#5}%
  \selectfont}%
\fi\endgroup%
\begin{picture}(4074,1570)(-761,-1319)
\put(1951,-136){\makebox(0,0)[lb]{\smash{\SetFigFont{12}{14.4}{\familydefault}{\mddefault}{\updefault}{\color[rgb]{0,0,0}$T_{j\pm 1}$}%
}}}
\put(2626,-1261){\makebox(0,0)[lb]{\smash{\SetFigFont{12}{14.4}{\familydefault}{\mddefault}{\updefault}{\color[rgb]{0,0,0}$j$}%
}}}
\put(1726,-811){\makebox(0,0)[lb]{\smash{\SetFigFont{12}{14.4}{\familydefault}{\mddefault}{\updefault}{\color[rgb]{0,0,0}$j\pm 1$}%
}}}
\put(2701,-211){\makebox(0,0)[lb]{\smash{\SetFigFont{12}{14.4}{\familydefault}{\mddefault}{\updefault}{\color[rgb]{0,0,0}$T_j$}%
}}}
\put(-224,-1261){\makebox(0,0)[lb]{\smash{\SetFigFont{12}{14.4}{\familydefault}{\mddefault}{\updefault}{\color[rgb]{0,0,0}$j$}%
}}}
\put(226,-661){\makebox(0,0)[lb]{\smash{\SetFigFont{12}{14.4}{\familydefault}{\mddefault}{\updefault}{\color[rgb]{0,0,0}$=$}%
}}}
\put(1276,-661){\makebox(0,0)[lb]{\smash{\SetFigFont{12}{14.4}{\familydefault}{\mddefault}{\updefault}{\color[rgb]{0,0,0}$+$}%
}}}
\put(-449,-211){\makebox(0,0)[lb]{\smash{\SetFigFont{12}{14.4}{\familydefault}{\mddefault}{\updefault}{\color[rgb]{0,0,0}$T_j$}%
}}}
\put(826,-886){\makebox(0,0)[lb]{\smash{\SetFigFont{12}{14.4}{\familydefault}{\mddefault}{\updefault}{\color[rgb]{0,0,0}$j$}%
}}}
\end{picture}

%% file: contour.pstex_t
\begin{picture}(0,0)%
\includegraphics{contour.pstex}%
\end{picture}%
\setlength{\unitlength}{3947sp}%
\begingroup\makeatletter\ifx\SetFigFont\undefined%
\gdef\SetFigFont#1#2#3#4#5{%
  \reset@font\fontsize{#1}{#2pt}%
  \fontfamily{#3}\fontseries{#4}\fontshape{#5}%
  \selectfont}%
\fi\endgroup%
\begin{picture}(3024,3024)(589,-2548)
\put(2476,-436){\makebox(0,0)[lb]{\smash{\SetFigFont{12}{14.4}{\familydefault}{\mddefault}{\updefault}{\color[rgb]{0,0,0}$r_n/8$}%
}}}
\put(2799,-1165){\makebox(0,0)[lb]{\smash{\SetFigFont{10}{12.0}{\familydefault}{\mddefault}{\updefault}{\color[rgb]{0,0,0}$1/8$}%
}}}
\put(976,-2161){\makebox(0,0)[lb]{\smash{\SetFigFont{12}{14.4}{\familydefault}{\mddefault}{\updefault}{\color[rgb]{0,0,0}${\C}_n^{(1)}$}%
}}}
\put(2701,-1486){\makebox(0,0)[lb]{\smash{\SetFigFont{12}{14.4}{\familydefault}{\mddefault}{\updefault}{\color[rgb]{0,0,0}${\C}_n^{(2)}$}%
}}}
\end{picture}

%% file: hankel.pstex_t
\begin{picture}(0,0)%
\includegraphics{hankel.pstex}%
\end{picture}%
\setlength{\unitlength}{3947sp}%
\begingroup\makeatletter\ifx\SetFigFont\undefined%
\gdef\SetFigFont#1#2#3#4#5{%
  \reset@font\fontsize{#1}{#2pt}%
  \fontfamily{#3}\fontseries{#4}\fontshape{#5}%
  \selectfont}%
\fi\endgroup%
\begin{picture}(6238,1583)(599,-1076)
\put(4918,-639){\makebox(0,0)[lb]{\smash{{\SetFigFont{10}{12.0}{\familydefault}{\mddefault}{\updefault}{\color[rgb]{0,0,0}$0$}%
}}}}
\put(4635, 69){\makebox(0,0)[lb]{\smash{{\SetFigFont{10}{12.0}{\familydefault}{\mddefault}{\updefault}{\color[rgb]{0,0,0}$i$}%
}}}}
\put(4115,-403){\makebox(0,0)[lb]{\smash{{\SetFigFont{10}{12.0}{\familydefault}{\mddefault}{\updefault}{\color[rgb]{0,0,0}$-1$}%
}}}}
\put(6241,-639){\makebox(0,0)[lb]{\smash{{\SetFigFont{12}{14.4}{\familydefault}{\mddefault}{\updefault}{\color[rgb]{0,0,0}$x_n\sim \log ^2 n$}%
}}}}
\put(5551,164){\makebox(0,0)[lb]{\smash{{\SetFigFont{12}{14.4}{\familydefault}{\mddefault}{\updefault}{\color[rgb]{0,0,0}$\Ha_n$}%
}}}}
\put(6376,-61){\makebox(0,0)[lb]{\smash{{\SetFigFont{12}{14.4}{\familydefault}{\mddefault}{\updefault}{\color[rgb]{0,0,0}$z_n$}%
}}}}
\put(1603,-639){\makebox(0,0)[lb]{\smash{{\SetFigFont{10}{12.0}{\familydefault}{\mddefault}{\updefault}{\color[rgb]{0,0,0}$0$}%
}}}}
\put(1320, 69){\makebox(0,0)[lb]{\smash{{\SetFigFont{10}{12.0}{\familydefault}{\mddefault}{\updefault}{\color[rgb]{0,0,0}$i$}%
}}}}
\put(800,-403){\makebox(0,0)[lb]{\smash{{\SetFigFont{10}{12.0}{\familydefault}{\mddefault}{\updefault}{\color[rgb]{0,0,0}$-1$}%
}}}}
\put(2251,164){\makebox(0,0)[lb]{\smash{{\SetFigFont{12}{14.4}{\familydefault}{\mddefault}{\updefault}{\color[rgb]{0,0,0}$\Ha$}%
}}}}
\end{picture}%

%% file: fepsilon.pstex_t
\begin{picture}(0,0)%
\includegraphics{fepsilon.pstex}%
\end{picture}%
\setlength{\unitlength}{3947sp}%
\begingroup\makeatletter\ifx\SetFigFont\undefined%
\gdef\SetFigFont#1#2#3#4#5{%
  \reset@font\fontsize{#1}{#2pt}%
  \fontfamily{#3}\fontseries{#4}\fontshape{#5}%
  \selectfont}%
\fi\endgroup%
\begin{picture}(7296,1719)(301,-1123)
\put(1726,-1036){\makebox(0,0)[lb]{\smash{{\SetFigFont{10}{12.0}{\familydefault}{\mddefault}{\updefault}{\color[rgb]{0,0,0}$\lambda+\epsilon$}%
}}}}
\put(1426,-1036){\makebox(0,0)[lb]{\smash{{\SetFigFont{10}{12.0}{\familydefault}{\mddefault}{\updefault}{\color[rgb]{0,0,0}$\lambda$}%
}}}}
\put(676,164){\makebox(0,0)[lb]{\smash{{\SetFigFont{10}{12.0}{\familydefault}{\mddefault}{\updefault}{\color[rgb]{0,0,0}$1$}%
}}}}
\put(2626,-811){\makebox(0,0)[lb]{\smash{{\SetFigFont{10}{12.0}{\familydefault}{\mddefault}{\updefault}{\color[rgb]{0,0,0}$x$}%
}}}}
\put(301,464){\makebox(0,0)[lb]{\smash{{\SetFigFont{10}{12.0}{\familydefault}{\mddefault}{\updefault}{\color[rgb]{0,0,0}$f_{\lambda,\epsilon}(x)$}%
}}}}
\put(4051,-1036){\makebox(0,0)[lb]{\smash{{\SetFigFont{10}{12.0}{\familydefault}{\mddefault}{\updefault}{\color[rgb]{0,0,0}$\lambda$}%
}}}}
\put(3301,164){\makebox(0,0)[lb]{\smash{{\SetFigFont{10}{12.0}{\familydefault}{\mddefault}{\updefault}{\color[rgb]{0,0,0}$1$}%
}}}}
\put(2926,464){\makebox(0,0)[lb]{\smash{{\SetFigFont{10}{12.0}{\familydefault}{\mddefault}{\updefault}{\color[rgb]{0,0,0}$g_{\lambda,\epsilon}(x)$}%
}}}}
\put(3601,-1036){\makebox(0,0)[lb]{\smash{{\SetFigFont{10}{12.0}{\familydefault}{\mddefault}{\updefault}{\color[rgb]{0,0,0}$\lambda-\epsilon$}%
}}}}
\put(6451,-1036){\makebox(0,0)[lb]{\smash{{\SetFigFont{10}{12.0}{\familydefault}{\mddefault}{\updefault}{\color[rgb]{0,0,0}$\lambda+\epsilon$}%
}}}}
\put(6151,-1036){\makebox(0,0)[lb]{\smash{{\SetFigFont{10}{12.0}{\familydefault}{\mddefault}{\updefault}{\color[rgb]{0,0,0}$\lambda$}%
}}}}
\put(7351,-811){\makebox(0,0)[lb]{\smash{{\SetFigFont{10}{12.0}{\familydefault}{\mddefault}{\updefault}{\color[rgb]{0,0,0}$x$}%
}}}}
\put(5401,164){\makebox(0,0)[lb]{\smash{{\SetFigFont{10}{12.0}{\familydefault}{\mddefault}{\updefault}{\color[rgb]{0,0,0}$1$}%
}}}}
\put(5026,464){\makebox(0,0)[lb]{\smash{{\SetFigFont{10}{12.0}{\familydefault}{\mddefault}{\updefault}{\color[rgb]{0,0,0}$h_{\lambda,\epsilon}(x)$}%
}}}}
\put(5701,-1036){\makebox(0,0)[lb]{\smash{{\SetFigFont{10}{12.0}{\familydefault}{\mddefault}{\updefault}{\color[rgb]{0,0,0}$\lambda-\epsilon$}%
}}}}
\end{picture}%

%% file: tree-decomp-binary.pstex_t
\begin{picture}(0,0)%
\includegraphics{tree-decomp-binary.pstex}%
\end{picture}%
\setlength{\unitlength}{3947sp}%
\begingroup\makeatletter\ifx\SetFigFont\undefined%
\gdef\SetFigFont#1#2#3#4#5{%
  \reset@font\fontsize{#1}{#2pt}%
  \fontfamily{#3}\fontseries{#4}\fontshape{#5}%
  \selectfont}%
\fi\endgroup%
\begin{picture}(3924,1420)(-461,-1319)
\put(2626,-1261){\makebox(0,0)[lb]{\smash{\SetFigFont{12}{14.4}{\familydefault}{\mddefault}{\updefault}{\color[rgb]{0,0,0}$j$}%
}}}
\put(1876,-136){\makebox(0,0)[lb]{\smash{\SetFigFont{12}{14.4}{\familydefault}{\mddefault}{\updefault}{\color[rgb]{0,0,0}$T_{j-1}$}%
}}}
\put(2851,-136){\makebox(0,0)[lb]{\smash{\SetFigFont{12}{14.4}{\familydefault}{\mddefault}{\updefault}{\color[rgb]{0,0,0}$T_{j+1}$}%
}}}
\put( 76,-1111){\makebox(0,0)[lb]{\smash{\SetFigFont{12}{14.4}{\familydefault}{\mddefault}{\updefault}{\color[rgb]{0,0,0}$j$}%
}}}
\put(376,-661){\makebox(0,0)[lb]{\smash{\SetFigFont{12}{14.4}{\familydefault}{\mddefault}{\updefault}{\color[rgb]{0,0,0}$=$}%
}}}
\put(1576,-661){\makebox(0,0)[lb]{\smash{\SetFigFont{12}{14.4}{\familydefault}{\mddefault}{\updefault}{\color[rgb]{0,0,0}$+$}%
}}}
\put(976,-661){\makebox(0,0)[lb]{\smash{\SetFigFont{12}{14.4}{\rmdefault}{\mddefault}{\updefault}{\color[rgb]{0,0,0}$\emptyset$}%
}}}
\end{picture}